\documentclass[letterpaper, 10pt, dvipsnames]{article}

\usepackage{bold-extra}
\usepackage[utf8]{inputenc}    
\usepackage[T1]{fontenc}
\usepackage[english]{babel}
\usepackage[normalem]{ulem}
\usepackage{dsfont, bm, pifont, mathrsfs}
\usepackage{stmaryrd}
\usepackage{comment}
\usepackage{bbm}

\usepackage[margin=30pt]{subcaption}
\usepackage{float, graphicx, caption, wrapfig, tikz}
\usetikzlibrary{shapes.misc}

\usepackage{amsthm, amsmath, amsfonts, amssymb, mathrsfs, mathtools}
\usepackage{amsrefs}

\usepackage{hyperref, xcolor, titlesec} 
\hypersetup{colorlinks = true, urlcolor = RoyalBlue!70!Black,linkcolor=BrickRed, citecolor=Plum, bookmarksopen = true}
\usepackage[nomarginpar]{geometry}
\numberwithin{equation}{section}

\usepackage{cleveref}
\renewcommand{\ge}{\geqslant}
\renewcommand{\le}{\leqslant}

\DeclareMathOperator{\Id}{Id}

\let \d \relax
\newcommand{\d}{\mathrm{d}}
\newcommand{\I}{\mathrm{i}}
\newcommand{\scp}[2]{\langle #1,#2\rangle}

\let\O\relax
\newcommand{\O}[1]{\mathcal{O}\left(#1\right)}
\newcommand{\E}[1]{\mathds{E}\left[#1\right]}

\newcommand{\unn}[2]{[\![#1,#2]\!]}
\let\Im\relax
\DeclareMathOperator{\Im}{Im}
\let\Re\relax
\DeclareMathOperator{\Re}{Re}
\DeclareMathOperator{\Tr}{Tr}
\DeclareMathOperator{\tr}{Tr}

\newcommand{\e}{\mathrm{e}}
\newcommand{\one}{\mathbbm{1}}
\renewcommand{\epsilon}{\varepsilon}
\renewcommand{\tilde}{\widetilde}
\newcommand{\vreg}{w}

\def\bet{\begin{thm}}
\def\eet{\end{thm}}
\def\bel{\begin{lem}}
\def\eel{\end{lem}}
\def\bas{\begin{ass}}
\def\eas{\end{ass}}
\def\bec{\begin{cor}}
\def\eec{\end{cor}}
\def\bed{\begin{defn}}
\def\eed{\end{defn}}
\def\bep{\begin{prop}}
\def\eep{\end{prop}}
\def\beq{\begin{equation}}
\def\eeq{\end{equation}}
\def\bea{\begin{equation*}}
\def\eea{\end{equation*}}
\def\bex{\begin{ex}}
\def\eex{\end{ex}}
\def\bp{\begin{proof}}
\def\ep{\end{proof}}

\def\1{{\mathbbm 1}}

\def\benr{\begin{enumerate}[label=(\roman*)]}
\def\eenr{\end{enumerate}}
\def\A{\mathcal{A}}
\def\B{\mathcal{B}}

\def\M{\mathcal{M}}
\def\D{\mathcal{D}}
\def\N{\mathbb{N}}

\def\R{\mathbb{R}}
\def\C{\mathbb{C}}
\def\P{\mathbb{P}}
\def\E{\mathbb{E}}
\def\q{\mathbf{q}}
\def\S{\mathbb{S}}
\newcommand{\cE}{\mathcal{E}}
\renewcommand{\u}{\mathbf{u}}

\def\msc{m_{\mathrm{sc} } }

\def\one{{\mathbbm 1}}
\def\eps{\varepsilon}

\def\scrho{\rho_{\mathrm{sc} }}
\newcommand{\bma}{\begin{bmatrix}}
\newcommand{\ema}{\end{bmatrix}}

\def \wt {\widetilde}
\def\phi{\varphi}
\def\matn{\mathrm{Mat}_{N} }
\DeclareMathOperator{\GOE}{GOE}
\DeclareMathOperator{\GUE}{GUE}
\DeclareMathOperator{\Ai}{Ai}

\renewcommand{\hat}{\widehat}
\newcommand{\ximax}{\widetilde{\bm{\xi}}}
\newcommand{\scc}{c_{\mathrm{sc}}}
\newcommand{\fa}{\mathfrak{a}}
\newcommand{\fb}{\mathfrak{b}}
\newtheorem{ccounter}{ccounter}[section]
\newtheorem{thm}[ccounter]{Theorem}
\newtheorem{lem}[ccounter]{Lemma}
\newtheorem{cor}[ccounter]{Corollary}
\newtheorem{defn}[ccounter]{Definition}
\newtheorem{prop}[ccounter]{Proposition}
\newtheorem{ass}[ccounter]{Assumption}
\newtheorem{ex}[ccounter]{Example}

\theoremstyle{definition}
\newtheorem{rmk}[ccounter]{Remark}

\titleformat{\paragraph}[runin]{\itshape\normalsize}{\theparagraph}{}{}
\titleformat{\subparagraph}[runin]{\itshape\normalsize}{\theparagraph}{0em}{}
\titleformat{\section}[block]{\normalfont\filcenter}{\Large\thesection .}{.7em}{\Large\scshape}
\titleformat{\subsection}[runin]{\normalfont}{\large \bf \thesubsection .}{.5em}{\large\bf}[.]
\titleformat{\subsubsection}[runin]{\normalfont}{\bf \thesubsubsection .}{.5em}{\bf}[.]

\begin{document}
\tikzset{every node/.style={draw,circle,fill=black, scale=.5}}

\title{\scshape\bfseries{Optimal delocalization for generalized Wigner matrices}}
\author{L. \textsc{Benigni}\\\vspace{-0.15cm}\footnotesize{\it{University of Chicago}}\\\footnotesize{\it{lbenigni@math.uchicago.edu}}\and P. \textsc{Lopatto\footnote{P.L.\ was partially supported by NSF grant DMS-1606305, NSF grant DMS-1855509, and the NSF Graduate Research Fellowship Program under grant DGE-1144152.}}\\\vspace{-0.15cm}\footnotesize{\it{Harvard University}}\\\footnotesize{\it{lopatto@math.harvard.edu}}}
\date{}
\maketitle

\renewcommand{\abstractname}{\normalsize\normalfont\scshape Abstract}
\renewenvironment{abstract}
 {\small
  \begin{center}
  \bfseries \abstractname\vspace{-.5em}\vspace{0pt}
  \end{center}
  \list{}{%
    \setlength{\leftmargin}{16mm}
    \setlength{\rightmargin}{\leftmargin}%
  }%
  \item\relax}
 {\endlist}

\begin{abstract}
\small{
We study the eigenvectors of generalized Wigner matrices with subexponential entries and prove that they delocalize at the optimal rate with overwhelming probability. We also prove high probability delocalization bounds with sharp constants. Our proof uses an analysis of the eigenvector moment flow introduced by \hyperlink{cite.bourgade2013eigenvector}{Bourgade and Yau (2017)} to bound logarithmic moments of eigenvector entries for random matrices with small Gaussian components. We then extend this control to all generalized Wigner matrices by comparison arguments based on a framework of regularized eigenvectors, level repulsion, and the observable employed by \hyperlink{cite.landon2018comparison}{Landon, Lopatto, and Marcinek (2018)} to compare extremal eigenvalue statistics. Additionally, we prove level repulsion and eigenvalue overcrowding estimates for the entire spectrum, which may be of independent interest.}
\end{abstract}

{
  \hypersetup{linkcolor=black}
  \tableofcontents
}

\section{Introduction}
Disordered quantum systems often exhibit one of two paradigmatic behaviors, existing in either a \emph{localized phase}, with insulating properties, or a \emph{delocalized phase}, with conducting properties. This work considers a fundamental example of a random Hamiltonian in the delocalized phase, the Wigner matrix, which is a real symmetric (or complex Hermitian) random matrix whose entries are independent up to the symmetry constraint. First introduced by Wigner in his pioneering investigations of heavy atomic nuclei, it represents a mean-field quantum system where transitions are possible between any two states, with uniformly comparable transition probabilities \cites{wigner1955characteristic,wigner1967random,wigner1958distribution}.

To illustrate the delocalization phenomenon, we consider the Gaussian Orthogonal Ensemble (GOE). We recall this ensemble is defined as the $N \times N$ real symmetric random matrix $\mathrm{GOE}_N = \{ g_{ij} \}_{1 \le i,j \le N}$ whose upper triangular entries $g_{ij}$ are mutually independent Gaussian random variables with variances $( 1 + \mathds{1}_{i =j}) N^{-1}$. It is well known that the distribution of $\mathrm{GOE}_N$ is invariant under conjugacy by orthogonal matrices. Therefore, its $\ell^2$-normalized eigenvectors are rotationally invariant and uniformly distributed on the unit sphere $\S^{N-1}$. Using this observation, it is straightforward to prove (see, for instance, \cite{o2016eigenvectors}*{Theorem 2.1}) that if $\u$ is any eigenvector of the $\mathrm{GOE}_N$, then

\begin{equation}
\label{e:delocintegrable}
\P\left(
	\Vert \u\Vert_\infty
	\geqslant
	\sqrt{
	\frac{2K^3\log N}{N}	
	}
\right)
\leqslant 2N^{1-K}+\exp\left( {-\frac{(K-1)^2}{4K}N}\right)
\end{equation} 
for any $K>1$. 
The bound \eqref{e:delocintegrable} shows that any eigenvector $\u$ of $\mathrm{GOE}_N$ is strongly delocalized, meaning its mass is distributed approximately equally among its coordinates. Both the growth rate $\sqrt{\frac{\log N}{N}}$ and the constant $\sqrt{2}$ are optimal; this can be seen by approximating entries of vectors uniformly distributed on the sphere by independent Gaussian random variables \cites{diaconis87dozen,jiang2005maxima,donoho2001uncertainty}. 

\begin{figure}[!ht]
	\centering
	\makebox[\textwidth][c]{
	\includegraphics[width=1.2\linewidth]{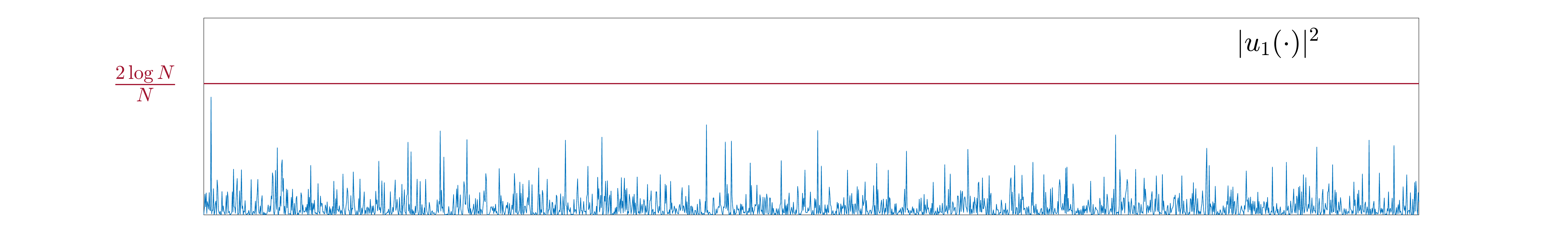}}\\
	\makebox[\textwidth][c]{
	\includegraphics[width=1.2\linewidth]{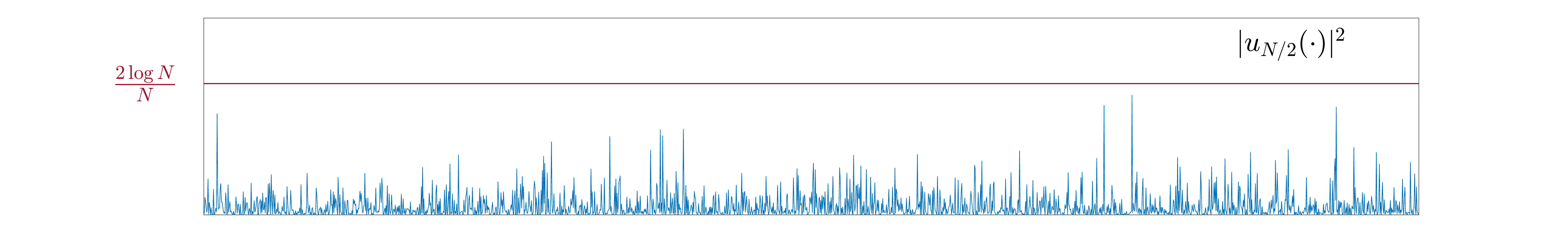}}
	\caption{The delocalization phenomenon for the first and middle eigenvectors of a symmetric $2000 \times 2000$ matrix with centered $\pm 1$ entries.}\label{f:deloc}
\end{figure}

It was conjectured in \cite{o2016eigenvectors} that delocalization with rate $ \sqrt{\frac{\log N}{N}}$ persists when the Gaussian entries $g_{ij}$ are replaced with any subexponential distribution. Figure \ref{f:deloc} illustrates the conjecture for a Bernoulli random matrix. For general distributions, it is no longer possible to appeal to rotational invariance, and an entirely different method is needed. The goal of this paper is to establish the optimal rate in this case, and additionally capture the optimal constant. We also prove isotropic versions of these results.

\subsection{Main results}

We first define generalized Wigner matrices with uniformly subexponential entries. 

\bed \label{d:wigner} A generalized Wigner matrix $H$ is a real symmetric or complex Hermitian $N\times N$ matrix whose upper triangular elements $\{h_{ij}\}_{i\le j}$ are independent random variables with mean zero and variances $\sigma_{ij}^2=\E(|h_{ij}|^2)$ that satisfy 

\beq \label{e:stochasticvar}
\sum_{i=1}^N \sigma_{ij}^2 =1 \quad \text{ for all } j\in \llbracket1,N\rrbracket\eeq
and
\beq 
\frac{c}{N} \le \sigma_{ij}^2 \le \frac{C}{N} \quad  \text{ for all } i,j \in \llbracket1,N\rrbracket\eeq
for some constants $C, c >0$. 
Further, there exists a constant $d >0$ such that, for all $1 \le i ,j \le N$ and $t > 1$, 
\beq \label{e:subexphypo}
\P\left( |h_{ij} | > t \sigma_{ij} \right) \le d^{-1} \exp\left( - t^{d} \right).
\eeq
\eed

Our first main result shows that delocalization with the optimal rate $ \sqrt{\frac{\log N}{N}}$ holds for these matrices with an error probability that may be taken as small as any power of $N$. It is analogous to \eqref{e:delocintegrable} with $K$ large. We also provide an isotropic version.

\bet\label{t:main1}
Let $H$ be a generalized Wigner matrix. Then for every $D>0$, there exists $C = C(D) >0$ such that 
\beq\label{e:isotropicmain2}
 \P\left(  \sup_{\ell \in \unn{1}{N} }  \| {\u_\ell} \|_\infty \ge  C \sqrt{\frac{\log N}{N}} \right) \le  CN^{-D}
\eeq
and
\beq\label{e:isotropicmain}
\sup_{\q \in \S^{N-1}} \P\left( \sup_{\ell \in \unn{1}{N} } \left| \scp{\q}{\u_\ell} \right| \ge   C \sqrt{\frac{\log N}{N}} \right) \le  CN^{-D}.
\eeq
\eet

Our second main result captures the optimal delocalization constant for the $\ell^\infty$ norm of an eigenvector. It corresponds to \eqref{e:delocintegrable} with $K$ close to $1$. We again obtain an isotropic version of this statement where we give an upper bound on any individual eigenvector projection. Since this is not an extreme statistics, we do not expect the constant to be optimal in this case.

\bet\label{t:main2}
Let $H$ be a generalized Wigner matrix. For every $\eps >0$, there exist constants $C=C(\eps)>0$ and $c=c(\eps)>0$ such that 
\begin{align}\label{e:t2claim1}
\P\left(
		\Vert\u_\ell \Vert_\infty
	\geqslant
	\sqrt{\frac{(2+\varepsilon)\log N}{N}}
\right)
\leqslant CN^{-c}
\end{align}
and
\begin{align}\label{e:t2claim2}
\sup_{\q\in\S^{N-1}}
\P\left(	
		\left| \scp{\q}{\u_\ell} \right|
	\geqslant
	\sqrt{\frac{(2+\varepsilon)\log N}{\e N}}
\right)
\leqslant CN^{-c}
\end{align}
for all $\ell \in \unn{1}{N}$ when $H$ is real symmetric. If $H$ is complex Hermitian, the $2+\eps$ is replaced by $1+\eps$.
\eet
Finally, for generalized Wigner matrices with smoother entry distributions, we are able to give an optimal form of delocalization for the maximal entry of the whole eigenbasis. The delocalization constants in this theorem are also sharp (see \cite{jiang2005maxima}*{Proposition 1}).
\bet\label{t:main3}
	Let $H$ be a generalized Wigner matrix and $\vartheta>0$. Assume that $\sqrt{N}H_{ij}$ has a density $\e^{-V_{N,i,j}(x)}\, \d x$ such that for any $k\geqslant 0$, there exists $C>0$ such that 
	\beq
	\vert V^{(k)}_{N,i,j}(x)\vert \leqslant CN^{k(1/6-\vartheta)}(1+\vert x\vert)^C
	\eeq 
	for all $i,j\in\unn{1}{N}$ and $x\in\mathbb{R}$. Then for every $\varepsilon>0$, there exist constants $C=C(\varepsilon)>0$ and $c=c(\varepsilon)>0$ such that 
	\begin{align}\label{heatconclude}
	\P\left(
		\sup_{\ell\in\unn{1}{N}}\Vert \u_\ell\Vert_\infty \geqslant 
		\sqrt{\frac{(4+\varepsilon)\log N}{N}}
	\right)\leqslant CN^{-c}
	\end{align}
	when $H$ is real symmetric. If $H$ is complex Hermitian, the $4+ \eps$ is replaced by $2+\eps$.
\eet

\subsection{Background}The literature on eigenvector delocalization is extensive and we review here only the previous works most relevant to our results. For a broader overview, the interested reader may consult a number of recent surveys \cites{o2016eigenvectors, erdos2017dynamical, bourgade2018survey}. 

The first strong delocalization bounds were achieved in the seminal papers \cites{erdos2009semicircle, erdos2009local, erdos2010wegner}. They showed that $\Vert \sqrt{N}\u_\ell\Vert_\infty\leqslant (\log N)^{4}$ with very high probability, for matrices with independent, identically distributed entries, but required that the entries posses a density satisfying a certain smoothness assumption, and that the eigenvectors correspond to eigenvalues in the bulk of the spectrum. Their method uses resolvent estimates to prove a \emph{local semicircle law} on small scales by a descending bootstrap argument. This law provides a quantitative rate for the convergence of the empirical spectral distribution to its limit on small  intervals (shrinking in $N$). Refinements of this approach to establishing delocalization appear in a number of later works. The smoothness assumption was relaxed in \cite{tao2011random}, which gave an upper bound of $(\log N)^{20}$, and delocalization was extended to eigenvectors corresponding to eigenvalues at the edge of the spectrum in \cite{tao2010random} with a bound of $(\log N)^{C}$ for some $C>0$. For generalized Wigner matrices with subexponential decay, as defined in \Cref{d:wigner}, delocalization was proved in \cite{erdos2012bulk} with rate $(\log N)^{13+6/d}$ for bulk eigenvectors. It was extended in \cite{erdos2012rigidity} to all eigenvectors with an upper bound of roughly $(\log N)^{\log\log N}$. 

The first optimal bound was proved in \cite{vu2015random} using a different technique, which centers on the study of random weighted projections in high dimensions. Using it, the authors were able to establish the local semicircle law on scales even smaller than those previously accessible. As a consequence, for matrices with independent, identically distributed, and bounded entries, they obtained the $\sqrt{\log N}$ rate for bulk eigenvectors and a $\log N$ rate for edge eigenvectors. The boundedness hypothesis can be relaxed to requiring only that the entries are sub-Gaussian by combining their method with \cite{rudelson2013hanson}, as remarked in \cite{o2016eigenvectors}*{Section 4}. For uniformly subexponential entries, as defined in \eqref{e:subexphypo}, the same technique  bounds $\Vert \sqrt{N}\u_\ell\Vert_\infty$ by $\log^{1/2+1/d}N$ for bulk eigenvectors and $\log^{1+2/d}N$  for eigenvectors close to the spectral edges  \cite{o2016eigenvectors}*{Corollary 4.4}. The method of \cite{vu2015random} was later applied in \cite{dumitriu2018sparse} to reach the optimal $\sqrt{\log N}$ delocalization rate for Wigner-type matrices with sub-Gaussian entries in the bulk. 

These results concerned delocalization for the extremal coordinates of eigenvectors, as in \eqref{e:isotropicmain2}. The first result regarding isotropic delocalization, as in \eqref{e:isotropicmain}, gave a bound of $(\log N)^{\log\log N}$ for Wigner matrices \cite{knowles2013isotropic}. This result was extended in \cite{bloemendal2014isotropic} to generalized Wigner matrices; they gave a bound of $N^\varepsilon$ for any $\epsilon>0$, for all eigenvectors and sufficiently large $N$

While all these results assume the existence of arbitrarily large moments of the entries (at least), recent works have aimed to weaken this assumption and consider distributions with heavier tails. In \cite{aggarwal2019bulk}, complete delocalization in the bulk was proved with rate $(\log N)^{C}$ under the existence of $2+\eps$ moments for generalized Wigner matrices. In \cites{gotze2018local, gotze2019local}, the optimal rate $\sqrt{\log N}$ for all eigenvectors was shown under the existence of a fourth moment, but with a weak probability bound and the requirement that the entries have identical variances.

From this overview, we  see that our improvements to the state of the art are threefold:
\begin{itemize}
	\item \emph{Optimal delocalization rate with very high probability throughout the spectrum:} \Cref{t:main1} is the first delocalization result for generalized Wigner matrices that controls \emph{all} eigenvectors with the optimal rate $\sqrt{\log N}$ and probability $1-N^{-D}$ for any $D>0$. This estimate permits a union bound over all eigenvectors to control the entire eigenbasis with rate $\sqrt{\log N}$ and very high probability; such a result is inaccessible with the weaker probability bounds obtained in \cites{gotze2018local, gotze2019local}. Another difference is that we \emph{do not} require that the entries are identically distributed or have equal variances.
	
	\item \emph{Optimal constants:} In \Cref{t:main2}, we prove that the constant $\sqrt{2}$ is a universal upper bound for the (appropriately rescaled) maximal entry of eigenvectors of symmetric generalized Wigner matrices. We also obtain the analogous result with constant $1$ for Hermitian matrices. With a more restrictive hypothesis on the entries, we also obtain the optimal constant for the upper bound on the entire eigenbasis in \Cref{t:main3}. Previously, sharp constants were available only for integrable random matrix ensembles. 
	
	\item \emph{Optimal isotropic delocalization:} We give the first isotropic delocalization estimate that achieves the optimal rate $\sqrt{\log N}$.

	\end{itemize}

While our results give an upper bound on the maximal coordinates of eigenvectors, there are other ways to describe the delocalization phenomenon. For instance, the asymptotic distribution of eigenvector entries was studied in \cites{knowles2013eigenvector, tao2012random, bourgade2013eigenvector, marcinek2020high}, where it was proved for generalized Wigner matrices that any fixed number of eigenvector entries are asymptotically independent Gaussian random variables. The entry distribution was also studied in \cite{aggarwal2020eigenvector} for L\'evy matrices, whose entries have infinite variance, in the delocalized phase of the spectrum. It was shown there that the asymptotic distribution of an entry arises from a one-parameter family of non-Gaussian laws, with the exact distribution determined by the location of the corresponding eigenvalue.

\emph{Quantum unique ergodicity} is another form of delocalization and consists of a concentration bound for an eigenvector's $\ell^2$ mass averaged over a subset of coordinates. Such a result for generalized Wigner matrices was first given in \cite{bourgade2013eigenvector} and was later strengthened in \cite{bourgade2018random} for Gaussian divisible ensembles, where it was used to prove universality results for random band matrices. A strong form of quantum unique ergodicity for all generalized Wigner matrices was subsequently shown in \cite{benigni2021fluctuations}.
\emph{No gaps delocalization} was proved in \cite{rudelson2016no}, where the authors show that any subset of coordinates with size linear in $N$ must contain some non-negligible $\ell^2$ mass. A related result establishing a lower bound on the smallest coordinate of an eigenvector can be found in \cite{o2016eigenvectors}*{Theorem 4.7}.

Finally, for the eigenvectors of non-symmetric matrices, \cite{rudelson2015delocalization}*{Theorem 1.1} proved a delocalization bound with rate $\log^{9/2} N$ for matrices with independent, mean zero entries with variances at least 1 using a novel geometric method. In particular, they do not require that the entries have identical distributions, or even the same variance, and the bound holds with error probability $N^{-D}$ for any $D>0$. More general models and isotropic delocalization were considered in \cites{alt2019spectral, alt2020inhomogeneous}.

\subsection{Proof strategy}
The proofs of our results are based on the dynamical approach to random matrix theory, which was introduced in \cite{erdos2011universality} for proving the universality of local eigenvalue statistics. The proofs all broadly consist of the following three steps; our novel contributions come in the second and third steps.

	\paragraph{{Step 1: Rigidity and isotropic local law.}}To begin our analysis, we require \emph{a priori} estimates on the eigenvalues and eigenvectors of our original matrix. The \emph{rigidity} of eigenvalues for generalized Wigner ensembles was proved in \cite{erdos2012rigidity} and states that with very high probability, all eigenvalues are close to their typical locations (which are deterministic). The \emph{isotropic local law} for generalized Wigner matrices, proved in \cites{knowles2013isotropic, bloemendal2014isotropic}, controls the resolvent of a generalized Wigner matrix as a quadratic form, and as a consequence bounds the inner products of eigenvectors with arbitrary unit vectors. We remark that we do not require that the local law holds down to the optimal scale $\frac{\log N}{N}$, since our optimal bounds instead arise from the next step; instead, we only need it on the scale $N^{-1 + \eps}$ for arbitrary $\eps >0$.

	\paragraph{{Step 2: Relaxation by the Dyson Brownian motion.}}We next use the \emph{a priori} estimates on  generalized Wigner matrices from the first step to study the eigenvectors of such matrices perturbed by a small additive Gaussian noise matrix, which we write as $\sqrt{t} \mathrm{GOE}_N$ for a parameter $t \ll 1$. It was shown in \cite{bourgade2013eigenvector} that if $t$ is viewed as a time parameter, the time evolution of moments of the eigenvector entries of $H + \sqrt{t} \mathrm{GOE}_N$ is governed by a system of differential equations with random coefficients (more precisely a random walk in a dynamic random environment) now known as the \emph{eigenvector moment flow}. While this flow has been used to derive numerous eigenvector statistics for different models \cites{bourgade2017eigenvector, benigni2017eigenvectors, bourgade2018random, benigni2019fermionic, aggarwal2020eigenvector, marcinek2020high}, we refine the previous analyses to handle growing moments for the first time; we permit moments which may grow as fast as $(\log N)^C$ for any $C>0$. Then, through a Markov inequality, we obtain the optimal delocalization estimates from these moment bounds. We emphasize that our strong control over these growing moments is the underlying mechanism that enables us to reach the optimal delocalization rate and optimal constants in our main theorems.
	
	\paragraph{{Step 3: Regularization and comparison.}} Finally, we extend the results from the previous step to show that the optimal delocalization bounds hold for \emph{all} generalized Wigner matrices, not just those with additive noise. This is the most technical step, to which the majority of the paper is devoted.
	
	At this juncture, the standard approach in the literature is the four moment method, which was first applied in random matrix theory in \cites{tao2011random,tao2010random}. The basic insight is that if two random matrices have entries whose first four moments match, then any sufficiently regular observable takes, in expectation, the same value for both matrices. To prove this, one replaces the entries of the first matrix with the second one by one, showing at each step that the given observable changes by a negligible amount, and then sums the error terms across all $\mathcal O (N^2)$ replacements. The replacements are accomplished by Taylor expanding the observable in the matrix entries; the four moment matching condition is exactly what is needed for the resulting error bounds to be effective. At a high level, this approach parallels exactly Lindeberg's proof of the central limit theorem.
	
	By standard results, given a generalized Wigner matrix $H$, it is possible to find a generalized Wigner matrix $M$ and a short time $t \ll 1$ such that the first four moments of the entries of $M + \sqrt{t} \mathrm{GOE}_N$ match those of $H$ (at least asymptotically). Given this, one hopes that the extension from the ensembles considered in the second step to all generalized Wigner matrices may be accomplished by the four moment method. However, we arrive at a central difficulty: the eigenvectors of a random matrix are highly singular observables, and it is difficult to control their derivatives with respect to matrix entries, which blocks control of the error terms in the Taylor expansion mentioned above. The same problem afflicts the other common approach to this step, the matrix continuity estimate of \cite{bourgade2013eigenvector}, which compares $H$ to $H + \sqrt{t} \mathrm{GOE}_N$.
	
	To overcome this obstacle, we introduce \emph{regularized eigenvectors}, which are smooth versions of the usual eigenvectors with the advantage that they are amenable to the four moment comparison. Our construction builds on the regularization introduced in \cite{knowles2013eigenvector}, and we discuss this work below, after completing our proof sketch.
	
	Let $H$ be a Wigner matrix with eigenvalues $\{ \lambda_i \}_{1 \le i \le N }$, and let $\u_i$ be the normalized eigenvector corresponding to $\lambda_i$. Fix $\ell\in \unn{1}{N}$, and let $\eta >0$ be a parameter smaller than the typical eigenvalue spacing around $\lambda_\ell$. Because the Poisson kernel integrates to $1$, we have
	\beq
	| u_\ell (i ) |^2   = \frac{\eta}{\pi} \int_{\R} \frac{ | u_\ell (i ) |^2\, \d E} {(E - \lambda_\ell)^2 + \eta^2} \approx 	\frac{\eta}{\pi} \int_{I} \frac{ | u_\ell (i ) |^2\, \d E} {(E - \lambda_\ell)^2 + \eta^2},
	\eeq
	where $I$ is an interval with length slightly larger than $\eta$ centered at $\lambda_\ell$. In the approximation, we removed a negligible portion of the integral, since the Poisson kernel is concentrated on scale $\eta$ around $E$. If $G(z) = (H -z)^{-1}$ is the resolvent for $H$, defined for $z\in \C$, then the spectral theorem gives 
	\beq\label{e:theapprox}
	\int_{I} \frac{ | u_\ell (i ) |^2\, \d E} {(E - \lambda_\ell)^2 + \eta^2} \le \int_{I} \sum_{j} \frac{ | u_j (i ) |^2\, \d E} {(E - \lambda_j)^2 + \eta^2} = \int_{I} \Im G_{ii} (E + \I \eta) \, \d E.
	\eeq
	It is therefore natural to use the last term in \eqref{e:theapprox} to define a  regularized eigenvector entry which approximates  $u_\ell(i)$. However, one must be careful, since in general this term could be much larger than $u_\ell(i)$ if several eigenvalues cluster in the small interval $I$ around $\lambda_\ell$ (as can be seen directly from the spectral expansion). Indeed, it is only possible to exclude eigenvalues other than $\lambda_\ell$ from $I$ with probability at most $1 - N^{-c}$ for a small constant $c>0$; while this suffices for \Cref{t:main2} in combination with the appropriate observable (described below), it is far from the high probability bounds we require for \Cref{t:main1}.
	
	However, we make the following observation: it is not necessary for the proof of \Cref{t:main1} that the regularized eigenvector be approximately equal to the usual one, only that it be \emph{comparable} up to some constant factor. As noted below, we are able to show that for any $D > 0$, there exists a constant $C$ such that there are no more than $C$ eigenvalues in $I$ with probability at least $1 - N^{-D}$. This enables the very high probability comparison for \Cref{t:main1}, which is based on comparing growing moments of regularized eigenvectors.
	
	While the high moment comparison we use to prove \Cref{t:main1} is effective for obtaining very high probability bounds, it is too crude to preserve the optimal constants. For \Cref{t:main2} we use the following smoothed maximal function as the observable in the four moment method, which was introduced in \cite{landon2018comparison} to study extremal statistics of eigenvalues:
\beq\label{e:freeenergy}
\frac{1}{\beta}\log\left(
	\sum_{i\in\unn{1}{N}}
	\exp(\beta u_\ell^2(i))
\right)
\approx
\Vert \u_\ell\Vert_\infty^2,
\eeq
where the approximation holds for large $\beta$ (see \eqref{e:thefirstestimate} below).
 The observable \eqref{e:freeenergy} should be thought of as a free energy functional with inverse temperature parameter $\beta$. The motivation for this choice comes from statistical mechanics, where it is well known that the free energy at low temperature is close to the ground state. 
	
	Finally, we note that the comparison method sketched for \Cref{t:main2} does not quite suffice for \Cref{t:main3}, since in the latter case it would require excluding eigenvalues in sub-microsopic intervals centered around every eigenvalue $\lambda_i$ simultaneously, which is impossible. We therefore resort to the reverse heat flow method introduced in \cite{EPR10}, which provides a strong comparison bound, but comes at the cost of requiring that the matrix entries have densities satisfying the smoothness hypothesis appearing in the statement of \Cref{t:main3}.
	
	\paragraph{} We now note two important technical points relevant to step 3 of the above outline. First, in order to control the number of eigenvalues in the small interval $I$ centered at $\lambda_\ell$, we require bounds known as \emph{level repulsion estimates}. Suppose for simplicity that $\lambda_\ell$ is in the bulk, $\lambda_\ell\in(-2+\kappa, 2-\kappa)$ for some $\kappa\in(0,2)$, so that the typical eigenvalue spacing is $N^{-1}$. To exclude eigenvalues other than $\lambda_\ell$ from $I$, we require that for small enough $\delta>0$, there exists $\alpha= \alpha(\delta)>0$ such that 
\beq\label{e:lre2intro}
\P\left(
	\#\{j\in\unn{1}{N}, \lambda_j\in(\lambda_\ell-cN^{-1-\delta},\lambda_\ell+cN^{-1-\delta})\}\geqslant 2
\right)\leqslant N^{-\alpha}
\eeq
for any constant $c>0$ and large enough $N$ (depending on $\delta$ and $c$). For the high probability bound, we must show that
for all $D>0$, there exists a $k\in\N$ such that 
\beq\label{e:overcrowd}
\P\left(
	\#\{j\in\unn{1}{N}, \lambda_j\in(\lambda_\ell-cN^{-1-\delta},\lambda_\ell+cN^{-1-\delta})\}>  k
\right)\leqslant N^{-D}
\eeq
for large enough $N$.
In the latter case, such bounds are also known as \emph{overcrowding estimates}. 

In the bulk, \eqref{e:overcrowd} was first proved in \cite{erdos2010wegner}*{Theorem 3.5} for Hermitian matrices, and the same argument was later adapted to symmetric matrices in \cite{bourgade2016fixed}*{Appendix B}. The estimate \eqref{e:lre2intro} follows from the bulk universality result in \cite{erdos2012gap}. Further, both estimates, for both symmetry classes, also follow from the main results of \cite{nguyen2018random}, which appeared later. At the edge, \cite{bourgade2014edge}*{Theorem 2.7} showed \eqref{e:lre2intro} for eigenvalues $\lambda_i$ such that $i \in \unn{1}{N^{1/4}} \cup \unn{N -N^{1/4}}{N}$, where we order the eigenvalues from least to greatest. To the best of our knowledge, level repulsion estimates for adjacent eigenvalues have not been established for the intermediate regime between the bulk and the edge eigenvalues covered by \cite{bourgade2014edge}, and overcrowding estimates have not been established outside of the bulk.

However, to complete step 3 of the outline for all eigenvectors, we require both kinds of estimates to hold throughout the spectrum. Therefore, in \Cref{sec:lre}, we prove \Cref{l:lr2}, a level repulsion estimate for adjacent eigenvalues which holds with no restriction on the location of the eigenvalues, and \Cref{p:lrmain}, an overcrowding estimate for sub-microscopic intervals anywhere in the spectrum. \Cref{l:lr2} is a straightforward consequence of the result \cite{bourgade2014edge} mentioned above and the analysis of Dyson Brownian motion given in \cite{bourgade2018extreme}, which gives gap universality throughout the spectrum.
The proof of \Cref{p:lrmain} is somewhat more involved. It proceeds by first establishing the result for generalized Wigner matrices with small Gaussian noise, then using a moment matching argument on (moments of) a smoothed eigenvalue counting observable to extend the bound to all generalized Wigner matrices. For the ensembles with Gaussian noise, we deal with the bulk of the spectrum by proving a slight generalization of a similar result from \cite{landon2017convergence}. At the edge, we first obtain the result for the $\mathrm{GOE}_N$ and its Hermitian counterpart using determinantal estimates, then use a dynamical result from \cite{bourgade2018extreme} to show it holds for any ensemble with a small Gaussian component.

The level repulsion estimate \Cref{l:lr2} also has an interesting consequence for the distribution of the eigenvectors. The proof of the asymptotic normality of the eigenvector entries given in \cite{bourgade2013eigenvector} relies on level repulsion estimates for neighboring eigenvalues. However, because of the lack of such an estimate in the intermediate regime noted above, it was only possible to establish asymptotic normality for $\lambda_i$ such that $i \in \unn{1}{N^{1/4}}\cup \unn{N^{1-c}}{ N - N^{1-c}} \cup \unn{N -N^{1/4}}{N} $, where $c>0$ is a small constant. Our \Cref{l:lr2} thus immediately implies the result for the entire spectrum after combining it with the argument in  \cite{bourgade2013eigenvector}; we state this as \Cref{c:asymptoticallynormal}.

Second, the regularization given in the heuristic \eqref{e:theapprox} is not actually smooth enough for our purposes. The integral there is over an interval $I$ centered at $\lambda_\ell$, and to differentiate this integral with respect to a matrix entry also involves differentiating $\lambda_\ell$. However, like the eigenvectors, the eigenvalues are too singular to permit the four moment comparison argument to go through. We therefore construct in \Cref{a:eigreg} \emph{regularized eigenvalues}, which serve as smooth counterparts to the usual eigenvalues (at least on a set of very high probability) and are suitable substitutes for the usual eigenvalues in the comparison. Such a construction was already performed for bulk eigenvalues in  \cite{landon2018comparison}*{Lemma 3.2}. Our method, which applies to all eigenvalues, is essentially the same, but the proof now requires a more careful treatment of the error terms to accommodate the varying inter-eigenvalue distances for eigenvalues outside the bulk.

\paragraph{} We would like to acknowledge the fundamental work \cite{knowles2013eigenvector}, which inspired the comparison strategy used in step 3 of the outline above. In \cite{knowles2013eigenvector}, the authors provide a moment matching scheme for eigenvector observables which also uses a resolvent-based eigenvector regularization and level repulsion estimates. Using this framework, they prove, \emph{inter alia}, that only two matching moments are needed for the comparison of edge eigenvectors, which implies universality for the entries of edge eigenvectors for generalized Wigner matrices. However, \cite{knowles2013eigenvector} considers a comparison with error probability $N^{-c}$ for some small $c>0$, instead of the $N^{-D}$ for any $D>0$ obtained in \Cref{t:main1}, and only eigenvectors corresponding to eigenvalues in the bulk and very near the edge (\cite{knowles2013eigenvector}*{Theorem 1.6} is stated for eigenvectors $\u_i$ for $i \in \unn{1}{(\log N)^{ C \log \log N}}$ for any fixed $C>0$). Further, only observables depending on a finite number of eigenvector entries are permitted. This contrasts with our delocalization results, which control $\mathcal O(N)$ entries simultaneously. We also note that results similar to those in \cite{knowles2013eigenvector} were proved soon after in \cite{tao2012random}, with the hypothesis of four matching moments.

\paragraph{} 
We have stated our results in the context of generalized Wigner matrices mainly for brevity. The essential inputs to the argument are the \emph{a priori} estimates listed in step 1 of the proof outline, and it is known that these hold for a broad class of mean-field ensembles \cites{aggarwal2019bulk, ajanki2017universality, aggarwal2018goe, EKYY13, che2019universality}, so it is likely our techniques apply more widely. In particular, the argument leading to the sharp constant in \eqref{t:main2} is quite general and does not require the subexponential decay hypothesis on the entries (essentially because the final comparison step uses the free energy functional, instead of controlling growing moments). We note that the constant $2$ is not correct in the general case and will depend on the asymptotic law of the spectrum for random matrices with a non-stochastic variance profile. Further, our results on level repulsion and overcrowding in \Cref{sec:lre}, our regularization of eigenvalues throughout the spectrum in \Cref{a:eigreg}, and the very high probability comparison scheme from \Cref{s:main} may also be of independent interest.

Finally, \Cref{t:main2} provides a first step to proving the universality of the distribution of the maximal entry of an eigenvector. For the $\mathrm{GOE}_N$, \cite{jiang2005maxima} proved the matching lower bound $\sqrt{N/\log N}\left\Vert \u \right\Vert_\infty \xrightarrow[]{\mathds{P}}\sqrt{2}$ and obtained
\[
\P\left(
	N\sup_{\ell\in\unn{1}{N}}\Vert \u_\ell\Vert_\infty^2-4\log N+\log\log N+\log(2\pi)\leqslant 2x
\right)\xrightarrow[N\to\infty]{}\e^{-\e^{-x}}.
\] 
We hope to address this question in future work.

\subsection{Outline of the paper} We state the preliminary estimates we need from previous works in \Cref{s:prem}, such as local laws, rigidity, and delocalization. In \Cref{sec:dynamics}, we define the eigenvector moment flow and perform the high moment analysis of the flow described in step 2 of our outline. In \Cref{s:eigvectreg}, we develop our eigenvector regularization scheme, assuming the eigenvalue regularization result proved in \Cref{a:eigreg} and the level repulsion estimates proved in \Cref{sec:lre}. Finally, \Cref{s:main} is dedicated to the proofs of \Cref{t:main1}, \Cref{t:main2}, and \Cref{t:main3}. \Cref{a:eigreg}, \Cref{sec:lre}, and \Cref{a:preliminary} prove the preliminary results necessary to regularize the eigenvectors.

\paragraph{Acknowledgments.} The authors thank P. Bourgade for suggesting the current problem and many helpful conversations. P.L. thanks A. Aggarwal, B. Landon, P. Sosoe, and H.-T. Yau for helpful conversations.

\section{Preliminaries}\label{s:prem}

In this section, we introduce some preliminary notions and results from previous works that are used throughout the rest of the paper.

We say $X \ll Y$ if there exists a small constant $c > 0 $ such that $N^c |X| \le Y$. 
We write $X = \O{Y}$ if there exists $C> 0$ such that $|X| \le C Y$. Here $X$ and $Y$ may depend on other parameters, but $C$ does not. We also say $X = \mathcal O_u(Y)$ if $|X| \le C_u |Y|$ for some constant $C_u >0$ depending only on a parameter $u$. The same notation with multiple subscripts denotes dependence on multiple parameters. The notation $\log N$ always denotes the natural logarithm.

Throughout this work, we suppress the dependence of various constants in our results on the constants in \Cref{d:wigner}. This dependence does not affect our arguments in any substantial way. Further, for concreteness we consider only real symmetric generalized Wigner matrices in what follows, except where noted, as the complex Hermitian case is entirely analogous.

Let $\matn$ be the set of $N\times N$ real symmetric matrices. We label the eigenvalues of matrices in increasing order, so that $\lambda_1 \le \lambda_2 \le \dots \le \lambda_N$ for the eigenvalues of a matrix $M \in \matn$. We define the Stieltjes transform of $M$ for $z \in \C \setminus \{ \lambda_1, \dots, \lambda_N\}$ by
\begin{align}
m_N (z) = \frac{1}{N} \sum_i \frac{1}{ \lambda_i - z}.
\end{align}
The semicircle law and corresponding Stieltjes transform are given by
\begin{align}
\scrho (E) = \frac{\sqrt{ (4 - E^2)_+ }}{2 \pi }, \qquad \msc (z) = \int_{\R} \frac{\scrho (x)\, \d x}{ x - z },
\end{align}
for $E \in \R$ and $z \in \C \setminus \R$, respectively.

The classical eigenvalue locations of the semicircle law are denoted by $\gamma_i$ and defined for $i \in \llbracket 1, N \rrbracket$ by
\beq\label{e:classical}
\frac{i}{N} = \int_{-2}^{\gamma_i} \scrho (x) \, \d x.
\eeq
The resolvent of $M\in \matn$ is given by $G(z) = ( M - z \Id)^{-1} $. We observe by the spectral theorem that
\beq\label{eq:defsresolv}
G (z) =
\sum_{i=1}^N
\frac{ \u_i  \u^*_i }{\lambda_i-z},
\eeq
where $\u_i$ is the $L^2$-normalized eigenvector, $\| \u_i \|_2 = 1$, corresponding to the eigenvalue $\lambda_i$ of $M$, and $\u^*_i$ denotes its transpose.\footnote{Note that $\u_i$ is ambiguous up to a choice of sign, which we select arbitrarily. It is straightforward to see that all of our definitions are invariant under this choice.}


\subsection{Dyson Brownian motion}\label{s:dbm}

The $N\times N$ real symmetric Dyson Brownian motion with initial value $H_0$ is the stochastic process $(H_s)_{0\le s \le 1}$ on the space of symmetric matrices defined by 
\beq\label{e:dysondyn}
 H_s =H_0 +  \frac{1}{\sqrt{N}} B(s),
\eeq
 where $B(s) \in \matn$ is a symmetric matrix such that $B_{ij}(s)$ and $B_{ii}(s)/\sqrt{2}$ are mutually independent standard Brownian motions for all indices $1 \le i < j\le N$. 

It is known that $(H_s)_{0\le s \le 1}$ has the same distribution as $\left( \u^s \bm \lambda(s) (\u^s)^* \right)_{0 \le s \le 1}$, where $\bm \lambda(s) = ( \lambda_1(s), \dots, \lambda_N(s))$ is a vector in $\R^N$ and $\u^s = \left( \u^s_1, \dots, \u^s_N \right) \in \R^{N \times N}$ are the stochastic processes that solve the system of equations

\begin{align}
\label{eq:dysonval}\d \lambda_k(s) &= \frac{\d \wt{B}_{kk}(s)}{\sqrt{N}}
+\left(
	\frac{1}{N}
	\sum_{\ell\neq k}
	\frac{1}{\lambda_k(s)-\lambda_\ell(s)}
\right)\d s, \\
\label{eq:dysonvect}\d \u^s_k &= \frac{1}{\sqrt{N}}\sum_{\ell\neq k}\frac{\d \wt{B}_{k\ell}(s)}{\lambda_k(s)-\lambda_\ell(s)}\u_l^s
-\frac{1}{2N}
\sum_{\ell\neq k}
\frac{\d t}{(\lambda_k(s)-\lambda_\ell(s))^2}\u_k^s,
\end{align}
with initial data $H_0 =\u^0\bm  \lambda(0) (\u^0)^*$, and $\wt{B}$ has the same distribution as $B$ \cite{bourgade2013eigenvector}*{Theorem 2.3}. 
We define $m_N^s(z)$ to be the Stieltjes transform of $H_s$, and $G^s(z)$ to be its resolvent.

\subsection{Local semicircle law}

Let $\mathbb{S}^{N-1} \subset \mathbb R^N$ be the set of vectors $\q \in \R^N$ such that $\| \q \|_2 =1$. For any $\omega >0$, we define the domain
\beq \D_\omega  =  \{ z = E + \I\eta \in \mathbb C : |E| < \omega^{-1}, N^{-1+ \omega} \le \eta \le \omega^{-1}  \}. 
\eeq

\bel[\cite{bourgade2013eigenvector}*{Lemma 4.2}]\label{l:goodset}
Fix  $\omega >  0$ and $\q\in\S^{N-1}$. 
Let $H$ be a $N\times N$ generalized Wigner matrix, and define $H_s$, $G^s$, $\u^s_k$, and $m^s_N(z)$ as above. Let $\mu$ be the measure on the space of joint eigenvalue and eigenvector trajectories  $(\bm \lambda(s), \u(s))_{0 \le s \le 1}$ induced by the Dyson Brownian Motion $(H_s)_{0\le s \le 1}$.
Then there exist constants $C_1(\omega), c_1(\omega)> 0 $  such that
\beq\label{e:exphigh}
\inf_{\q \in\mathbb{S}^{N-1}}
\P_\mu \left(
	\mathcal{A}(\q,\omega)
\right)
\geqslant1-C_1\exp\left({-c_1(\log N)^{c
_1\log\log N}}\right),
\eeq
where $\mathcal{A}(\q,\omega)$ is the set of trajectories $(\bm \lambda(s), \u(s))_{0 \le s \le 1}$  where all of the following statements hold.

\begin{enumerate}
\item For all $z=E+\I \eta\in \D_\omega$,
\beq
\label{e:isolaw}
\sup_{s\in[0,1]}
\left\vert
	\scp{\q}{G^s(z)\q}
	-
	\msc(z)
\right\vert
\leqslant
N^{\omega}\left(
	\sqrt{\frac{\Im \msc(z)}{N\eta}}
	+
	\frac{1}{N\eta}
\right).
\eeq
\item For all $z=E+\I \eta\in \D_\omega$,
\beq\label{e:sclaw}
\sup_{s\in [0,1]}
\left\vert
	m^s_N(z)-\msc(z)
\right\vert
\leqslant
\frac{N^\omega}{N\eta}, \qquad
\sup_{s\in[0,1]} \left| G^s_{ij}(z) - \one_{i = j} \msc(z) \right| \le N^{\omega}\left(
	\sqrt{\frac{\Im \msc(z)}{N\eta}}
	+
	\frac{1}{N\eta}
\right).
\eeq
\item For all $k\in\unn{1}{N}$,
\beq\label{e:rigidity}
\sup_{s\in[0,1]}
\left\vert
	\lambda_k(s)-\gamma_k
\right\vert
\leqslant
N^{-2/3+\omega}\left[\min(k, N - k +1) \right]^{-1/3}, \qquad
\sup_{s\in[0,1]}
\scp{\q}{ \u_k^s}^2
\leqslant
N^{-1+\omega}.
\eeq
\end{enumerate}
\eel
\begin{rmk}
This lemma was proved in \cite{bourgade2013eigenvector} in the case where the matrix entry distributions have finite moments, with weaker probability bounds. The proof relies on the isotropic semicircle law stated in \cite{bloemendal2014isotropic}*{Theorem 2.12}, which is the origin of the probability bounds. However, our assumption that the entry distributions are uniformly subexponential enables the stronger probability estimates given in \Cref{l:goodset}, as noted after \cite{bloemendal2014isotropic}*{Theorem 4.1}. See also \cite{knowles2013isotropic}*{Theorem 2.12} for the case where all entries have the same variance.
\end{rmk}

Finally, we require some estimates on the Stieltjes transform of the semicircle law.
\bel[\cite{erdos2017dynamical}*{Lemma 6.2}]\label{l:msclower}
There exists a constant $\scc>0$ such that for all $z = E + \I \eta$ with $\eta \in (0, 10]$, when $E \in [-2, 2]$ we have
\beq\label{e:msclower}
\scc \sqrt{ \left| | E| - 2 \right| +  \eta} \le \Im \msc(z) \le  \scc^{-1} \sqrt{ \left| | E| - 2 \right| +  \eta},
\eeq
and when $E \in [-20, 20] \setminus [-2, 2]$ we have 
\beq\label{e:msclower2}
\frac{\scc  \eta}{\sqrt{ \left| | E| - 2 \right| +  \eta} } \le \Im \msc(z) \le \frac{  \scc^{-1} \eta}{\sqrt{ \left| | E| - 2 \right| +  \eta} }.
\eeq
\eel

\section{Relaxation by Dyson Brownian motion}\label{sec:dynamics}

In this section, we control large moments of eigenvector entries for generalized Wigner matrices with small additive Gaussian perturbations. These perturbations are given by Dyson Brownian motion for times $s \ll 1$. Our main result is \Cref{p:highmom}, which bounds the moments of the eigenvector entries when $s \gg N^{-1/3}$ for moments as large as $( \log N)^C$. Then, in \Cref{c:delocdyson}, we apply \Cref{p:highmom} to give high-probability delocalization estimates for the perturbed ensembles.

\subsection{Eigenvector moment flow}

Let $H$ be a generalized Wigner matrix, and let $H_s$ be Dyson Brownian motion with initial condition $H_0 = H$, as defined in \Cref{s:dbm}. It was shown in \cite{bourgade2013eigenvector} that moments of the eigenvector entries of $H_s$ obey a parabolic equation with random coefficients, called the \emph{eigenvector moment flow}, which we now describe.

To any index set $\{ (i_1, j_1), \dots , (i_m, j_m)\}$ with distinct $i_k \in \llbracket 1, N \rrbracket$ and positive $j_k \in \N$, we may associate the vector $\bm \xi = (\xi_1, \xi_2, \dots, \xi_N) \in \mathbb N^N$ with $\xi_{i_k} = j_k$ for $1 \le k \le m$ and $\xi_p =0$ for $p \notin \{ i_1, \dots, i_m\}$. We think of $\bm \xi$ as a particle configuration on the integers, with $j_k$ particles at site $i_k$ for all $k$ and zero particles on the sites not in $\{ i_1, \dots, i_m\}$. 
The configuration $\bm \xi^{ij}$ is defined as the result of moving one particle in $\bm \xi$ from $i$ to $j$, if this is possible. That is, if $i \ne j$ and $\xi_i >0$, then $\xi^{ij}_k$ equals $\xi_k + 1$ if $k=j$, $\xi_k - 1$ if $k=i$, and  $\xi_k$ if $k \notin \{ i, j \}$. When $\xi_i =0$, we set $\bm \xi^{ij} = \bm \xi$.

The moment observable $f_s\colon \N^N \rightarrow \R$ is defined as follows, where we recall that $\bm \lambda(s)$ and $\u^s$ represent the eigenvalue and eigenvector processes associated to $H_s$ by \eqref{eq:dysonval} and \eqref{eq:dysonvect}. Given $\q\in\S ^{N-1}$ and a path $\bm \lambda = (\bm \lambda(s))_{0\le s \le1}$, we set
\beq
\label{eq:deffs}
f_{s}(\bm{\xi})  = f_{\bm \lambda,s}(\bm{\xi}) 
=
\frac{1}{\M(\bm{\xi})}
\E\left[
	\prod_{k=1}^N \scp{\q}{\sqrt{N}\u^s_k}^{2\xi_k}
	\middle\vert
	\bm{\lambda}
\right],
\quad\text{where}\quad
\M(\bm{\xi})
=
\prod_{k=1}^N (2\xi_k -1 )!!.
\eeq
The normalization factor $\mathcal M(\bm \xi)$ is chosen because $(2\xi_k -1 )!!$ is the $2\xi_k$-th moment of a standard Gaussian.
\begin{figure}[!ht]
	\centering
	\begin{tikzpicture}[scale = 1.4]
		\draw[-] (-1,0) -- (7,0);
		\node[fill=White] at (0,0) {};  
 		\node[label={below:\LARGE $k_1$}] at (1,0) {};
		\node[fill=White] at (2,0) {};
		\node[label=below:\LARGE $k_2$] at (3,0) {};
		\node at (3,0.4) {};
		\node at (3,.8) {}; 
		\node[label=below:\LARGE $k_3$] at (4,0) {};
		\node at (4,.4) {};
		\node[fill=White] at (5,0) {};
		\node[fill=White] at (6,0) {};
	\end{tikzpicture}
	\caption{An example configuration corresponding to the moment $\scp{\q}{u_{k_1}}^2\scp{\q}{u_{k_2}}^6\scp{\q}{u_{k_3}}^4$.}
\end{figure}
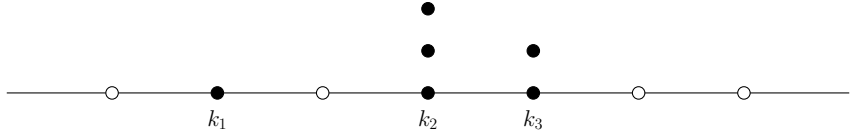

The time evolution of the observable $f_s(\bm{\xi})$ is given by the parabolic equation 
in the following theorem.
\bet[\cite{bourgade2013eigenvector}*{Theorem 3.1}]
Let $\q \in \S^{N-1}$ be a unit vector. For all $s\in (0,1)$, the moment observable $f_s$ defined in \eqref{eq:deffs} satisfies the equation
\beq\label{eq:emf}
\partial_s f_s(\bm{\xi})
=
\sum_{k\neq \ell}
2\xi_k(1+2\xi_\ell)
\frac{f_s(\bm{\xi}^{k,\ell})-f_s(\bm{\xi})}{N(\lambda_k(s)-\lambda_\ell(s))^2}.
\eeq
\eet 
\noindent  The equation \eqref{eq:emf} can be seen as a multi-particle random walk in a random environment given by the eigenvalues $\bm{\lambda}$. 
 We interpret it as saying that particles jump from site $k$ to site $\ell$ with rate $\frac{2\xi_k(1+2\xi_\ell)}{N(\lambda_k-\lambda_\ell)^2}$.

\subsection{Eigenvector moments of Gaussian divisible ensembles}
We now use the eigenvector moment flow to obtain a bound on large moments of eigenvectors of the dynamics \eqref{eq:dysonvect}. 
\noindent For the next proposition, we recall that $f_{\bm \lambda,s}(\bm{\xi})$ is a function of $\q \in \S^{N-1}$, although this is suppressed in the notation. 
\bep\label{p:highmom}
Fix $C, \delta>0$, $\theta\in(0,1/3)$, and a vector $\q\in\S^{N-1}$.  Then there exists $N_0 = N_0(C, \delta , \theta)$ such that for $N\geqslant N_0 $, $t \in[N^{-1/3+\theta},1]$,  and $0 \le n\leqslant \lceil (\log N)^{C} \rceil$, we have
\beq
\sup_{\bm{\xi}:\,\vert\bm{\xi}\vert=n}
\E\left[ f_{\bm \lambda,t}(\bm{\xi})\right]
\leqslant
(1+\delta)^n.
\eeq
\eep
\begin{rmk}
	\Cref{p:highmom} can be compared with \cite{bourgade2013eigenvector}*{Theorem 4.3}. However, in \cite{bourgade2013eigenvector} the authors obtain the asymptotic limit for finite moments, while we bound moments growing in $N$. 
\end{rmk}
\begin{proof}


Recall the notation of \Cref{l:goodset}, and let $\omega>0$ be a small parameter that will be determined later. Let $\mathcal A_1 (\q, \omega)$ be the set of paths $\bm \lambda = (\bm \lambda(s))_{0 \le s \le 1}$ such that the statements \eqref{e:isolaw}, \eqref{e:sclaw}, and \eqref{e:rigidity} hold with probability at least 
$ 1 - C_1^{1/2} \exp\left( -\frac{c_1}{2} (\log N)^{c_1 \log \log N}\right)$
with respect to the marginal distribution on paths $\bm \lambda$ induced by the measure $\P_\mu$ (defined in \Cref{l:goodset}). 
Then \Cref{l:goodset} and Fubini's theorem imply
\beq\label{e:Croothigh}
\P_\mu\left( \mathcal A_1 (\q, \omega) \right) \ge 1 - C_1^{1/2} \exp\left( -\frac{c_1}{2} (\log N)^{c_1 \log \log N}\right).
\eeq
We now restrict our attention to paths $\bm \lambda \in \mathcal A_1 (\q, \omega)$.

We define $\ximax = \ximax^{(n)}= \ximax^{(n)}(s)$ to be the maximizer among all configurations with $n$ particles for the moment observable \eqref{eq:deffs}, so that
\beq
f_s(\ximax) 
=
\sup_{\bm{\xi}:\,\vert\bm{\xi}\vert=n}
f_s(\bm{\xi}).
\eeq
When there are multiple maximizers, we pick one arbitrarily, subject to the constraint that $\widetilde{\bm \xi} (s)$ remains piecewise constant in $s$.  We let $(k_1,\dots,k_p)$ be the sites where $\ximax$ has at least one particle. Denoting the number of particles of $\widetilde {\bm \xi}$ at site $i$ by $\widetilde \xi_i$, these are the indices in $\unn{1}{N}$ such that $\widetilde \xi_{k_i}>0$. We emphasize that the implicit constants in each occurrence of the $\mathcal O$ notation in this proof will be independent of $p$ and all parameters in the statement of the proposition, including $\q$ and $n$.

Using \eqref{e:keyinequality}, we now proceed by induction on the number of particles $n$. 
Set $\eta =  N^{-2/3+\theta/2}$, and
consider the set of times $t_k = t  - \left(\lceil (\log N)^C \rceil - k \right) N^\omega \eta^{1/2}$ for $0 \le k \le  \lceil (\log N)^C \rceil $. Our induction hypothesis is that for every $\bm \lambda \in \mathcal A_1 (\q, \omega)$, we have
\beq\label{e:inducthypo}
\sup_{ s \in [ t_n , 1]}\sup_{\bm{\xi}:\,\vert\bm{\xi}\vert=n}
f_{\bm \lambda,t}(\bm{\xi})
\leqslant
(1+\delta)^n
\eeq
for $N > N_0$.
Here $N_0 = N_0(C, \delta, \theta)$ will be chosen in the course of the following computation.
For the base case $n=0$, \eqref{e:inducthypo} is trivial, since $f_s(\bm \xi)=1$. Next, for the induction step, fix $n \le \lceil (\log N)^C \rceil$ and suppose that \eqref{e:inducthypo} holds for $n-1$.

To complete the induction, we begin by deriving a differential inequality for $f_s(\ximax)=f_s(\ximax^{(n)}(s))$.\footnote{
At times when the maximum is obtained by two or more indices, $f_s(\ximax)$ may not be differentiable. But the following reasoning goes through with the redefinition $\partial_s f_s(\ximax) = \limsup_{u\rightarrow s} \frac{ f_s(\ximax) - f_u(\ximax)}{s -u}$.
}
 Fix some path $\bm \lambda \in \mathcal A_1 (\q, \omega)$. 
Using the eigenvector moment flow \eqref{eq:emf}, we see that for any $s \in [0,1]$, 
\begin{align}\label{eq:bounddyn0}
\partial_s f_s(\ximax)
&=
\sum_{i=1}^p \sum_{\ell\neq k_i}
2\widetilde\xi_{k_i}(1+2\widetilde\xi_{\ell})
\frac{f_s(\ximax^{k_i,\ell})-f_s(\ximax)}{N(\lambda_{k_i}(s)-\lambda_\ell(s))^2}\\
&\leqslant
\frac{2}{N\eta}
\sum_{i=1}^p \sum_{\ell\neq k_i}
(f_s(\ximax^{k_i,\ell})-f_s(\ximax))\frac{\eta}{(\lambda_{k_i}(s)-\lambda_\ell(s))^2+\eta^2}\\
&=\frac{2}{\eta} \sum_{i=1}^p
\Im
\sum_{\ell\neq k_i}
\frac{f_s(\ximax^{k_i,\ell})}{N(\lambda_\ell(s)-z_{k_i})} 
 - 
 \frac{2}{\eta}\sum_{i=1}^p
\Im
\sum_{\ell\neq k_i}
\frac{f_s(\ximax) }{N(\lambda_\ell(s)-z_{k_i})}.\label{e:secondtermdynamics}
\end{align}
In the last line, we defined $z_{k_i}=\lambda_{k_i}+\I\eta$. In the inequality, we used the fact that $f_s(\ximax^{k_i,\ell})\leqslant f_s(\ximax)$ by the definition of $\ximax$, so that $f_s(\ximax^{k,\ell})-f_s(\ximax) \le 0$ for all $k,\ell$. We also used that $\widetilde\xi_l \ge 0$ for all $l \in \llbracket 1 , N \rrbracket$, and $\widetilde\xi_{k_i}>0$.

We control the second term in \eqref{e:secondtermdynamics} by
\begin{align}
\sum_{i=1}^p
\Im
\sum_{\ell\neq k_i}
\frac{f_s(\ximax)}{N(\lambda_\ell(s) - z_{k_i})}
&=
f_s(\ximax)
\left(
\sum_{i=1}^p
\Im m_N^s(z_{k_i})
-\frac{p}{N\eta}\right)
\\
&=
f_s(\ximax)\left(
\sum_{i=1}^p
\Im \msc(z_{k_i})
+
\O{\frac{pN^\omega}{N\eta}}\right),\label{eq:bounddyn1}
\end{align}
where we used \eqref{e:sclaw} and $\eta \ge N^{-1+\omega}$ in the last equality. For the first term in \eqref{e:secondtermdynamics}, we use $f_s(\ximax) \ge  f_s(\ximax^{k,\ell}) \ge 0$, which holds by the definition of $f_s(\ximax)$, to see that
\beq\label{e:316}
\sum_{i=1}^p
\Im
\sum_{\ell\neq k_i}
\frac{f_s(\ximax^{k_i,\ell})}{N(\lambda_\ell(s)-z_{k_i})}
=
\sum_{i=1}^p
\Im
\sum_{\ell \neq k_1,\dots,k_p}
\frac{f_s(\ximax^{k_i,\ell})}{N(\lambda_\ell(s)-z_{k_i})}
+
\O{
f_s(\ximax)\frac{p^2}{N\eta}
}.
\eeq
Now, for the first term in \eqref{e:316}, fix an index $i\in\unn{1}{p}$, and observe
\begin{multline}
\Im 
\sum_{\ell\neq k_1,\dots,k_p}
\frac{f_s(\ximax^{k_i,\ell})}{N(\lambda_\ell(s)-z_{k_i})}\\
=
\frac{1}{\M\left(\ximax\setminus\{k_i\}\right)}
\E\left[
	\prod_{k}
	\scp{\q}{\sqrt{N}\u_k^s}^{2(\xi_k - \one_{ k = k_i})}
	\Im 
	\sum_{\ell\neq k_1,\dots,k_p}
	\frac{\scp{\q}{\sqrt{N}\u_\ell}^2}{N(\lambda_\ell(s)-z_{k_i})}
	\middle\vert
	\bm{\lambda}
\right],
\end{multline}
where we denote by $\ximax\setminus\{k_i\}$ the configuration $\ximax$ with a particle removed from the site $k_i$. We also used the identity $\M\left(\ximax^{k_i,\ell}\right)=\M\left(\ximax\setminus\{k_i\}\right)$, which holds because $\ell\notin \{k_1,\dots,k_p\}$. The last sum can be related to the resolvent $G^s(z_{k_i})$ using \eqref{eq:defsresolv}. Using this fact in conjunction with \eqref{e:isolaw}, we find
\beq\label{3144}
\Im
\sum_{\ell\neq k_1,\dots,k_p}
\frac{\scp{\q}{\u_\ell}^2}{(\lambda_\ell(s)-z_{k_i})}
\leqslant
\Im \scp{\q}{G^s(z_{k_i})\q}
=
\Im \msc(z_{k_i})
+
\O{
N^\omega\sqrt{\frac{\Im \msc(z_{k_i})}{N\eta}}+\frac{N^\omega}{N\eta}
}
\eeq
holds with exponentially high probability, by our assumption that $\bm \lambda \in \mathcal A_1 (\q, \omega)$.
We therefore obtain
\beq\label{eq:bounddyn2}
\Im \sum_{\ell\neq k_1,\dots,k_p}
\frac{f_s(\ximax^{k_i,\ell})}{N(\lambda_\ell(s)-z_{k_i})}
\leqslant
\left(
	\Im \msc(z_{k_i})
	+
	\O{N^\omega\sqrt{\frac{\Im \msc(z_{k_i})}{N\eta}}+\frac{N^\omega}{N\eta}}
\right)f_s(\ximax\setminus\{k_i\}) + \O{N^{-2}},
\eeq
where the $\O{ N^{-1}}$ term comes from using the trivial bounds 
$\scp{\q}{\u_\ell}^2 \le 1$ and
$\Im \frac{1}{(\lambda_\ell(s)-z_{k_i})} \le  \eta^{-1}$ on the exceptional event where
\eqref{e:isolaw} fails. We also used the fact that the probability of this event is at most $C_1^{1/2} \exp\left( -\frac{c_1}{2} (\log N)^{c_1 \log \log N}\right)$, by the definition of $\mathcal A_1$.

Observe that  $f_s(\ximax\setminus\{k_i\})\le f_s(\ximax^{(n-1)})$ by the definition of $\ximax^{(n-1)}$. Then inserting \eqref{eq:bounddyn1}, \eqref{e:316}, and \eqref{eq:bounddyn2} in \eqref{eq:bounddyn0} gives
\begin{align}
\partial_s f_s(\ximax^{(n)})
&\leqslant
-\frac{2f_s(\ximax^{(n)})}{\eta}\left(\sum_{i=1}^p\Im \msc(z_{k_i})+\O{\frac{p^2N^\omega}{N\eta}}\right)
\\
&+
\frac{2 f_s(\ximax^{(n-1)} )}{\eta} 
\left(
\sum_{i=1}^p\Im \msc(z_{k_i})+\O{N^\omega\sqrt{\frac{\Im \msc(z_{k_i})}{N\eta}}+\frac{N^\omega}{N\eta}}
\right)
+ \O{N^{-2}}.
\end{align}
Rearranging this yields
\begin{align}\label{e:rearrange}
\partial_s f_s(\ximax^{(n)}) \le& - \frac{2}{\eta}\left(\sum_{i=1}^p\Im \msc(z_{k_i})+\O{\frac{p^2N^\omega}{N\eta}}\right)\\  &\times \left(  f_s(\ximax^{(n)})  -  f_s(\ximax^{(n-1)} ) \frac{\sum_{i=1}^p\Im \msc(z_{k_i})+\O{N^\omega\sqrt{\frac{\Im \msc(z_{k_i})}{N\eta}}+\frac{N^\omega}{N\eta}}}{\sum_{i=1}^p\Im \msc(z_{k_i})+\O{\frac{p^2N^\omega}{N\eta}}} \right) + \O{N^{-2}} .
\end{align}

We now restrict our attention to times  $s \ge N^{-1/3+\theta}/2$, so that $s \gg \sqrt{\eta}$. Set $\omega = \theta/100$, and recall from \eqref{e:msclower} that $\Im \msc(z) \ge \scc \sqrt{\eta}$ for some $\scc >0$. Then

\beq\label{e:coefflimit}
\lim_{N\rightarrow \infty} \frac{\sum_{i=1}^p\Im \msc(z_{k_i})+\O{N^\omega\sqrt{\frac{\Im \msc(z_{k_i})}{N\eta}}+\frac{N^\omega}{N\eta}}}{\sum_{i=1}^p\Im \msc(z_{k_i})+\O{\frac{p^2 N^\omega}{N\eta}}}
 = 1,
\eeq
where we used our assumption that $1\le p\leqslant n\leqslant  \lceil (\log N)^C \rceil$. By \eqref{e:coefflimit} and \eqref{e:rearrange}, there exists $N_0 = N_0(C, \delta, \theta)$ such that, for $N \ge N_0$, 

\beq\label{e:rearrange2}
\partial_s f_s(\ximax) \le - \frac{2}{\eta}\left(\sum_{i=1}^p\Im \msc(z_{k_i})+\O{\frac{p^2 N^\omega}{N\eta}}\right)\left(  f_s(\ximax)  -  f_s(\ximax^{(n-1)} ) (1 + \delta/2) \right)
+\O{N^{-2}}.
\eeq
Note that $N_0$ does not depend on $n$ since we can bound the error $\frac{p^2 N^\omega}{N\eta}$ by $\frac{(\log N)^{3C} N^\omega}{N\eta}$. Using $\Im \msc(z) \ge \scc \sqrt{\eta}$ from \eqref{e:msclower}, $p \ge 1$, and the definition of $\omega$, we find after possibly adjusting $N_0$ upward that, for $N \ge N_0$, 
\beq\label{e:rearrange3}
\partial_s  f_s(\ximax)  \le - \frac{c}{2\sqrt{\eta}} \left(  f_s(\ximax)  -  f_s(\ximax^{(n-1)} ) (1 + \delta/2) \right) + \frac{ c }{2 N}
\eeq
for $c =  \scc /10 > 0$.

We now consider times $s \ge t_{n-1}$ and use the induction hypothesis $\eqref{e:inducthypo}$ to bound $f_s(\ximax^{(n-1)} )$ in \eqref{e:rearrange3}. This gives 
\beq\label{e:rearrange3b}
\partial_s  \left( f_s(\ximax) - ( 1 + \delta)^{n-1} (1 + \delta/2) \right)  \le - \frac{c}{\sqrt{\eta}} \left(  f_s(\ximax)  -  ( 1 + \delta)^{n-1} (1 + \delta/2) \right).
\eeq
If $f_r(\ximax)  \le  (1+\delta)^{n-1} (1 + \delta/2)$ for some $r$, then $f_s(\ximax)  \le  (1+\delta)^{n-1} (1 + \delta/2)$ holds for all $s > r$, since $f_s(\ximax)$ is decreasing by \eqref{eq:bounddyn0}.
We therefore assume that $f_s(\ximax)  \ge (1+\delta)^{n-1} (1 + \delta/2)$ for all $s \in [t_{n-1}, t_n]$; 
 otherwise, the induction hypothesis \eqref{e:inducthypo} holds for $n$, and the induction step is complete.
Because $f_s(\ximax)  -  (1+\delta)^{n-1} (1 + \delta/2) \ge 0$, \eqref{e:rearrange3b} implies
\beq\label{e:rearrange4}
\partial_s \log \left( f_s(\ximax)  -  ( 1 + \delta)^{n-1} (1 + \delta/2) \right) \le - \frac{c }{\sqrt{\eta}}.
\eeq

Integrating \eqref{e:rearrange4} on the interval $[t_{n-1}, t_n]$ gives 
\beq
\log \left( f_{t_n} (\ximax)  -   ( 1 + \delta)^{n-1} (1 + \delta/2) \right) \le  \log \left( f_{t_{n-1}} (\ximax)  - ( 1 + \delta)^{n-1} (1 + \delta/2)  \right) - \frac{c }{\sqrt{\eta}}(t_n -t_{n-1}).
\eeq
Exponentiating yields
\begin{align}
f_{t_n} (\ximax)  & \le ( 1 + \delta)^{n-1} (1 + \delta/2)  + \left( f_{t_{n-1}}(\ximax)  -  ( 1 + \delta)^{n-1} (1 + \delta/2) \right) \exp\left( - \frac{c }{\sqrt{\eta}}({t_n}-t_{n-1}) \right)\\
& \le ( 1 + \delta)^{n-1} (1 + \delta/2)   +  f_{t_{n-1}} (\ximax)   \exp\left( - \frac{c }{\sqrt{\eta}}(t_n- t_{n-1}) \right).
\end{align}
Observe that because we assumed $\bm \lambda \in \mathcal A_1 (\q, \omega)$, using delocalization \eqref{e:rigidity}, the definition \eqref{eq:deffs} of $f^{(n)}_s(\bm \xi)$, and 
$n \le \lceil (\log N)^C \rceil$, we have that
$f^{(n)}_s(\ximax) \le N^{\omega} f^{(n-1)}_s(\ximax) + N^{-1}$
for all $s\in[0,1]$ and $N \ge N_0$. The $N^{-1}$ comes from the exceptional set where delocalization does not hold, as in \eqref{eq:bounddyn2}.
 Therefore 
\begin{align}\label{e:keyinequality}
f_{t_n}(\ximax)   &\le  ( 1 + \delta)^{n-1} (1 + \delta/2)  + \left( f_{t_{n-1} }(\ximax^{(n-1)})   N^{\omega}   + N^{-1}
\right)\exp\left( - \frac{c }{\sqrt{\eta}}(t_n- t_{n-1} ) \right).
\end{align}

By putting the induction hypothesis \eqref{e:inducthypo} for $n-1$ into \eqref {e:keyinequality}, we obtain for $N > N_0$ that
\begin{align}\label{e:keyinequality2}
f_{t_n}(\ximax^{(n)}) 
&   \le (1+\delta)^{n-1} (1 + \delta/2)   + \left((1+\delta)^{n-1}    N^{\omega} + N^{-1} \right)  \exp\left( - \frac{c }{\sqrt{\eta}}( N^\omega \eta^{1/2} ) \right)\\
&   \le (1+\delta)^{n-1} (1 + \delta/2)   + (1+\delta)^{n-1}  (\delta/2) \le (1 + \delta)^n. \label{e:keyinequality3}
\end{align}
In the last line we used that $N^{\omega}  \exp\left( - \frac{c }{\sqrt{\eta}}( N^\omega \eta^{1/2} ) \right) < \delta/3$ for $N > N_0$, where we increased $N_0$ if necessary. 
By our previous observation that $f_{s}(\ximax^{(n)})$ is decreasing, we deduce from \eqref{e:keyinequality3} that for each $\bm \lambda$,
\beq f_{s}(\ximax^{(n)}) \le (1 + \delta)^n
\eeq
for all $s \ge t_{n}$. This completes the induction step.

We have shown that \eqref{e:inducthypo} holds for all $0 \le n \le \lceil (\log N)^C \rceil$, and therefore that
\beq\label{e:emfgood}
\sup_{\bm{\xi}:\,\vert\bm{\xi}\vert=n}
\E\left[\one_{\mathcal A_1} f_{\bm \lambda,t}(\bm{\xi})\right]
\leqslant
(1+\delta)^n.
\eeq
Using \eqref{e:Croothigh}, we also have the trivial bound
\begin{align}\label{e:emfbad}
 \E\left[\one_{\mathcal A^c_1} f_{\bm \lambda,t}(\bm{\xi})\right]& \le N^n C_1^{1/2} \exp\left( -\frac{c_1}{2} (\log N)^{c_2 \log \log N}\right)\\ &\le N^{(\log N)^C} C_1^{1/2} \exp\left( -\frac{c_1}{2} (\log N)^{c_2 \log \log N}\right)
 \le \delta
\end{align}
for $N \ge N_0(\delta)$ after possibly increasing $N_0$.  

Combining \eqref{e:emfgood} and \eqref{e:emfbad}, we obtain
\beq 
 \E\left[ f_{\bm \lambda,t}(\bm{\xi})\right] \le (1+\delta)^n + \delta  \le (1+\delta)^n + \delta (1+\delta)^{n-1} \le (1 + 2\delta)(1 + \delta)^{n-1} \le (1 + 2\delta)^n.
\eeq
The claim follows after redefining $\delta$. 
\end{proof}

Using \Cref{p:highmom}, we now show the optimal delocalization bound for eigenvectors of  generalized matrices perturbed by small Gaussian noise.\bec\label{c:delocdyson}
Fix $\theta\in(0, 2/3 )$. For $s\in[N^{-1/3+\theta},1]$ and any $\delta, \epsilon>0$, there exists $C = C (\theta, \delta, \epsilon)$ such that
\beq\label{e:gaussdiv2}
\sup_{k \in\unn{1}{N}}
\P\left(
		\Vert\u^s_k \Vert_\infty
	\geqslant
	\sqrt{\frac{(2+\varepsilon)\log N}{N}}
\right)
\leqslant C( \log N)^{1/2} N^{-\log\left(\frac{2+\varepsilon}{2(1+\delta)^2}\right)},
\eeq
\beq\label{e:gaussdiv1}
\sup_{\q\in\S^{N-1}}
\P\left(	\sup_{k\in\unn{1}{N}}
		\scp{\q}{\u_k^s}
	\geqslant
	\sqrt{\frac{(2+\varepsilon)\log N}{N}}
\right)
\leqslant C( \log N)^{1/2} N^{-\log\left(\frac{2+\varepsilon}{2(1+\delta)^2}\right)},
\eeq
and
\beq\label{e:gaussdiv3}
\P\left(	
	\sup_{k\in\unn{1}{N}}
	\Vert\u^s_k \Vert_\infty
	\geqslant
	\sqrt{\frac{(4+\varepsilon)\log N}{N}}
\right)
\leqslant C( \log N)^{1/2} N^{-2\log\left(\frac{4+\varepsilon}{4(1+\delta)^2}\right)}.
\eeq
\eec
\begin{proof}
We only prove \eqref{e:gaussdiv1}, since the proofs of \eqref{e:gaussdiv2} and \eqref{e:gaussdiv3} are similar. Fix $\q\in \S^{N-1}$ and set $J =  \lfloor \log N \rfloor$. By Markov's inequality applied with the $2J$-th moment and a union bound over all $k\in\unn{1}{N}$,
\beq\label{eq:markov}
\P\left(
	\sup_{k\in\unn{1}{N}}
	\scp{\q}{\u_k^s}
	\geqslant
	\sqrt{\frac{(2+\varepsilon)\log N}{N}}
\right)
\leqslant
\frac{N^{1+J}}{(2+\varepsilon)^{J}(\log N)^{J}}
\E\left[
	\scp{\q}{\u_1^s}^{2 J }
\right].
\eeq
By \Cref{l:goodset}, we have
\begin{align}
\E\left[
	\scp{\q}{\u_1^s}^{2J}
\right]
&\le
\frac{(1+\delta)^{2J}}{N^{J}}(2J-1)!!\label{e:gd3}
\end{align}
for $N \ge N_0(\delta, \theta)$. 
Stirling's formula shows that $(2J-1)!! \le C  (2 J)^{(2 J+1)/2} N^{-1}$ for some $C > 0$.
Inserting this bound into \eqref{e:gd3},  we obtain using \eqref{e:gd3} together with \eqref{eq:markov} that
\beq
\P\left(
	\sup_{k\in\unn{1}{N}}
	\scp{\q}{\u_k^s}
	\geqslant
	\sqrt{\frac{(2+\varepsilon)\log N}{N}}
\right)
\leqslant
C ( \log N)^{1/2} N^{-\log\left(\frac{2+\varepsilon}{2(1+\delta)^2}\right)},
\eeq
after adjusting $C = C(\theta, \delta, \eps) >0$.
Since all bounds are uniform in $\q$, we can finish by taking the supremum over all possible $\q\in\S^{N-1}$.
This completes the proof. 
\end{proof}
\section{Eigenvector regularization}\label{s:eigvectreg}

We now construct regularized versions of the eigenvectors of any $M \in \matn$. These are used in the next section to implement our comparison arguments extending the optimal delocalization bounds derived in the previous section to all generalized Wigner matrices. 

In \Cref{s:constructre}, we state some preliminary results and give the definition of the regularized eigenvector projections in \Cref{d:regvect}. Then, in \Cref{s:estimatesre}, we provide estimates relating the regularized eigenvector projections to the usual eigenvectors, and also control their derivatives with respect to any matrix entry.

\subsection{Construction of regularized eigenvector projections} \label{s:constructre}

 We begin by stating a result on the regularization of eigenvalues. For any $w \in [0,1]$, $M = \left( m_{ij} \right)_{1\le i,j \le N} \in \matn$, and indices $a,b \in \unn{1}{N}$, we let $\Theta^{(a,b)}_w M  \in \matn$ be defined as follows. Set $\Theta^{(a,b)}_w M$ to be the $N\times N$ matrix whose $(i,j)$ entry is equal to $m_{ij}$ if $(i,j) \notin \{ (a,b), (b,a) \}$. If $(i,j) \in \{ (a,b), (b,a) \}$, then set the $(i,j)$ entry equal to $w m_{a,b} = w m_{b,a}$. Further, for $k\in\unn{1}{N}$, we denote $\hat{k}=\min(k,N+1-k)$. We also set $\Theta^{(a,b)}_w  G  = (\Theta^{(a,b)}_w M - z)^{-1}$. 
 
 In the case of bulk eigenvalues, a version of the following proposition appeared as  \cite{landon2018comparison}*{Lemma 3.2}. The proof was based on an explicit construction of $\widetilde{\lambda}_i$ using the Helffer--Sj\"{o}strand formula. With additional technical effort, a version of the same construction can also be applied to the edge. We prove this proposition in \Cref{a:eigreg}.

\bep\label{p:regulareigval}
Fix $\delta, \eps >0$. For all $i \in \unn{1}{N}$, there exist functions $\widetilde{\lambda}_{i, \delta, \eps} \colon \matn \rightarrow \R$ such that the following holds, where we write $\tilde \lambda_i = \tilde \lambda_{i,\delta, \eps}$. For any generalized Wigner matrix $H$, there exists an event $\mathcal A = \mathcal A (\delta, \eps)$ such that for all  $j\in \unn{1}{5}$ and $a,b,c,d \in \unn{1}{N}$,
\begin{equation}\label{e:regev}
	\one_{\mathcal A} \left\vert	
		\widetilde{\lambda}_i(H)-\lambda_i(H)
	\right\vert	
	\le  \frac{N^{\eps-\delta}}{N^{2/3}\hat{i}^{1/3}}, \qquad 
	\sup_{0\leqslant w \leqslant 1}	\one_{\mathcal A}
	\left\vert
		\partial_{ab}^j \widetilde{\lambda}_i(\Theta^{(c,d)}_w H)
	\right\vert
	\le
	\frac{N^{j (\eps+\delta)}}{N^{2/3} \hat i^{1/3}},
\end{equation}
and 
\beq\label{e:regevpbound}
\P( \mathcal A^c) \le C_1 \exp\left({-c_1(\log N)^{c_1 \log\log N}}\right),
\eeq
for some constants $C_1(\delta, \eps), c_1(\delta, \eps)>0$.
Further, there exists $C >0$ such that
\beq\label{e:detbounds}
\sup_{0\leqslant w \leqslant 1}
\left\vert
	\partial^j_{ab} \widetilde{\lambda}_i(\Theta^{(c,d)}_wH)
\right\vert
\leqslant C N^{Cj}.
\eeq
\eep

We now use \Cref{p:regulareigval} to construct regularized eigenvector projections. 
Suppose $M \in \matn$ has eigenvalues $\bm{\lambda}=(\lambda_1,\dots\lambda_N)$, labeled in increasing order. Fix $\delta_1, \eps_1 >0$ and let the corresponding regularized eigenvalues $\tilde \lambda_i (M) =\tilde \lambda_{i, \delta_1, \eps_1}(M)$ be given by \Cref{p:regulareigval}.
For any $i \in \unn{1}{N}$ and $\delta_2 >0$, we define the intervals
\beq\label{eq:definte1}
I_{\delta_2} (x) =\left[
	x- \frac{ N^{-\delta_2}}{N^{2/3} \hat{i}^{1/3}},x + \frac{ N^{-\delta_2}}{N^{2/3} \hat{i}^{1/3}}\right], \qquad
\hat I_{\delta_2} (x) =\left[
	x - \frac{ N^{-\delta_2}}{2N^{2/3} \hat{i}^{1/3}}, x+ \frac{ N^{-\delta_2}}{2N^{2/3} \hat{i}^{1/3}}\right], 
\eeq
for $x\in \R$.
Given an interval $I$, we also denote the counting function for the eigenvalues $\bm{\lambda}$ by
\beq\label{e:lambdaI1}
\mathcal{N}_{\bm{\lambda}}(I)=
\left\vert
	\left\{
		i \in\unn{1}{N}\mid \lambda_i \in I
	\right\}
\right\vert.
\eeq

\bed\label{d:regvect}
Suppose $M \in \matn$ and fix $\delta_1, \eps_1, \delta_2, \eps_2 >0$. Let $\tilde \lambda_i = \tilde \lambda_{i, \delta_1, \eps_1}(M)$ denote the regularized eigenvalues of \Cref{p:regulareigval}. For $\ell\in\unn{1}{N}$ and $\q\in\S^{N-1}$, we define the regularized eigenvector projections of $M$ by
\beq\label{e:regvect}
v_\ell(\q) = v_\ell(M, \q) = v_\ell(M, \q, \delta_1, \eps_1, \delta_2,\varepsilon_2) = \frac{1}{\pi}
\int_{\hat I_{\delta_2}(\tilde{\lambda}_\ell)}\Im \scp{\q}{G(E+\I\eta_\ell)\q}\, \d E,
\qquad
\eta_\ell = \frac{N^{-\varepsilon_2}}{N^{2/3} \hat \ell^{1/3}}.
\eeq 
 Given $s\ge 0$, we let $v^s_\ell(\q)$ denote the regularized eigenvector projections for the Dyson Brownian motion $H_s$.
\eed
\begin{rmk}
In \Cref{s:main}, we will choose parameters in the previous definition such that $\delta_1 > \eps_1$ and $\delta_ 1 > \eps_2 > \delta_2$. 
The inequality $\delta_1 > \eps_1$ implies that the regularized eigenvalues approximate the standard ones at a scale smaller than the average inter-particle distance, by \eqref{e:regev}. Using \Cref{l:boundregeig} below, the inequalities $\delta_ 1 > \eps_2 > \delta_2$ are necessary to make $v_\ell$ approximate the eigenvector projections $\scp{\q}{\u_p}^2$.
\end{rmk}

We also introduce an event where all the estimates we need in order to control the regularized eigenvector projections hold. 

\bed\label{d:goodset}
Suppose $M \in \matn$. For $\omega,\delta_1, \varepsilon_1,\delta_2, \eps_2 > 0$, $\q \in\mathbb{S}^{N-1}$, and $k\in\N$, we define the events
\beq
\B(\q,\omega,\delta_1,\varepsilon_1,\delta_2, \eps_2)
=
\B_1(\q,\omega)\cap \B_1(\q,\eps_2/8) \cap \B_2(\delta_1,\varepsilon_1),
\eeq
\beq
\tilde \B(\q,\omega,\delta_1,\varepsilon_1,\delta_2,\eps_2, k, \ell)  = \B(\q,\omega,\delta_1,\varepsilon_1,\delta_2, \eps_2) \cap \B_3(\delta_2,k, \ell),
\eeq 
where
\begin{itemize}
\item $\B_1(\q,\omega)$ is the event defined in \Cref{l:goodset} where the isotropic local law \eqref{e:isolaw}, the local semicircle law \eqref{e:sclaw}, and rigidity and isotropic delocalization \eqref{e:rigidity} hold for all $\Theta^{(a,b)}_w M$ uniformly in $a,b \in \unn{1}{N}$ and $w \in [0,1]$,\footnote{More precisely, we demand that these equations hold with $\Theta^{(a,b)}_w G$ replacing $G^s$, $\lambda_k( \Theta^{(a,b)}_w M)$ replacing $\lambda_k$, the eigenvectors $\Theta^{(a,b)}_w \u_k$ of $\Theta^{(a,b)}_w M$ replacing $\u^s_k$, and $\sup_{a,b \in \unn{1}{N}}\, \sup_{w \in [0,1]}$ replacing $\sup_{s \in [ 0,1]}$.}


\item $\B_2(\delta_1,\varepsilon_1)$ is the event where 
\beq\label{e:regcloseB}
  \left\vert	
		\widetilde{\lambda}_i(M)-\lambda_i(M)
	\right\vert	
	\le  \frac{N^{\eps_1}}{N^{2/3+{\delta_1}}\hat{i}^{1/3}}
	\eeq
	and
\beq
\sup_{i \in \unn{1}{N}}\sup_{a,b,c,d \in \unn{1}{N} } \sup_{0\leqslant w \leqslant 1} 
	\left\vert
		\partial_{ab}^j \widetilde{\lambda}_i(\Theta^{(c,d)}_w M)
	\right\vert
	\le
	\frac{N^{j (\delta_1 + \eps_1)}}{N^{2/3} \hat i^{1/3}}
\eeq
hold for all $j\in \unn{1}{5}$, where the $\tilde \lambda_i$ are the regularized eigenvalues given by \Cref{p:regulareigval},
\item $\B_3(\delta_2, k, \ell)$ is the event on which
$\mathcal{N}_{\bm{\lambda}}(I_{\delta_2}(\lambda_\ell ))\leqslant k$,
where $I_{\delta_2}(\lambda_\ell )$ is defined in \eqref{e:lambdaI1} and $\bm \lambda$ is the vector of eigenvalues of $M$.
\end{itemize}
\eed 

\begin{figure}[!ht]
	\centering
	\includegraphics[width=.75\linewidth]{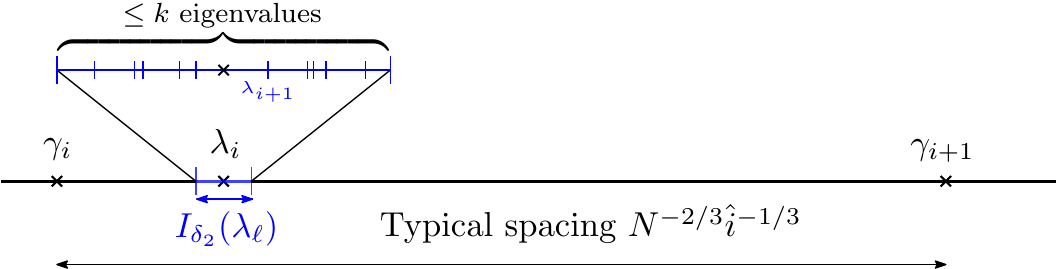}
	\caption{Illustration of the event $\B_3(\delta_2,k,\ell)$: on this event, we cannot fit more than $k$ eigenvalues on a sub-microscopic interval (which is asymptotically smaller than the typical spacing between eigenvalues).}
\end{figure}
\begin{rmk}
The presence of $\B_1(\q,\eps_2/8)$ in the above definition is a technical convenience. This set is used only in the proof of \Cref{l:belowG}, stated below.
\end{rmk}

The proofs of the following two lemmas are deferred to \Cref{a:preliminary}. In particular, the proof of \Cref{l:good2} is based on eigenvalue overcrowding estimates proved in \Cref{sec:lre}.

\bel\label{l:good1}
Let $H$ be a generalized Wigner matrix. For $\omega, \delta_1, \epsilon_1, \delta_2,\eps_2 >0$ and $\q \in\mathbb{S}^{N-1}$, there exist constants $C_1 = C_1(\omega, \delta_1, \eps_1, \delta_2, \eps_2), c_1 = c_1(\omega, \delta_1, \eps_1, \delta_2, \eps_2)>0$ such that
\beq\label{e:unionthis1}
\P \left(\B^c \left(\q,\omega,\delta_1,\varepsilon_1,\delta_2,\eps_2 \right) \right)  \le C_1 \exp\left({-c_1(\log N)^{c_1 \log\log N}}\right).
\eeq
\eel

\bel\label{l:good2}
Let $H$ be a generalized Wigner matrix and fix $\ell \in \unn{1}{N}$.
For any $D>0$ and $\omega, \delta_1, \epsilon_1, \delta_2 >0$ with $1/100 > \delta_2 > \eps_1 > 0$, 
there exists $k_0 = k_0(D, \delta_2)\in\mathbb{N}$ such that for any $k \ge k_0$,
\beq
\P \left(\tilde \B^c \left(\q,\omega,\delta_1,\varepsilon_1,\delta_2, \eps_2, k, \ell \right) \right)\le N^{-D}
\eeq
for all $\q\in\S^{N-1}$ and $N \ge N_0(D, \omega, \delta_1, \epsilon_1, \delta_2, \eps_2, k)$.
\eel


Before proceeding to our main estimates, we require the following lemma, which estimates Green's functions and their derivatives below the natural scale $N^{-1}$. The proof can be found in \Cref{a:preliminary}.

\bel\label{l:belowG}
Fix $\omega, \delta_1, \eps_1, \delta_2, \eps_2>0$ such that $\delta_1 > \eps_1$, $\q\in\S^{N-1}$, and $\ell \in \unn{1}{N}$.  
Let $H$ be a generalized Wigner matrix, and let $\B = \B(\q,\omega,\delta_1,\varepsilon_1,\delta_2, \eps_2)$ be the set from \Cref{d:goodset}. Let $\eta_\ell =\frac{N^{-\varepsilon_2}}{N^{2/3} \hat \ell^{1/3}}$.
Then there exists a constant $C= C(\eps_2)>0$ such that for all $E \in I_{\delta_2}(\lambda_\ell) \cup I_{\delta_2}(\tilde \lambda_\ell)$,
\beq\label{e:belowG}
\sup_{a,b \in \unn{1}{N}} \sup_{w \in [0,1]} \one_{\B}\Im\scp{\q}{\Theta^{(a,b)}_wG(E+\I\eta_\ell)\q} \le C N^{4\epsilon_2} \left( \frac{\hat \ell}{N}\right)^{1/3}
\eeq
and
\beq\label{e:extendedderiv}
\sup_{ a,b, c,d \in \unn{1}{N}} \sup_{j \in \unn{1}{5}} \sup_{w \in [0,1]} \one_{\B}\, \partial^j_{cd} \Im \langle \q , \Theta^{(a,b)}_w G(E + \I \eta_{\ell} ) \q \rangle \le C N^{10 j \eps_2}
\left( \frac{\hat \ell}{N}\right)^{1/3}.
\eeq
\eel

We also have the following lemma, whose proof is also deferred to \Cref{a:preliminary}.
\bel\label{l:eigsqsum}
Let $H$ be a generalized Wigner matrix. Fix $k\in\N$, $\q\in\S^{N-1}$, and $\omega,\delta_1,\varepsilon_1,\delta_2, \eps_2> 0$. Then there exists $C = C(\omega) >0$ such that, for for all $\ell\in\unn{1}{N}$, we have
\beq\label{e:eigsqsum}
\one_{\B} \sum_{p: \vert p-\ell\vert>N^{2\omega}}
\frac{1}{(\lambda_p-\lambda_\ell)^2}
\leqslant
C N^{4/3 + 2\omega} {\hat \ell}^{2/3},
\eeq
where $\B = \B(\q,\omega,\delta_1,\varepsilon_1,\delta_2,\eps_2)$. 
\eel

\subsection{Estimates on regularized eigenvector projections}\label{s:estimatesre}
In the next two lemmas, we show that the eigenvectors and their regularizations are comparable, in a certain sense, on a set of high probability. The following lemma shows that the eigenvectors can be used to bound their regularizations.

\bel\label{l:boundregeig}
Let $H$ be a generalized Wigner matrix. Fix $k\in\N$, $\q\in\S^{N-1}$, and $\omega,\delta_1,\varepsilon_1,\delta_2, \eps_2> 0$. Then 
for all $\ell\in\unn{1}{N}$, we have
\beq\label{e:boundedreclaim}
\one_{\tilde \B} v_\ell(\q )\leqslant \sum_{p=\ell-k + 1}^{\ell+k -1}\scp{\q}{\u_p}^2
+\mathcal{O}_{\omega, \eps_2}\left(
\frac{1}{N}\left(
		N^{3\omega + \delta_2 - \eps_2}
		+
		N^{\varepsilon_1+5\varepsilon_2 - \delta_1}
	\right)	
\right)
\eeq
where $\tilde \B = \tilde \B(\q,\omega,\delta_1,\varepsilon_1,\delta_2, \eps_2, k, \ell)$. 
\eel
\begin{proof}
We first give an estimate on a variant of the regularized eigenvector projections, where the resolvent is integrated on an interval centered on the actual eigenvalues instead of their regularized counterparts. Denote
\beq\label{eq:tilderegeig}
\vreg_\ell(\q) 
=
\frac{1}{\pi}\int_{\hat I_{\delta_2}(\lambda_\ell)}\Im\scp{\q}{G(E+\I\eta_\ell)\q}\, \d E.
\eeq
First,  we have
\beq\label{eq:regvect1}
\vreg_\ell(\q)
=
\frac{1}{\pi}\int_{\hat I_{\delta_2}(\lambda_\ell)}\sum_{p=1}^{N}\frac{\eta_\ell\scp{\q}{\u_p}^2}{(\lambda_p-E)^2+\eta_\ell^2}\d E
=
\frac{1}{\pi}\int_{\hat I_{\delta_2}(\lambda_\ell)}\left(\sum_{\vert p-\ell\vert \leqslant N^{2\omega}}+\sum_{\vert p-\ell\vert> N^{2\omega}}\right)
\frac{\eta_\ell \scp{\q}{\u_p}^2}{(\lambda_p-E)^2+\eta_\ell^2}\d E.
\eeq
In this decomposition, the second sum should be interpreted as an error term. We now work exclusively on the event $\tilde \B$ and drop this from our notation. On $\tilde \B$, we have the isotropic delocalization estimate $\scp{\q}{\u_p}^2\leqslant N^{-1+\omega}$, and thus
\begin{align}
\frac{1}{\pi}\sum_{\vert p-\ell\vert>N^{2\omega}}
\int_{\hat I_{\delta_2}(\lambda_\ell)}
\frac{\eta_\ell \scp{\q}{\u_p}^2}{(\lambda_p-E)^2+\eta_\ell^2}\d E
& \leqslant 
\sum_{\vert p-\ell\vert>N^{2\omega}}
\int_{\hat I_{\delta_2}(\lambda_\ell)} \frac{\eta_\ell N^{-1+\omega}}{(\lambda_p-E)^2}\d E
\\
&\leqslant \frac{4 \eta_\ell N^{-1+\omega - \delta_2}}{ N^{2/3} \hat \ell^{1/3} }
\sum_{\vert p-\ell\vert>N^{2\omega}}
\frac{1}{(\lambda_p-\lambda_\ell)^2}\label{e:insertpr1}
\end{align}
for large enough $N \ge N_0( \omega)$.
In the final inequality we used  the fact $|\lambda_p - E| >  \frac{1}{2} | \lambda_p - \lambda_\ell|$ for $N \ge N_0( \omega)$, which holds by the following calculation. By rigidity \eqref{e:rigidity} and the definition of $\hat I_{\delta_2}(\lambda_\ell)$, 
\beq
| \lambda_p - E |  \ge  | \lambda_p - \lambda_\ell | - |E - \lambda_\ell| \ge  | \lambda_p - \lambda_\ell |  - N^{- 2/3 - \delta_2} {\hat \ell}^{-1/3},
\eeq
so we just need to show that $| \lambda_p - \lambda_\ell | \gg N^{- 2/3 - \delta_2} {\hat \ell}^{-1/3}$. Observe
\beq
 | \lambda_p - \lambda_\ell |  \ge  | \gamma_p - \gamma_\ell | - |E - \lambda_\ell| \ge |\gamma_p - \gamma_\ell| - N^{-2/3} {\hat p}^{-1/3} - N^{-2/3} {\hat \ell}^{-1/3}.
\eeq
The case $\ell \in [ N/4, 3N/4]$ is trivial, so we consider $\ell \le N/4$; the case $\ell \ge 3N/4$ will follow by symmetry. 
When $p > \ell$, the claim now follows by using the definition \eqref{e:classical} to compute that $| \gamma_p - \gamma_\ell | \ge c  ( p - \ell) N^{-2/3}{\hat p}^{-1/3}$ for some $c>0$, which implies $|\gamma_p - \gamma_\ell| \gg N^{-2/3} {\hat p}^{-1/3} + N^{-2/3} {\hat \ell}^{-1/3}$ when $| p - \ell | > N^{2\omega}$.
The case $\ell > p$ is similar.

Now, inserting \eqref{e:eigsqsum} in \eqref{e:insertpr1} and using \eqref{eq:regvect1} shows that
\beq\label{eq:regvect3}
\vreg_\ell(\q)
\le
\frac{1}{\pi}\sum_{\vert p-\ell\vert\leqslant N^{2\omega}}
\int_{\hat I_{\delta_2}(\lambda_\ell)}
\frac{\eta_\ell \scp{\q}{\u_p}^2}{(\lambda_p-E)^2+\eta_\ell^2}\d E
+
C N^{-1 + 3\omega - \delta_2 - \eps_2} .
\eeq
We next remove the terms in the sum in \eqref{eq:regvect3} corresponding to the eigenvalues that do not lie in sub-microscopic interval $I_{\delta_2}(\lambda_\ell)$. We write
\beq
\sum_{\vert p-\ell\vert\leqslant N^{2\omega}}
\int_{\hat I_{\delta_2}(\lambda_\ell)}
\frac{\eta_\ell\scp{\q}{\u_p}^2}{(\lambda_p-E)^2+\eta_\ell^2}
\d E
=
\left(
\sum_{p:\lambda_p\in I_{\delta_2}(\lambda_\ell)}
+
\sum_{\substack{p:\lambda_p\not\in I_{\delta_2}(\lambda_\ell)\\\vert p-\ell\vert\leqslant N^{2\omega}}}
\right)
\int_{\hat I_{\delta_2}(\lambda_\ell)}
\frac{\eta_\ell\scp{\q}{\u_p}^2}{(\lambda_p-E)^2+\eta_\ell^2}
\d E.
\eeq
Using the isotropic delocalization, $| I_{\delta_2}(\lambda_\ell)  | =\eta_\ell N^{\eps_2 - \delta_2}$, and $|\lambda_p - E| \ge \frac{1}{2} | I_{\delta_2}(\lambda_\ell)  |$, we bound the second sum as follows:
\beq\label{eq:regvect4}
\sum_{\substack{p:\lambda_p\not\in I_{\delta_2}(\lambda_\ell)\\\vert p-\ell\vert\leqslant N^{2\omega}}}
\int_{\hat I_{\delta_2}(\lambda_\ell)}
\frac{\eta_\ell\scp{\q}{\u_p}^2}{(\lambda_p-E)^2+\eta_\ell^2}
\d E
\leqslant  N^{2\omega}  \eta_\ell N^{\eps_2 - \delta_2}   \frac{\eta_\ell}{ (\frac{1}{2} \eta_\ell N^{\eps_2 - \delta_2}  )^2 } N^{-1 + \omega} = 4 N^{-1 + 3\omega  + \delta_2 - \epsilon_2}.
\eeq
Putting the estimate \eqref{eq:regvect4} in \eqref{eq:regvect3},  we obtain that
\beq\label{eq:regvect5}
\vreg_\ell(\q)
= 
\frac{1}{\pi}\sum_{p:\lambda_p\in I_{\delta_2}(\lambda_\ell)}
\int_{\hat I_{\delta_2}(\lambda_\ell)}
\frac{\eta_\ell \scp{\q}{\u_p}^2}{(\lambda_p-E)^2+\eta_\ell^2}
+
\mathcal{O}_{\omega}\left(N^{-1 + 3  \omega  + \delta_2 - \eps_2}\right) .
\eeq

Observe that by \eqref{e:regcloseB},
\beq\label{e:symmdiffI}
 | \hat I_{\delta_2} (\lambda_\ell) \triangle \hat I_{\delta_2}( \tilde \lambda_\ell) | \le  | \lambda_\ell -\tilde \lambda_\ell | 
\le  \frac{N^{\eps_1}}{N^{2/3 + \delta_1} {\hat \ell}^{1/3}}. \eeq
Combining \eqref{e:symmdiffI} with \eqref{e:belowG} and using the definitions of $v_\ell(\q)$ and $\vreg_\ell(\q)$, we obtain on $\tilde \B$ that
\beq\label{eq:regvect6}
v_\ell(\q) = \vreg_\ell(\q)
+
\mathcal O_{\eps_2} \left( \frac{N^{\eps_1+ 5 \epsilon_2 - \delta_1}}{N}\right).
\eeq

By definition, on the event $\tilde \B$ there are at most $k$ eigenvalues in $I_{\delta_2}(\lambda_\ell)$. Using that the Poisson kernel integrates to 1,
\beq\label{e:lrlater}
\sum_{p:\lambda_p\in I_{\delta_2}(\lambda_\ell)}
\int_{\hat I_{\delta_2}(\lambda_\ell)}
\frac{\eta_\ell \scp{\q}{\u_p}^2}{(\lambda_p-E)^2+\eta_\ell^2}
\leqslant 
\sum_{p: \lambda_p\in I_{\delta_2}(\lambda_\ell)}
\scp{\q}{\u_p}^2
\leqslant \sum_{p=\ell -k + 1}^{\ell+k -1}\scp{\q}{\u_p}^2
\eeq
Finally, combining the estimates \eqref{eq:regvect5} and \eqref{eq:regvect6} with \eqref{e:lrlater}, we obtain \eqref{e:boundedreclaim}.
\end{proof}

The next lemma complements \Cref{l:boundregeig} and bounds the eigenvectors by their regularizations.

\bel\label{l:trivialbound}
Fix $\q\in\S^{N-1}$ and  $\omega,\delta_1,\varepsilon_1, \delta_2, \eps_2>0$ such that $\epsilon_2 > \delta_2$. Then for any $\ell \in\unn{1}{N}$, we have
\beq\label{e:onlyprove}
\one_{ \B}  \scp{\q}{\u_\ell}^2
\leqslant 
v_\ell(\q,\delta_2,\varepsilon_2)
+
\mathcal O_{\eps_2} \left(
	\frac{N^{\omega+\delta_2-\varepsilon_2}}{N}
	+
	\frac{N^{\eps_1+ 5 \epsilon_2 - \delta_1}}{N}
\right)
\eeq
where $ \B =  \B(\q,\omega,\delta_1,\varepsilon_1,\delta_2, \eps_2)$.
\eel
\begin{proof}
By directly integrating the Poisson kernel using $\int \frac{a \, \d x}{x^2 + a^2}  = \arctan\left( \frac{x}{a} \right)$ for $a > 0$, and delocalization, we have on $\B$ that
\beq
\scp{\q}{\u_\ell}^2
=
\frac{\eta_\ell}{\pi}\int_{\R}\frac{\scp{\q}{\u_\ell}^2}{(\lambda_\ell-E)^2+\eta_\ell^2}\d E
=
\frac{1}{\pi}
\int_{\hat I_{\delta_2}(\lambda_\ell)}
\frac{\eta_\ell\scp{\q}{\u_\ell}^2}{(\lambda_\ell-E)^2+\eta_\ell^2}\d E
+
\O{\frac{N^{\omega+\delta_2-\varepsilon_2}}{N}}.
\eeq
Because 
\beq
\frac{\eta_\ell \scp{\q}{\u_\ell}^2}{(\lambda_\ell-E)^2+\eta_\ell^2}\leqslant
\sum_{p=1}^N
\frac{\eta_\ell \scp{\q}{\u_p}^2}{(\lambda_p-E)^2+\eta_\ell^2}
=
\Im \scp{\q}{G(E+\I\eta_\ell)\q},
\eeq
we obtain
\beq\label{e:weobp}
\scp{\q}{\u_\ell}^2
\leqslant \tilde v_\ell(\q)
+
\mathcal O_{\eps_2} \left( \frac{N^{\omega+\delta_2-\varepsilon_2}}{N}\right),
\eeq
where $\tilde v_\ell(\q)$ was defined in \eqref{eq:tilderegeig}.
Finally, we obtain the claimed bound \eqref{e:onlyprove} by using \eqref{eq:regvect6} to replace $\tilde v_\ell(\q)$  with $v_\ell(\q)$ in \eqref{e:weobp}.
\end{proof}

Fix $\ell\in\unn{1}{N}$ and $\delta_2>0$.
We now regularize the indicator function for the set where a level repulsion estimate holds around the regularized eigenvalue $\tilde \lambda_\ell$.  Let $q_{\ell,{\delta_2}}$ be a smooth function such that $q_{\ell,{\delta_2}}(x)=1$ for $x\in  \left[
	\tilde \lambda_\ell - \frac{ 3N^{-\delta_2}}{2N^{2/3} \hat{\ell}^{1/3}},\, \tilde \lambda_\ell+ \frac{ 3N^{-\delta_2}}{2N^{2/3} \hat{\ell}^{1/3}}
\right]$, 
$q_{\ell,{\delta_2}}(x)=0$ for $x\notin \left[
	\tilde \lambda_\ell - \frac{ 2N^{-\delta_2}}{N^{2/3} \hat{\ell}^{1/3}},\, \tilde \lambda_\ell + \frac{ 2N^{-\delta_2}}{N^{2/3} \hat{\ell}^{1/3}}
\right]$, and there exists $C>0$ (independent of $\delta_2$ and $\ell$) such that both $\left|q^{(j)}_{\ell,{\delta_2}}(x) \right| \le CN^{(2/3 + {\delta_2})j } \hat \ell^{j/3}$  and $q_{\ell,{\delta_2}}(x)\leqslant 1$ hold for all $x\in\R$ and $j \in \unn{1}{5}$. We now fix $\nu>0$ and further define
\beq
F(M) = F_{\ell,{\delta_2}, \nu}(M)=\sum_{p: \vert p-\ell\vert \leqslant N^\nu} q_{\ell,{\delta_2}}\left(\tilde{\lambda}_p\left(M \right)\right)
\eeq
for any $M \in \matn$.

Fix $k\in\N$ and consider a smooth function $r=r_k$ such that $r(x)=1$ for $x\leqslant k - 1 $, $r(x)=0$ for $x\geqslant k$, and  $| r^{(j)}(k)(x) | \le C$ for all $x \in \R$ and $j \in \unn{1}{5}$, for some $C>0$. 
Given $M \in \matn$, we define
\beq\label{e:chidef}
\chi_k = \chi_k(M)= r_k(F_{\ell,{\delta}, \nu}(M)).
\eeq
The function $\chi_k(M)$ is a smooth version of the indicator function for the event that there are no more than $k$ eigenvalues of $M$ in the sub-microscopic interval $ I_{\delta_2}(\tilde \lambda_\ell)$. 

 For $M \in \matn$, $\ell \in \unn{1}{N}$, and $m \in \N$, we define 
\beq\label{e:Tdef}
T(M) = T(M, \ell) = N \chi_k(M, \ell)  v_\ell(M, \q), \qquad T_m(M) = T(M)^m.
\eeq 
We now control high moments of the observable $T(M)$ using \Cref{l:boundregeig} and \Cref{p:highmom}. 

\bel \label{l:dynamichighmom}
Fix $\delta_1 \in (0,1)$,  $C_1, \theta$, $\nu >0$, and $\q\in\S^{N-1}$. For any generalized Wigner matrix $H$, $s\in  [N^{-1/3 + \theta}, 1]$, $\ell\in\unn{1}{N}$,  positive integer $n\leqslant (\log N)^{C_1}$, and any $\delta,D>0$, there exist constants $C = C(\delta_1)>0$, $\sigma = \sigma(\delta_1)>0$ such that
\beq
\E\left[
T_n(H_s)
\right]
\leqslant 
(2k (1+\delta))^{n}
(2n)!!
\left(
	1
	+
	CN^{-\sigma}
	4^n
\right)
+N^{-D}
\eeq
for $N \ge N_0 ( C_1,D, \delta, \theta, \delta_1)$, where $T(H_s)$ is defined as in \eqref{e:Tdef} using 
$\eps_2 = \delta_1/10^2$, $ \delta_2 = \delta_1/10^3$,  $\eps_1 = \delta_1/10^4$,
$\nu = \delta_1/ 10^5$, and any choice of $k \in \N$.

\eel
\begin{proof}
Set $\omega =  \nu$ and $\tilde \B= \tilde \B(\q,\omega,\delta_1,\varepsilon_1,\delta_2,k, \ell)$. 
We write

\beq\label{e:hpdecomp}
\E\left[ T_n(M) \right]
=
\E\left[
	T_n(M)
	\left(
		\one_{\tilde \B}
		+
		\one_{\tilde \B^c}
	\right)
\right].
\eeq
By \Cref{l:boundregeig}, we have
\beq\label{eq:momentreg1}
\one_{\tilde \B} Nv_\ell(\q)
\leqslant \sum_{p=\ell-k+1}^{\ell+k-1} N\scp{\q}{\u_p^s}^2 + \mathcal{O}_{\delta_1} ({N^{-\sigma}}).
\eeq
for some $\sigma = \sigma(\delta_1) > 0$. This lemma applies because $(1+ s^2)^{-1/2} H_s$ is a generalized Wigner matrix, and its eigenvectors are invariant under rescaling. Because $0\le \chi_k \le 1$ and $v_\ell \ge 0$, we have that 
\beq
\E\left[
 T_n(M) 
	\one_{\tilde \B}
\right]
\le
\E\left[
	\left(
		Nv^s_\ell(\q)
	\right)^n
	\one_{\tilde \B}
\right].
\eeq
Using \eqref{eq:momentreg1}, we have
\beq\label{e:rcp1}
\E\left[
	T_n(M) \one_{\tilde \B}
\right]
\leqslant 
\E\left[
	\left(
		\sum_{p=\ell-k+1}^{\ell+k-1}
		N\scp{\q}{\u_p^s}^2
		+
		\mathcal{O}_{\delta_1} (N^{-\sigma})
	\right)^n
\right].
\eeq

By applying \Cref{p:highmom}, we have for $ m\le n\leqslant (\log N)^{C_1}$ that
\beq\label{e:insertdyn}
\E\left(
	\sum_{p=\ell-k+1}^{\ell+k-1}
	N\scp{\q}{\u_p^s}^2
\right)^m
=
\E\left[
	\sum_{p_1,\dots,p_m=\ell-k+1}^{\ell+k-1}
	\prod_{q=1}^m N\scp{\q}{u^s_{p_q}}^2
\right]
\leqslant 
(2k - 1)^{m}
(1+\delta)^m
(2m)!!
\eeq
for $N \ge N_0(C_1, \delta, \theta)$, since $s\in [ N^{-1/3 + \theta}]$. Inserting \eqref{e:insertdyn} into \eqref{e:rcp1} gives
\begin{align}
\E\left[
	\left(
		\sum_{p=\ell-k+1}^{\ell+k-1}
		N\scp{\q}{\u_p^s}^2
		+
		\O{N^{-\sigma}}
	\right)^n
\right]
&\leqslant 
(2k - 1)^n(1+\delta)^n(2n)!!\left(
	1
	+
	CN^{-\sigma}\sum_{i=1}^{2n}\binom{2n}{i}
\right)\\
&\leqslant 
(2k)^n(1+\delta)^n(2n)!!\left(
	1
	+
	CN^{-\sigma}2^{2n}
\right),
\end{align}
for some $C= C(\delta_1)>0$, where in the last inequality we used $\sum_{i=1}^{2n}\binom{2n}{i}\leqslant 2^{2n}$.

Let $\mathcal C$ be the event where rigidity \eqref{e:rigidity} holds for $H$ with parameter $\omega' = \nu/10$. We now claim that there exists $N_0(\delta_1) >0$ such that
\beq\label{e:1keystep}
\one_{\mathcal C\cap \B_2(\delta_1, \eps_1)} F(M) =\one_{\mathcal C \cap \B_2(\delta_1, \eps_1)} \sum_{p: \vert p-\ell\vert \leqslant N^\nu} q_{\ell,{\delta_2}}\left(\tilde{\lambda}_p\left(M \right)\right) = \one_{\mathcal C\cap \B_2(\delta_1, \eps_1)}\sum_{p \in \unn{1}{N}} q_{\ell,{\delta_2}}\left(\tilde{\lambda}_p\left(M \right)\right).
\eeq
for $N > N_0$.
(We recall $\B_2(\delta_1, \eps_1)$ was defined in \eqref{e:regcloseB}.)
It suffices to show that $q_{\ell,{\delta_2}}\left(\tilde{\lambda}_p\left(M \right)\right) =0$ for $| p - \ell | > N^\nu$.  Suppose first that $\ell < N/4$, and consider two cases: $p - \ell \ge N^\nu$ and $\ell - p \ge N^\nu$.

 In the first case, for $p - \ell \ge N^\nu$, for $p \le N/2$ we have
\begin{align}
\left|  \tilde{\lambda}_p - \tilde{\lambda}_\ell  \right| &\ge \left|  \gamma_p - \gamma_\ell  \right|  - \left|  \gamma_\ell - \lambda_\ell \right|  - \left|  \lambda_\ell - \tilde{\lambda}_\ell \right| - \left|  \gamma_p - \lambda_p \right|  - \left|  \lambda_p - \tilde{\lambda}_p \right| \\
& \ge c N^{ \nu  - \nu/3 - 2/3} \ell^{-1/3} - N^{ \nu/10 - 2/3} \ell^{-1/3} - N^{(\eps_1 - \delta_1)  - 2/3} \ell^{-1/3}\\ & -N^{ \nu/10 - 2/3} \ell^{-1/3} - N^{ (\eps_1 - \delta_1) - 2/3} \ell^{-1/3}\\
&\ge (c/2) N^{ 2\nu/3 - 2/3} \ell^{-1/3}\label{e:theend1}
\end{align}
on $\mathcal C \cap \mathcal B_2(\delta_1, \eps_1)$, for some constant $c>0$ and $N > N_0(\delta_1)$.  We used the definition \eqref{e:classical} and the fact that the gaps $\gamma_{p+1} - \gamma_{p}$ decrease as $p$ increases from $1$ to $N/2$, so that 
\beq
\left|  \gamma_p - \gamma_\ell  \right| \ge \left|  \gamma_{\ell + N^{\nu}} - \gamma_\ell  \right|  \ge N^{\nu} (\gamma_{\ell + N^\nu + 1} - \gamma_{\ell + N^\nu}) \ge c  N^{ \nu - 2/3} ( \ell + N^{\nu} ) ^{-1/3} \ge c   N^{  \nu  - \nu/3 - 2/3}  \ell^{-1/3},
\eeq 
where we decreased the value of $c$ in the last inequality.
We also used rigidity \eqref{e:rigidity}, \eqref{e:regcloseB}, $\delta_1 > \eps_1$, and $N^{ - 2/3} \ell^{-1/3} \ge N^{ - 2/3} p^{-1/3}$. We therefore obtain \eqref{e:1keystep} from \eqref{e:theend1} after recalling that $q_{\ell, \delta_2}$ is supported on $\left[
	\tilde \lambda_\ell - \frac{ 2N^{-\delta_2}}{N^{2/3} \hat{\ell}^{1/3}},\, \tilde \lambda_\ell + \frac{ 2N^{-\delta_2}}{N^{2/3} \hat{\ell}^{1/3}}
\right]$. When $p \ge N/2$, the argument is simpler, and we omit it. The second case, where $\ell - p \ge N^\nu$, is proved similarly.

We have established the claim for $\ell < N/4$. When $\ell \ge 3N /4 $, the conclusion follows by symmetry, and for $N /4 \le \ell \le 3N /4$ the argument is both similar and easier, so we omit it. We therefore see that \eqref{e:1keystep} holds for all choices of $\ell \in \unn{1}{N}$.

Let $\B_3 = \B_3 (\delta_2, k, \ell) = \{\mathcal{N}_{\bm{\lambda}}(I_{\delta_2}( \lambda_\ell ))\leqslant k\}$ be the event where there are at most $k$ eigenvalues in $ I_{\delta_2} (  \lambda_{\ell} )$, as defined in \Cref{d:goodset}.
Using \eqref{e:1keystep}  and the definition of $q_{\ell, \delta_2}$, we see that $\one_{\mathcal C\cap \B_2(\delta_1, \eps_1)} F(M) \le k$ implies there are at most $k$ regularized eigenvalues $\tilde \lambda_p$ in $\left[
	\tilde \lambda_\ell - \frac{ 3N^{-\delta_2}}{2N^{2/3} \hat{\ell}^{1/3}},\, \tilde \lambda_\ell + \frac{ 3N^{-\delta_2}}{2N^{2/3} \hat{\ell}^{1/3}}
\right]$. This in turn implies, using \eqref{e:regcloseB}, that there are most $k$ eigenvalues $\lambda_p$ in $I_{\delta_2}( \lambda_\ell)$ for $N \ge N_0(\delta_1)$. Therefore, we have
\beq \left\{ \one_{\mathcal C\cap \B_2(\delta_1, \eps_1)} F(M) \le k \right\} \subset \B_3\eeq 
for $N \ge N_0(\delta_1)$. Using the definition of $\chi_k$, we see $\chi_k \neq 0$ implies $F(M) \le k$, so this gives $ {\mathcal C\cap \B_2 } \cap {\left\{ \chi_k \neq 0  \right\}} \subset \B_3$. 
We now observe
\begin{align}
\left\{ \chi_k \neq 0  \right\} \cap \tilde \B^c &\subset \left\{ \chi_k \neq 0  \right\} \bigcap \left( \left( \left( \mathcal C \cap \B_2 \right) \cap \tilde \B^c\right)  \cup  \left( \left( \mathcal C \cap \B_2 \right)^c \cap \tilde \B^c\right) \right)\\
&\subset    \left( \left\{ \chi_k \neq 0  \right\} \cap  \left( \mathcal C \cap \B_2 \right) \cap \tilde \B^c\right)  \cup  \left( \left( \mathcal C \cap \B_2 \right)^c \cap \tilde \B^c\right) \\
&\subset  \left( \B_3 \cap \tilde \B^c\right)  \cup  \big( \left( \mathcal C \cap \B_2 \right)^c \big)\\
&\subset  \B^c  \cup  \mathcal C^c,
\end{align}
where in the last line we used $\B_3 \cap \tilde \B^c  \subset \B^c$ and $\B_2^c\subset \B^c$, which hold by the definitions of $\B$ and $\tilde\B$ in \Cref{d:goodset}.

Using \Cref{l:good1} and  \Cref{l:goodset} to bound $\P(\B^c)$ and $\P(\mathcal C^c)$, we find 
\beq\label{e:newbadset}
\P\left( \B^c \right) + \P\left( \mathcal C^c \right) \le C_2 \exp\left({-c_1(\log N)^{c_1\log\log N}}\right), 
\eeq
for some constants $C_2(\delta_1), c_1(\delta_1)>0$.  

Therefore, for the second term in \eqref{e:hpdecomp}, we have
\begin{align}
\E\left[
	T_n(M)
	\one_{\tilde \B^c}
\right]
& \le
\E\left[
	T_n(M)
	(\one_{\mathcal \B^c} + \one_{\mathcal C^c} )
	\right]
\\
& \le
N^{5 (\log N)^{C_1}}
C_2 \exp\left({-c_1(\log N)^{c_1\log\log N}}\right)
\leqslant N^{-D}
\end{align}
for $N \ge N_0(D, \delta_1)$. Here we used the fact that $v^s_\ell$ is always bounded by $N^{3 + \eps_2} \le N^4$  using the trivial bound $\vert G_{ij}(E+\I\eta)\vert\leqslant \eta^{-1}$ for all $i,j \in \unn{1}{N}$ in the definition \eqref{e:regvect}. 
\end{proof}

We conclude this section by bounding the derivatives of $T(M)$. 
\bel
Let $H$ be a generalized Wigner matrix. Fix $\eps_1, \delta_1, \eps_2, \delta_2, \nu \in (0,1)$ and $n \le 10 \log N$, and set $\omega = \nu$, $\B= \B(\q,\omega,\delta_1,\varepsilon_1,\delta_2)$.
Then there exists $C = C(\eps_2) >0$ such that for all $\ell \in \unn{1}{N}$ and $j \in \unn{1}{5}$, we have
\beq\label{e:Tderiv}
\sup_{a,b,c,d \in \unn{1}{N} } \sup_{0\leqslant w \leqslant 1}  \one_{\mathcal B} \left| \partial_{ab}^j T_n(\Theta^{(c,d)}_w H) \right| 
\le C N^{ 25 j  (\nu + \delta_1 + \eps_1 + \delta_2 + \eps_2)},  \eeq
and
\beq\label{e:Tderiv2}
 \sup_{a,b,c,d \in \unn{1}{N} } \sup_{0\leqslant w \leqslant 1}   \left| \partial_{ab}^j T_n(\Theta^{(c,d)}_w H) \right| \le C N^{Cj},
\eeq
where we recall $T(M)$ was defined in \eqref{e:Tdef}.
\eel
\bp 
We use the chain rule in the definition of $T(M)$, \eqref{e:Tdef}, to compute its derivatives. Since $T_n(M) = T(M)^n$, there are powers of $n$ coming from the derivatives of $x\mapsto x^n$. Because $n \le 10 \log N$ by hypothesis, these powers are logarithmic in size, and they can be absorbed into the powers of $N$ appearing in the rest of this proof to produce the final bound \eqref{e:Tderiv}. (The $20j$ from \eqref{eq:combinderiv} below becomes $25j$ in the final estimate \eqref{e:Tderiv} to account for these logarithmic factors.) Then, by the product rule and $T(M) = N \chi_k(M, \ell)  v_\ell(M, \q)$, it is enough to bound the derivatives of $\chi = \chi_k(\Theta^{(c,d)}_w H)$ and $v_\ell = v_\ell(\q,\Theta^{(c,d)}_w H)$ separately. 

We begin by bounding the derivatives of $v_\ell$, starting with $\partial_{ab} v_\ell$. 
We claim that
\beq\label{e:vderiv}
\one_{\mathcal B} \left| \partial^{j}_{ab} v_\ell \right| \le C N^{-1 + (j+1) (\delta_1 + \eps_1) + 10j\eps_2}, \qquad \left| \partial^{j}_{ab} v_\ell \right| \le C N^{Cj}
\eeq
for $j \in \unn{1}{5}$ and some $C= C(\eps_2)>0$. 
By the Leibniz integral rule, we have
\begin{align}\label{e:deriv1}\pi  \partial_{ab} v_\ell &= \Im \scp{\q}{G( \tilde \lambda_\ell 
 + N^{-{\delta_2} - 2/3} {\hat \ell}^{-1/3}+\I\eta_\ell)\q} \left( \partial_{ab}  \tilde \lambda_\ell \right) \\ 
\label{e:deriv2} & - \Im \scp{\q}{G(  \tilde \lambda_\ell  - N^{-{\delta_2} - 2/3} {\hat \ell}^{-1/3}+\I\eta_\ell)\q} \left( \partial_{ab}  \tilde \lambda_\ell \right)   \\
\label{e:deriv3} &+ \int_{\hat I_{\delta_2}(\tilde{\lambda}_\ell)}\Im \partial_{ab} \scp{\q}{G(E+\I\eta_\ell)\q}\, \d E.
\end{align}
Using the derivative estimate \eqref{e:regev} and Green's function bound \eqref{e:belowG}, we find that
\beq\label{e:1std1}
 \left| \eqref{e:deriv1} + \eqref{e:deriv2}  \right| \le C N^{-1 + \delta_1 + \eps_1 + 4\eps_2}.
\eeq
Using $| \hat I_\delta(\tilde \lambda_\ell) | \le  \frac{ N^{-{\delta_2}}}{N^{2/3} \hat{\ell}^{1/3}}$ and \eqref{e:extendedderiv}, we find
\beq\label{e:1std2}
 | \eqref{e:deriv3} |  \le C N^{10  \eps_2} \left( \frac{\hat \ell}{N}\right)^{1/3}   \frac{ N^{-\delta_2}}{N^{2/3} \hat{\ell}^{1/3}} \le C  N^{-1 + 10 \eps_2 }. 
\eeq
Together, \eqref{e:1std1} and \eqref{e:1std2} imply
\beq
\left| \partial_{ab} v_\ell \right| \le C N^{-1 + \delta_1 + \eps_1 + 10\eps_2}
\eeq
For the second derivative, we obtain
\begin{align}
\label{e:deriv1A}\pi  \partial^2_{ab} v_\ell &= \Im \scp{\q}{G( \tilde \lambda_\ell  + N^{-\delta_2 - 2/3} {\hat \ell}^{-1/3}+\I\eta_\ell)\q} \left( \partial^2_{ab}  \tilde \lambda_\ell \right) 
\\ &+ \partial_{ab}\Im \scp{\q}{G( \tilde \lambda_\ell  + N^{-\delta_2 - 2/3} {\hat \ell}^{-1/3}+\I\eta_\ell)\q} \left( \partial_{ab}  \tilde \lambda_\ell \right)\\ 
\label{e:deriv2A} & - \Im \scp{\q}{G(  \tilde \lambda_\ell  - N^{-\delta_2 - 2/3} {\hat \ell}^{-1/3}+\I\eta_\ell)\q} \left( \partial^2_{ab}  \tilde \lambda_\ell \right) 
\\ &  - \partial_{ab} \Im \scp{\q}{G(  \tilde \lambda_\ell  - N^{-\delta_2 - 2/3} {\hat \ell}^{-1/3}+\I\eta_\ell)\q} \left( \partial_{ab}  \tilde \lambda_\ell \right) \\
\label{e:deriv3A} &+ \int_{\hat I_{\delta_2}(\tilde{\lambda}_\ell)}\Im \partial^2_{ab} \scp{\q}{G(E+\I\eta_\ell)\q}\, \d E,
\end{align}
and we deduce using \eqref{e:regev} and \eqref{e:belowG} that 
\beq
\left| \partial^{2}_{ab} v_\ell \right| \le C N^{-1 + 2 (\delta_1 + \eps_1) + 20\eps_2}
\eeq
Continuing in this way, we obtain the first bound in \eqref{e:vderiv}. The second bound follows from the trivial estimate $|G_{ij}(E + \I \eta) | \le \eta^{-1}$. 

Using \eqref{e:belowG} and the definition of $v_\ell$, we also note the bounds
\beq\label{e:0thv}
\one_{\mathcal B} \left|  v_\ell \right| \le C N^{-1 + 4\eps_2},
 \qquad \left| v_\ell \right| \le N^2
\eeq
for some $C = C(\eps_2) >0$. 

We now bound the derivatives of $\chi_k$. Recall that $\chi_k(M) = r_k(F(M))$, and that $r_k$ and its first five derivatives are bounded in absolute value by a constant. It therefore suffices to bound the derivatives of $F(M)$. 

Using the derivative estimate \eqref{e:regev} and the definition of $q$, we have for $j\in \unn{1}{5}$ that
\begin{align}\label{e:derivativeF1}
\one_{\mathcal B} \left| \partial_{ab}^j F(M)  \right|
& \le   \sum_{p : | p - \ell|  \le N^{\nu} } \left| \partial_{ab}^j q\left( \tilde \lambda_p \right) \right|\\ 
&\le \sum_{m=1}^j C N^{\nu} \left( N^{(2/3 + \delta_2)m} {\hat p}^{m/3} \right) \left( N^{-m(2/3) + (j+1)(\eps_1 + \delta_1)} \hat{p }^{-m/3} \right)
 \le C N^{(j+1) (\nu + \delta_1 + \eps_1 + \delta_2)}.
\end{align}
Additionally, using \eqref{e:detbounds}, for $j\in \unn{1}{5}$ we have the trivial bound
\beq\label{e:derivativeFtrivial1}
\one_{\mathcal B} \left| \partial_{ab}^j F(M)  \right|
 \le   \sum_{p : | p - \ell|  \le N^{\nu} } \left| \partial_{ab}^j q\left( \tilde \lambda_p \right) \right| \le C N^{C j  }.
\eeq
This yields 
\beq\label{e:chideriv}
\one_{\mathcal B} \left| \partial_{ab}^j \chi_k \right| \le C N^{(j+1) (\nu + \delta_1 + \eps_1 + \delta_2)} , \qquad \left| \partial_{ab}^j \chi_k \right| \le C N^{j ( C + \nu)}.
\eeq
Combining \eqref{e:vderiv}, \eqref{e:0thv}, and \eqref{e:chideriv} yields
\beq \label{eq:combinderiv}
\one_{\mathcal B} \left| \partial_{ab}^j \left( \chi_k v_{\ell} \right) \right| 
\le C N^{-1 + 20 j  (\nu + \delta_1 + \eps_1 + \delta_2 + \eps_2)},  
\qquad  \left| \partial_{ab}^j \left( \chi_k v_{\ell} \right) \right| \le C N^{Cj}
\eeq
for $j\in \unn{1}{5}$ and some $C = C(\eps) > 0$. We used above the fact that $|\chi_k | \le 1$ by definition. This completes the proof.
\ep

\section{Proofs of main results} \label{s:main}

Given the proceeding sections, we are ready to prove our main results. \Cref{s:main1} contains the proof of \Cref{t:main1},  \Cref{s:main2} contains the proof of \Cref{t:main2}, and \Cref{s:main3} contains the proof of \Cref{t:main3}.

In this section only, given some $\delta_1 >0$ we fix the choice of parameters
\beq\label{e:paramchoice}
\eps_2 = \delta_1/10^2, \quad \delta_2 = \delta_1/10^3, \quad 
\epsilon_1 = \delta_1/10^4,
\quad \nu = \omega =\delta_1/ 10^5,
\eeq
which appeared in the statement of \Cref{l:dynamichighmom}. With these choices, $T(H)$ has two free parameters, $\delta_1$ and $k \in \N$.

\subsection{Proof of \texorpdfstring{\Cref{t:main1}}{Theorem 1.2}}\label{s:main1}

\bel \label{l:newcomparison}
Let $H$ be a generalized Wigner matrix. Then there exists $\delta_1>0$ such that the following holds. For all $k\in \N$, there exists $C = C(k) >0$ 
such that
\beq
 \E \left[ T_n \left( H \right)\right]  \le C^n (2n)!!
\eeq
for all $n \in \N$ satisfying $ n \le \lceil  \log N \rceil$, where the parameters in the definition of $T_n$ are chosen as in \eqref{e:paramchoice}
\eel
\bp
We first recall the dynamics $H_s$ defined in \eqref{e:dysondyn} for any initial matrix $H_0$, and set $s_1 = N^{-1/4}$. 
By the moment matching lemma \cite{erdos2017dynamical}*{Lemma 16.2},
 there exists a generalized Wigner matrix $H_0$ such that the matrix $R= (1 +  s_1^2 )^{-1/2}H_{s_1}$ satisfies $\E[h_{ij}^k] = \E[r_{ij}^k]$ for $k\in\unn{1}{3}$ and $\left| \E[h_{ij}^4 ]  -  \E[r_{ij}^4 ] \right| \le C N^{-2} s_1$ for some constant $C>0$ depending only on the constants used to verify that \Cref{d:wigner} holds for $H_0$. 
We observe that $R$ is also a generalized Wigner matrix according to definition \Cref{d:wigner}; in particular, it has subexponential entry distributions, and the scaling factor $(1 +  s_1^2 )^{-1/2}$ is chosen so that its variance matrix satisfies condition \eqref{e:stochasticvar}.

Fix any bijection
\beq
\phi \colon \{ ( i, j) : 1 \le i \le j \le N\} \rightarrow \unn{1}{\gamma_N},
\eeq
where $\gamma_N = N ( N + 1) /2$, and define the matrices $H^1, H^2, \dots , H^{\gamma_N}$ by 
\beq
h_{ij}^\gamma = 
\begin{cases}
h_{ij} & \text{if } \phi(i,j) \leq \gamma
\\
r_{ij} & \text{if } \phi(i,j) > \gamma
\end{cases}
\eeq
for $i \le j$.

Fix some $\gamma \in \unn{1}{\gamma_N}$ and consider the indices $(i,j)$ such that $\phi(i,j) = \gamma$. For any $m \ge 1$, we may Taylor expand $T_m \left(H^\gamma\right)$ in the $(i,j)$ entry, write $\partial = \partial_{ij}$, and find

\begin{align}
\label{e:taylora11} T_m \left( H^\gamma \right) - T_m \left(  \Theta^{(i,j)}_0 H^\gamma \right) &= \partial T_m \left(  \Theta^{(i,j)}_0 H^\gamma \right) h_{ij} + \frac{1}{2!} \partial^2 T_m \left(  \Theta^{(i,j)}_0 H^\gamma \right) h_{ij}^2 
+ \frac{1}{3!}\partial^3 T_m \left(  \Theta^{(i,j)}_0 H^\gamma \right) h_{ij}^3 \\
& \label{e:taylora21} + \frac{1}{4!}\partial^4 T_m \left(  \Theta^{(i,j)}_0 H^\gamma \right) h_{ij}^4 + \frac{1}{5!}\partial^5 T_m \left(  \Theta^{(i,j)}_{w_1(\gamma)} H^\gamma \right) h_{ij}^5,
\end{align}
where $w_1(\gamma) \in [0,1]$ is a random variable depending on $h_{ij}$. 

Subtracting the analogous expansion of $T_m \left( H^{\gamma -1} \right)$ in the $(i,j)$ entry and taking expectation, we find
\begin{align}
 \label{e:fourthorder1} \E \left[ T_m \left( H^{\gamma } \right)  \right]- \E \left[ T_m \left(   H^{\gamma -1} \right) \right] &= \frac{1}{4!} \E \left[ \partial^4 T_m \left(  \Theta^{(i,j)}_0 H^\gamma \right) h_{ij}^4 \right] - \frac{1}{4!}\E \left[ \partial^4 T_m \left(  \Theta^{(i,j)}_0 H^\gamma \right) r_{ij}^4 \right] \\
  \label{e:fifthorder1} &+ \frac{1}{5!} \E \left[ \partial^5 T_m \left(  \Theta^{(i,j)}_{w_1(\gamma)} H^\gamma \right) h_{ij}^5\right]
  - \frac{1}{5!} \E \left[ \partial^5 T_m \left(  \Theta^{(i,j)}_{w_2(\gamma)} H^\gamma \right) r_{ij}^5 \right],\end{align}
where we used that $\Theta^{(i,j)}_0 H^\gamma$ is independent from $h_{ij}$ and $r_{ij}$, and that $\E [ h^k_{ij} ] = \E [ r^k_{ij} ]$ for $k \in \unn{1}{3}$. Here $w_2(\gamma) \in [0,1]$ is a random variable depending on $r_{ij}$. 

By \Cref{l:dynamichighmom} applied with the parameters $\theta = 1/100$, $\delta = 1$, $D=1$, and  $C_1 = 2$, we have for for $m \in \N$ with $m \le \lceil 2 \log N \rceil$ that
\beq\label{e:TRlogN}
\E \left[ T_{ m } (R) \right] \le
L^m (2m)!!
\eeq
for some $L = L(k, \delta_1) > 0$ and all $N \ge N_0(\delta_1)$. We set $g(m) = L^m (2m)!!$.

We now use \eqref{e:TRlogN} and the expansion in \eqref{e:fourthorder1} and \eqref{e:fifthorder1} to show that $\E \left[ T_{ m } (H) \right] \le 3 g(m)$ for all $m \le \lceil 2 \log N \rceil$.
Our argument proceeds by induction, with the induction hypothesis at step $m \in \N$ being that 
\beq\label{e:inductionhypo1}
 \E \left[ T_n \left(  \Theta^{(a,b)}_{w } H^\gamma \right)\right]  \le 3 g(n)
\eeq
holds for all $0 \le n \le m \le \lceil 2 \log N \rceil$ and choices of  $w \in [ 0 ,1 ]$ and $(a,b) \in \unn{1}{N}^2$.

 The base case $m = 0$ is trivial. Assuming the induction hypothesis holds for $m -1$, we will show it holds for $m$. Using the independence of $h_{ij}$ and $r_{ij}$ from $ \Theta^{(i,j)}_0 H^\gamma$, we may rewrite the first two terms terms on the right side of \eqref{e:fourthorder1} as 
\beq\
\E \left[ \partial^4 T_m \left(  \Theta^{(i,j)}_0 H^\gamma \right) h_{ij}^4 \right] - \E \left[ \partial^4 T_m \left(  \Theta^{(i,j)}_0 H^\gamma \right) r_{ij}^4 \right] =  \E \left[ \partial^4 T_m \left(  \Theta^{(i,j)}_0 H^\gamma \right)  \right]  \E \left[  h_{ij}^4 - r_{ij}^4  \right].
\eeq
For the second factor, we recall that $\left|  \E \left[ h_{ij}^4\right] - \E \left[ r_{ij}^4 \right]  \right|  \le C N^{-2}s_1 =  C N^{-2 - 1/4}$. For the first, we compute 
\begin{align}\label{e:4thcompute1}
\partial^4 T_m  = \partial^4 \left( T^m \right) &= m T_{m-1} T^{(4)}  + 3 m (m - 1 ) T_{m-2} (T^{(2)})^2 + m (m-1) (m-2) (m-3) T_{m-4} (T')^4 \\
& + 4 m (m-1) T_{m-2} T^{(1)} T^{(3)} + 6 m (m-1) (m-2) T_{m-3} (T')^2 T^{(2)}.
\end{align}
We write $\B= \B(\q,\omega,\delta_1,\varepsilon_1,\delta_2)$. Using the induction hypothesis \eqref{e:inductionhypo1} for $ n \le m -1$, $m \le \lceil 2 \log N \rceil$, the fact that $T_m \ge 0$, and \eqref{e:Tderiv}, we find that
\beq\label{e:4analogy1}
\left|  \E \left[ \one_{\mathcal B} \partial^4 T_m \left(  \Theta^{(i,j)}_0 H^\gamma \right)  \right]  \right| \le C ( \log N)^4 N^{100(\nu + \delta_1 + \eps_1 + \delta_2 + \eps_2)} g(m-1),
\eeq
where $C =  C(\delta_1) >0$ is a constant.

 Further, by \eqref{e:Tderiv2} and \Cref{l:good1}, we find 
\beq\label{e:4analogytrivial1}
\left|\E \left[  \one_{\mathcal B^c}  \partial^4 T_m \left(  \Theta^{(i,j)}_0 H^\gamma \right)  \right]  \right| \le C N^{-2}
\eeq
for some constant $C = C(\delta_1)  > 0$.\footnote{We note that the constants in the probability bound given by \Cref{l:good1} do not depend on the choice of $\gamma$, since the $H^\gamma$ verify \Cref{d:wigner} simultaneously for the appropriate choice of constants. Therefore, the $C$ in \eqref{e:4analogytrivial1} is uniform in $\gamma$. Analogous remarks apply to our other uses of the four moment method in this work.}

 It follows from \eqref{e:4analogy1} and \eqref{e:4analogytrivial1} that if $\delta_1$ is chosen small enough, so that\footnote{This bound is slightly stronger than necessary for the fourth order term, but will be needed for the fifth order remainder term.}
\beq\label{e:setparams}
130(\nu + \delta_1 + \eps_1 + \delta_2 + \eps_2)<  1/8,
\eeq
 then 
\beq\label{e:4final1}
\left| \E \left[ \partial^4 T_m \left(  \Theta^{(i,j)}_0 H^\gamma \right) h_{ij}^4 \right] - \E \left[ \partial^4 T_m \left(  \Theta^{(i,j)}_0 H^\gamma \right) r_{ij}^4 \right] \right| \le  C ( \log N)^4 N^{-2-1/8} g(m-1).
\eeq
holds for all $m \le \log N$, for some constant $C = C(\delta_1)  > 0$; we now fix $\delta_1$ so that it satisfies this condition. Therefore, if $ N \ge N_0 (\delta_1)$, then
\beq\label{e:4final2}
\left| \E \left[ \partial^4 T_m \left(  \Theta^{(i,j)}_0 H^\gamma \right) h_{ij}^4 \right] - \E \left[ \partial^4 T_m \left(  \Theta^{(i,j)}_0 H^\gamma \right) r_{ij}^4 \right] \right| \le  \frac{1}{2} N^{-2 } g(m-1).
\eeq

Let $\mathcal D$ be the event where $\sup_{i,j} | r_{ij} | + | h_{ij} |\le C N^{-1/2 + \delta_1}$ holds. Since the variables $r_{ij}$ and $h_{ij}$ are subexponential, we have
\beq\label{e:Csubexp}
\P\left( \mathcal D^c \right) \le D_1 \exp\left({-d_1(\log N)^{d_1\log\log N}}\right), 
\eeq
for some constants $D_1(\delta_1), d_1(\delta_1)>0$.

For the terms in \eqref{e:fifthorder1}, we compute
\begin{align}
\left| \E \left[ \partial^5 T_m \left(  \Theta^{(i,j)}_{w_1(\gamma)} H^\gamma \right) h_{ij}^5\right] \right| 
&\le \left| \E  \left[  \one_{\mathcal D} \partial^5 T_m \left(  \Theta^{(i,j)}_{w_1(\gamma)} H^\gamma \right) h_{ij}^5\right] \right| + \left| \E \left[ \one_{\mathcal D^c} \partial^5 T_m \left(  \Theta^{(i,j)}_{w_1(\gamma)} H^\gamma \right) h_{ij}^5\right] \right|\\
&\le C N^{-5/2 + 5 \delta_1} \left( \E \left[ \left| \partial^5 T_m \left(  \Theta^{(i,j)}_{w_1(\gamma)} H^\gamma \right) \right| \right]  + 1\right),
\end{align}
for some $C(\delta_1) > 0$. In the last line, we used \eqref{e:Csubexp} and the inequality of \eqref{e:Tderiv2}.

An analogous bound holds for the second term in \eqref{e:fifthorder1}. Then repeating the previous argument for the fourth order term given in \eqref{e:4analogy1} and \eqref{e:4analogytrivial1}, and using \eqref{e:setparams}, we find that there exists $C(\delta_1)>0$ such that
\begin{align}\label{e:5finala}
\Bigg|
\E \left[ \partial^5 T_m \left(  \Theta^{(i,j)}_{w_1(\gamma)} H^\gamma \right) h_{ij}^5\right]
  -   \E &\left[ \partial^5 T_m \left(  \Theta^{(i,j)}_{w_2(\gamma)} H^\gamma \right) r_{ij}^5 \right]
  \Bigg|\\
&\le  C (\log N)^5 N^{-5/2 + 5\delta_1 + 125(\nu + \delta_1 + \eps_1 + \delta_2 + \eps_2)} g(m-1)\\
& \le C (\log N)^5 N^{-2 - 1/8} g(m-1).
\end{align}
Therefore, if $ N \ge N_0 (\delta_1)$, then
\beq\label{e:5final2a}
\left|
\E \left[ \partial^5 T_m \left(  \Theta^{(i,j)}_{w_1(\gamma)} H^\gamma \right) h_{ij}^5\right]
-   \E \left[ \partial^5 T_m \left(  \Theta^{(i,j)}_{w_2(\gamma)} H^\gamma \right) r_{ij}^5 \right]
\right|
\le \frac{1}{2} N^{-2} g(m-1)
\eeq

Combining \eqref{e:4final2} and \eqref{e:5final2a} yields
\beq\label{e:gammainc}
\left| \E \left[ T_m \left( H^{\gamma } \right)  \right]- \E \left[ T_m \left(   H^{\gamma -1} \right) \right] \right| \le N^{-2}g(m-1) \le N^{-2} g(m),
\eeq
and summing \eqref{e:gammainc} over all $\gamma_N$ pairs $(i,j)$, we find
\beq\label{e:4momentconclude1}
\left| \E \left[ T_m \left( R \right)  \right]- \E \left[ T_m \left( H^\gamma \right) \right] \right| \le g(m)
\eeq
for any $\gamma$. Using \eqref{e:TRlogN} and \eqref{e:4momentconclude1}, we deduce that 
\beq
\E\left[ T_m\left( H^\gamma \right) \right] \le 2g(m).
\eeq
 This verifies the induction hypothesis \eqref{e:inductionhypo1} when $w = 1$. 

To address other values of $w$, we consider the following expansion:
\begin{align}
\label{e:taylorc1a} T_m \left( H^\gamma \right) - T_m \left(  \Theta^{(a,b)}_w H^\gamma \right) &= \partial T_m \left(  \Theta^{(a,b)}_0 H^\gamma \right) h_{ij} + \frac{1}{2!} \partial^2 T_m \left(  \Theta^{(a,b)}_0 H^\gamma \right) h_{ij}^2 
+ \frac{1}{3!}\partial^3 T_m \left(  \Theta^{(a,b)}_0 H^\gamma \right) h_{ij}^3 \\
& \label{e:taylorc2a} + \frac{1}{4!}\partial^4 T_m \left(  \Theta^{(a,b)}_0 H^\gamma \right) h_{ij}^4 + \frac{1}{5!}\partial^5 T_m \left(  \Theta^{(a,b)}_{\tau(w)} H^\gamma \right) h_{ij}^5,
\end{align}
Here $\tau(w) \in [ 0, 1]$ is a random variable. The same argument that gave the bound \eqref{e:4final2} shows that the right side of \eqref{e:taylorc1a} and the first term of \eqref{e:taylorc2a} may be bounded in absolute value by $\frac{1}{2} g(m)$. The second term of \eqref{e:taylorc2a} is also bounded in absolute value by $\frac{1}{2}g(m)$ by the reasoning leading to \eqref{e:5finala}. We  conclude
\beq
\sup_{w \in [0,1]} \sup_{a,b \in \unn{1}{n}} \E \left[ T_m \left(  \Theta^{(a,b)}_{w } H^\gamma \right)\right]  \le 3 g(m).
\eeq
This completes the induction and therefore concludes the proof.
\ep

\bel\label{l:markovclaim}
Let $H$ be a generalized Wigner matrix. Then there exists $\delta_1>0$ such that the following holds. For every $D >0$, $k\in \N$, and $\ell \in \unn{1}{N}$, there exists $C = C(D, k) >0$ such that
\beq\label{e:markovclaim}
\P\left( T(H, \ell) \ge C \log N  \right) \le C N^{-D}.
\eeq
\eel
\bp
 Set $J = \lceil \log N \rceil$. By Markov's inequality applied to the $J$-th moment of $T(H)$ and \eqref{l:newcomparison}, we have for any $K>0$ that
\beq\label{e:markovcompare}
\P\left( T(H) \ge 2K \log N \right) \le \P\left( T(H) \ge K J \right) \le \frac{\E \left[ T_{J} (H)  \right]}{(K J)^{ J} } \le \frac{  L^{J} (2 J-1)!! }{(K J )^{J} },
\eeq
where $L = L(k)$ and the $\delta_1$ defining $T(H)$ are given by \Cref{l:newcomparison}. By Stirling's formula,  
\beq (\lceil 2\log N -1\rceil)!! \le B  (2 \lceil  \log N \rceil)^{(2\lceil \log N\rceil +1)/2} N^{-1}\eeq for some constant $B > 0$. We then deduce from \eqref{e:markovcompare} that
\beq
\P\left( T(H) \ge KJ \right) \le \frac{ B L^{ J } (2J)^{(2 J+1)/2}}{(K J)^{J} N } =  \left( \frac{L}{K}\right)^{J } \frac{B \sqrt{J}}{N}.
\eeq
The conclusion follows after taking $C=\max(2K,B)$  and $ K \ge L \exp ( 2D ) $.
\ep

Before proceeding to the proof of \Cref{t:main1}, we require the following level repulsion estimate, which is proved in \Cref{sec:lre}.

\bep\label{p:lrmain}
 Let $H$ be a generalized Wigner matrix, and let $\bm{\lambda}$ be the vector of its eigenvalues. Fix $D>0$, $\fa>0$, and $\delta_2 \in (0, 1/100)$. Then there exist $k = k(D, \fa, \delta_2)\in \N$ and a constant $C = C(D, \fa , \delta_2) > 0$ such that for any $\ell \in \unn{1}{N} $,

	\beq
	\P \left(
	 \mathcal{N}_{\bm{\lambda}  }(I(\delta_2, \fa, \lambda_\ell)  )\geqslant k
	\right)
	\leqslant  C N^{- D},
	\eeq
	where  
	\beq I(\delta, \fa, \lambda_\ell) = \left[
	\tilde \lambda_\ell - \frac{ \fa N^{-\delta_2}}{N^{2/3} \hat{\ell}^{1/3}},\, \tilde \lambda_\ell+ \frac{ \fa N^{-\delta_2}}{N^{2/3} \hat{\ell}^{1/3}} \right].
	\eeq

\eep

\bp[Proof of \Cref{t:main1}]  Let $\mathcal C$ be the event where rigidity \eqref{e:rigidity} holds for $H$ with parameter $\omega' = \nu/10$. 
It was shown in the proof of \Cref{l:dynamichighmom} as \eqref{e:1keystep} that
\beq\label{e:1keystep2}
\one_{\mathcal C\cap \B_2(\delta_1, \eps_1)} F(M) = \one_{\mathcal C\cap \B_2(\delta_1, \eps_1)}\sum_{p \in \unn{1}{N}} q_{\ell,{\delta_2}}\left(\tilde{\lambda}_p\left(M \right)\right).
\eeq
For $k\in \N$, define
\beq I = \left[
	\tilde \lambda_\ell - \frac{ 3N^{-\delta_2}}{N^{2/3} \hat{\ell}^{1/3}},\, \tilde \lambda_\ell+ \frac{ 3N^{-\delta_2}}{N^{2/3} \hat{\ell}^{1/3}} \right], \qquad
	\mathcal A =\mathcal A(k) = \left\{ { \mathcal N}_{\bm \lambda}( I) \le k - 1 \right\},
	\eeq
where $\bm \lambda$ is the vector of eigenvalues of $H$. Then by the definition of $\chi_k$, we have for any $k \in \N$ that
\beq\label{e:chiequal1}
\A(k) \cap  \B_2(\delta_1, \eps_1) \cap \mathcal C \subset \{ F(M) \le k -1  \} \subset \{ \chi_k = 1 \}.
\eeq

Define $T(H)$ using the choice of $\delta_1 >0$ provided by \Cref{l:markovclaim}, the parameter choices \eqref{e:paramchoice}, and a parameter $k\in \N$, which will be chosen later. We  obtain using \eqref{e:chiequal1} and the definition $T(M) = N \chi_k(M) v_\ell(M,\q)$ given in \eqref{e:Tdef} that
$\one_{\A \cap \B_2 \cap \mathcal C} N v_\ell(\q) \le T(M)$. Therefore,  \eqref{e:markovclaim} yields
\beq\label{e:input1}
\P\left(  \one_{\A \cap \B_2 \cap \mathcal C} N  v_\ell(\q) \ge K \log N \right) \le K N^{-D - 2}
\eeq
for some $K = K(D,k) > 0$.
By \Cref{p:lrmain}, \Cref{p:regulareigval}, and \eqref{e:rigidity}, there exists $k_0 = k_0(D, \delta_1)$ such that
\beq\label{e:input2}
\P\left( \left(\A \cap \B_2 \cap \mathcal C \right)^c \right) \le N^{-D -2}
\eeq
for $k \ge k_0$ and $N \ge N_0(k, \delta_1)$.
We conclude from \eqref{e:input1} and \eqref{e:input2} that 
\beq\label{e:markov1}
\P\left(  N  v_\ell(\q) \ge K \log N \right) \le  (K+1) N^{-D - 2} 
\eeq
for $k \ge k_0$ and $N \ge N_0$. Set $\mathcal B = \B(\q,\omega,\delta_1,\varepsilon_1,\delta_2, \eps_2 )$. By \Cref{l:trivialbound}, we have
\beq\label{e:trivb}
\one_{ \B} \scp{\q}{\u_\ell}^2
\leqslant 
v_\ell(\q)
+
\O{
	N^{-1 -c }
}
\eeq
for some $c = c(\delta_1) >0$. Then combining \eqref{e:trivb} with \eqref{e:markov1} gives
\beq\label{e:markov12}
\P\left( \one_{\B}N  \scp{\q}{\u_\ell}^2 \ge  2K \log N \right) \le  (K+1) N^{-D - 2} 
\eeq
for $N\ge N_0 (k, \delta_1)$, after increasing $N_0$. By \Cref{l:good2}, there exists $k_1 \in \N$ such that for $k = \max (k_0, k_1)$, we have
\beq\label{e:markov2}
\P\left( N \scp{\q}{\u_\ell}^2 \ge  2K \log N \right) \le  N^{-D -1 }
\eeq
for $N \ge N_0(D)$. This proves \eqref{e:isotropicmain} after taking taking $C=C(D)$ large enough. (Recall that $\delta_1$ was fixed earlier). Finally, \eqref{e:isotropicmain2} follows from \eqref{e:markov2} after taking $\q = \mathbf e_i$ for $i \in \unn{1}{N}$ and applying a union bound over all $N$ choices of $i$, where the $\mathbf e_i$ are the standard basis vectors, and increasing the value of $C$.
\ep

\subsection{Proof of \texorpdfstring{\Cref{t:main2}}{Theorem 1.3}}\label{s:main2}

Given a vector $\bm w \in \R^N$ and $\beta >0$, we introduce
\beq
A  = A _{\beta}(\bm w) = \frac{1}{\beta} \log \left( \sum_{i \in \unn{1}{N}} \exp \left( \beta  w_i \right)\right).
\eeq
The following lemma is elementary and its proof is omitted.
\bel\label{l:entropy} For any $\bm w \in \R^N$, we have 
\beq\label{e:entropy}
 \left|  \sup_{i \in \unn{1}{N} } w_i  -  A_{\beta}(\bm w) \right|  < \frac{2 \log N}{\beta}.\eeq
\eel

We let $\{ \mathbf e_i \}_{i=1}^N$ denote the standard basis vectors for $\R^N$, defined by $\mathbf e_a(b) = \one_{a = b}$. For any $M \in \matn$ and $\ell \in \unn{1}{N}$, we set $\bm v_\ell  = \left( N  v_\ell (\mathbf e_i ) \right)_{i=1}^N \in \R^N$, where $v_\ell$ is the regularized eigenvector projection for $M$ defined in \Cref{d:regvect}. 

\bel \label{l:Abeta}
Let $H$ be a generalized Wigner matrix and fix $\ell \in \unn{1}{N}$. Fix $\delta_1>0$ and make the parameter choices \eqref{e:paramchoice}. Define $\beta = N^{\delta_1}$, and set 
\beq 
 \mathcal E = \bigcap_{i \in \unn{1}{N}}  \B(\mathbf e_i ,\omega,\delta_1,\varepsilon_1,\delta_2,\eps_2), \qquad  
\tilde {\mathcal E}=  \bigcap_{i \in \unn{1}{N}}  \tilde \B(\mathbf e_i ,\omega,\delta_1,\varepsilon_1,\delta_2,\eps_2,1, \ell).\eeq
First, there exists $c = c(\delta_1) >0$ such that
\beq\label{e:thefirstestimate}
\one_{\tilde \cE}  \left|  \sup_{i \in \unn{1}{N} } N \langle \u_\ell ,\mathbf e_i \rangle^2 - A_\beta ( \bm v_\ell(H) )\right| \leq 
C N^{-c}.
\eeq
Second,  for all $j \in \unn{1}{5}$,
\beq \label{e:TF}
\sup_{a,b,c,d \in \unn{1}{N} } \sup_{w \in [0,1]}  \one_{\mathcal E} \left|  \partial_{ab}^j A_\beta \left(  \bm v_\ell ( \Theta^{(c,d)}_w H) \right) \right| \leq C N^{10 j (2 \delta_1 + \eps_1 +\eps_2 ) }.
\eeq
Third, we have for $j \in \unn{1}{5}$ that
\beq \label{e:as}\sup_{a,b,c,d \in \unn{1}{N} }\sup_{w \in [0,1]}  \left|  \partial_{ab}^j A_\beta \left(  \bm v_\ell (  \Theta^{(c,d)}_w H) \right)\right| \leq C N^{Cj}. \eeq
Here $C=C(\delta_1)>0$ is a constant. \eel

\begin{rmk}
	Note that in the definition of the event $\tilde{\mathcal{E}}$, we set $k=1$ as a parameter for $\tilde{\mathcal{B}}.$ In other words, we consider the event where $\lambda_\ell$ is the unique eigenvalue in the sub-microscopic interval $I_{\delta_2}(\lambda_\ell)$.
\end{rmk}

\bp By \eqref{e:boundedreclaim}, we have for any $i \in \unn{1}{N}$ that
\beq\label{e:bre1}
\one_{\hat  \cE} v_\ell(\mathbf e_i) \le \langle \mathbf e_i, \u_\ell \rangle^2 + \mathcal O_{\delta_1} \left( N^{-1 + 3\omega + \delta_2 - \eps_2} + N^{-1+  \eps_1 + 5\eps_2 - \delta_1 } \right).
\eeq
By \eqref{e:onlyprove}, 
\beq\label{e:bre2}
\one_{\cE} \langle \mathbf e_i, \u_\ell \rangle^2 \le v_\ell(\mathbf e_i) +  \mathcal O \left( N^{-1 + \omega + \delta_2 - \eps_2 } + N^{-1 + \eps_1 + 5 \eps_2 - \delta_1} \right).
\eeq
We deduce from \eqref{e:bre1} and \eqref{e:bre2} that
\beq\label{e:bre3}
\one_{\tilde \cE} \left| N \langle \mathbf e_i,\u_\ell \rangle^2  - N v_\ell(\mathbf e_i) \right| \le c^{-1} N^{-c}
\eeq
for some $c = c(\delta_1) > 0$. Since \eqref{e:bre3} is independent of $i$, we deduce \eqref{e:thefirstestimate} after taking a supremum over all $i\in \unn{1}{N}$ and using \eqref{e:entropy}.

We now claim that the partial derivatives of $A_\beta (\bm w)$ with respect to the entries of the vector $\bm w \in \R^N$ satisfy
\beq \label{e:bre4}
 \sum_{\underline{j}} \left|\frac{\partial^d A(\bm w)}{\partial_{j_1}\dots \partial_{j_d}}\right| \le C \beta^{d-1} ,
 \eeq
for any $d \in \N$. Here the sum runs over all multi-indices $\underline{j} = (j_1, \cdots, j_d )$ with values in $\unn{1}{N}^d$,  $\partial_j = \partial_{v_j}$, and $C = C(d)>0$ is a constant. This inequality follows by straightforward differentiation, and complete details are given in \cite{landon2018comparison}*{Lemma 3.4}. Using the chain rule, \eqref{e:bre4} and the first inequality in \eqref{e:vderiv} imply \eqref{e:TF}. Similarly, \eqref{e:bre4} and the second inequality in \eqref{e:vderiv} imply \eqref{e:as}.
\ep

Before proceeding, we require the following proposition, which is proved in \Cref{sec:lre}.

\bep\label{l:lr2}
Let $H$ be a generalized Wigner matrix. Then there exists $\delta_0 >0$ such that for all $\delta \in (0, \delta_0)$, there exist constants $C = C(\delta)$ and $c=c(\delta) > 0$ such that for any $i \in \unn{1}{N -1 }$, 
\beq\label{e:optlr}
\P\left( \lambda_{i+1} - \lambda_i < \frac{ N^{-\delta} }{N^{2/3} {\hat i}^{1/3} } \right) < C(\delta) N^{-\delta - c}.
\eeq
\eep

We also introduce some notation. For any $\eps > 0$, we let $f =f_\eps\colon \R \rightarrow \R$ be a smooth, weakly increasing function such that $|f(x)| \le 1$ for all $x\in \R$, $f(x) = 0$ for $x \le (2 + \eps )\log N$, $f(x) = 1$ for $x \ge (2 + 2\eps )\log N$, and the derivatives of $f(x)$ satisfy $|f^{(j)}(x) | \le C (\log N)^j$ for all $ j \in \unn{1}{5}$ and some constant $C=C(\eps) > 0$. For $M\in \matn$ and $\ell\in\unn{1}{N}$, we also define
\beq
S(M) = S_{\beta,\eps} (M)=  f_\eps( A_\beta(\bm v_\ell(M))).
\eeq

\bel\label{l:hsS}
Let $H_0$ be a generalized Wigner matrix and fix $\ell\in\unn{1}{N}$. Fix $\delta_1, \eps >0$ and $s \in [ N^{-1/4}, 1]$, and choose parameters according to \eqref{e:paramchoice}. Then there exist constants $C = C(\delta_1, \eps)>0$ and $c =c(\eps)>0$ such that
\beq
\E \left[ S_{\beta, \eps} (H_s) \right] \le C N^{-c},
\eeq
where $\beta = N^{\delta_1}$.
\eel
\bp
Define $\tilde \cE$ as in \Cref{l:Abeta} and observe that $\P\left( \tilde \cE^c\right) \le CN^{-c}$ for some constants $C= C(\delta)>0$ and $c =c (\delta_1)>0$ by \eqref{e:optlr} and \Cref{l:good1}. Then we have, using the definition of $S(H_s)$ and \eqref{e:thefirstestimate}, that
\begin{align}
\E \left[ S(H_s) \right] & =  \E \left[ \one_{\tilde \cE} f_\eps( A_\beta(\bm v_\ell(M)))\right] +  \E \left[ \one_{\tilde \cE^c} f_\eps( A_\beta(\bm v_\ell(M)))\right] \\ & \le  \E \left[ f_\eps( N \| \u_\ell^s \|^2_\infty   + C_1 N^{-c_1} ) \right]  +CN^{-c} \\ 
&\le 
\P\left( N \| \u_\ell^s \|^2_\infty   + C_1 N^{-c_1} \ge (2 + \eps) \log N  \right)  +CN^{-c}
\end{align}
for some constants $C_1 = C_1 ( \delta_1) >0$ and $c_1(\delta_1) > 0$. For $N \ge N_0(\delta_1)$, we obtain
\beq
\E \left[ S(H_s) \right] \le \P\left( N \| \u_\ell^s \|^2_\infty   \ge \left(2 + \frac{\eps}{2} \right) \log N  \right)+CN^{-c} \le C_2 (\log N)^{1/2} N^{-c_2}+CN^{-c}
\eeq
for some constants $C_2 =C(\eps)>0$ and $c_2 = c_2(\eps) >0$, where in the last inequality we used \eqref{e:gaussdiv2} with $\delta = \eps/10$ and $\theta = 1/10$.  This completes the proof.
\ep

\bel
Let $H$ be a generalized Wigner matrix and fix $\ell \in \unn{1}{N}$ and $\eps >0$. Set  $\delta_1 = 10^{-4}$, $\beta = N^{\delta_1}$, and define $S_{\beta, \eps}$ according to the choice of parameters \eqref{e:paramchoice}. Then we have
\beq\label{e:optconstmoment}
\E \left[S_{\beta, \eps} (H) \right] \le C N^{-c}
\eeq
for constants $C = C(\eps) >0$ and $c = c(\eps) >0$.
\eel

\bp
Define $\mathcal E$ as in the statement of  \Cref{l:Abeta}. By the definition of $f_\eps(x)$ and \Cref{l:Abeta}, we have
\beq\label{e:Sderiv}
\sup_{a,b,c,d \in \unn{1}{N} } \sup_{w \in [0,1]}  \one_{\mathcal E} \left|  \partial_{ab}^j S \left(  \Theta^{(c,d)}_w H \right) \right| \leq C (\log N )^j N^{25 j (2 \delta_1 + \delta_2 + \eps_1 +\eps_2 + \nu) },
\eeq
 and almost surely that
\beq \label{e:Sas}\sup_{a,b,c,d \in \unn{1}{N} }\sup_{w \in [0,1]}  \left|  \partial_{ab}^j S \left(  \Theta^{(c,d)}_w H \right)  \right| \leq C N^{Cj}\eeq
for all $j\in\unn{1}{5}$, where $C > 0$ is a constant. 
We now define $s_1$, $\phi$, $H^\gamma$, and $R$ to have the same meanings as in the proof of \Cref{l:newcomparison}. For a fixed choice of $\gamma = \phi(i,j)$, following that proof gives
\begin{align}
 \label{e:fourthorder2} \E \left[ S \left( H^{\gamma } \right)  \right]- \E \left[S \left(   H^{\gamma -1} \right) \right] &= \frac{1}{4!} \E \left[ \partial^4S \left(  \Theta^{(i,j)}_0 H^\gamma \right) h_{ij}^4 \right] - \frac{1}{4!}\E \left[ \partial^4S \left(  \Theta^{(i,j)}_0 H^\gamma \right) r_{ij}^4 \right] \\
  \label{e:fifthorder2} &+ \frac{1}{5!} \E \left[ \partial^5S \left(  \Theta^{(i,j)}_{w_1(\gamma)} H^\gamma \right) h_{ij}^5\right]
  - \frac{1}{5!} \E \left[ \partial^5S \left(  \Theta^{(i,j)}_{w_2(\gamma)} H^\gamma \right) r_{ij}^5 \right],\end{align}
where $w_1(\gamma), w_2(\gamma) \in[0,1]$ are random variables, and $\partial = \partial_{ij}$. 

Using \eqref{e:Sderiv} and \eqref{e:paramchoice}, we have
\beq \label{e:boundongoodset4}
 \left| \E \left[ \one_{\mathcal E} \partial^4S \left(  \Theta^{(i,j)}_0 H^\gamma \right)  \right] \right| + \left|\E \left[ \one_{\mathcal E} \partial^4 S \left(  \Theta^{(i,j)}_0 H^\gamma \right)  \right] \right| \le C (\log N )^4 N^{100 (2 \delta_1 + \delta_2 + \eps_1 +\eps_2 + \nu) }.
 \eeq
By \Cref{l:good1} and a union bound, we have
\beq\label{e:newbadset1}
\P\left( \mathcal E^c\right) \le C_1 \exp\left({-c_1(\log N)^{c_1 \log\log N}}\right), 
\eeq
for some constants $C_1(\delta_1), c_1(\delta_1)>0$. Using \eqref{e:newbadset1} and \eqref{e:Sas}, we find
\beq\label{e:4thbad}
\left| \E \left[\one_{\mathcal E^c} \partial^4S \left(  \Theta^{(i,j)}_{0} H^\gamma \right) \right] \right|
  +  \left| \E \left[ \one_{\mathcal E^c} \partial^4S \left(  \Theta^{(i,j)}_{0} H^\gamma \right)  \right]\right| \le C N^{ - 10 },
 \eeq
for some constant $C>0$. Then \eqref{e:boundongoodset4}, \eqref{e:4thbad}, and the choice  $\delta_1 = 10^{-4}$ give
\beq\label{e:5thfinal1}
 \left| \E \left[  \partial^4S \left(  \Theta^{(i,j)}_0 H^\gamma \right) h_{ij}^4 \right] -\E \left[  \partial^4 S \left(  \Theta^{(i,j)}_0 H^\gamma \right) r_{ij}^4 \right] \right| \le C N^{ -2 -c }
 \eeq
for some constants $C, c>0$. We also used the independence of $h_{ij}$ and $r_{ij}$ from $ \Theta^{(i,j)}_0 H^\gamma $.

Let $\mathcal D$ be the event where $\sup_{i,j} | r_{ij} | + | h_{ij} |\le C N^{-1/2 + \delta_1}$ holds. Since the variables $r_{ij}$ and $h_{ij}$ are subexponential, we have
\beq\label{e:Dsubexp}
\P\left( \mathcal D^c \right) \le D_1 \exp\left({-d_1(\log N)^{d_1\log\log N}}\right), 
\eeq
for some constants $D_1, d_1 >0$. On $\mathcal D$, we therefore have $\sup_{i,j} | r_{ij} |^5 + | h_{ij} |^5 \le C N^{-5/2 + 5\delta_1}$. 
Using the Cauchy--Schwarz inequality, \eqref{e:Sas} and \eqref{e:Dsubexp}, and that $|h_{ij}|$ and $|r_{ij}|$ have finite 10th moments, we find
\beq\label{e:5thall}
\left| \E \left[\one_{\mathcal D^c}\partial^5S \left(  \Theta^{(i,j)}_{w_1(\gamma)} H^\gamma \right) h_{ij}^5\right] \right|
  +  \left| \E \left[ \one_{\mathcal D^c} \partial^5S \left(  \Theta^{(i,j)}_{w_2(\gamma)} H^\gamma \right) r_{ij}^5 \right]\right| \le C N^{-10}
 \eeq
for some constant $C>0$. We also have

\begin{align}\label{e:5thall2}
\left| \E \left[\one_{\mathcal D}\partial^5S \left(  \Theta^{(i,j)}_{w_1(\gamma)}  H^\gamma \right) h_{ij}^4\right] \right|
  &+  \left| \E \left[ \one_{\mathcal D} \partial^5S \left(  \Theta^{(i,j)}_{w_2(\gamma)} H^\gamma \right) r_{ij}^4 \right]\right|
  \\ \le N^{-5/2 + 5\delta_1} &\left( \E \left[ \left| \partial^5S \left(  \Theta^{(i,j)}_{w_1(\gamma)} H^\gamma \right) \right| \right] 
  +   \E \left[ \left| \partial^5S \left(  \Theta^{(i,j)}_{w_2(\gamma)} H^\gamma \right) \right| \right] \right)\le C N^{-2 -c }\label{e:5thfinal2}
 \end{align}
for some constants $C, c>0$, where the last inequality follows from reasoning similar to \eqref{e:boundongoodset4} and \eqref{e:4thbad}.

Combining \eqref{e:fourthorder2}, \eqref{e:fifthorder2}, \eqref{e:5thfinal1}, and \eqref{e:5thfinal2}, and summing over all $N(N+1)/2$ choices of $\gamma$, we obtain
\beq\label{e:theaboveR}
\left| \E \left[S_{\beta, \eps} (H) \right]  - \E \left[S_{\beta, \eps} (R) \right] \right| \le C N^{-c}.
\eeq 
 The conclusion follows after applying \Cref{l:hsS} to $ \E \left[S_{\beta, \eps} (R) \right]$ in \eqref{e:theaboveR}.\ep

\bp[Proof of \Cref{t:main2}]
Observe \eqref{e:optconstmoment} implies that 
\beq
\P\left( A_{\beta}( \bm v(H) )  \ge (2 + 2\eps) \log N  \right) \le C N^{-c}
\eeq
for some constants $C = C(\eps) >0$ and $c =c (\eps) >0$. Together with \eqref{e:thefirstestimate}, this implies 
\beq
\P\left( \one_{\tilde \cE} \sup_{i \in \unn{1}{N} } N \langle\u_\ell ,\mathbf e_i \rangle  \ge (2 + 2\eps) \log N  + C N^{-c} \right)  \le C N^{-c}
\eeq
after adjusting the values of $C$ and $c$. For $N \ge N_0(\eps)$, we obtain
\beq
\P\left( \one_{\tilde \cE} \sup_{i \in \unn{1}{N} } N \langle\u_\ell ,\mathbf e_i \rangle  \ge (2 + 3\eps) \log N  \right)  \le C N^{-c}.
\eeq
We now observe that there exist constants $C_1, c_1 >0$ such that $\P\left( \tilde \cE^c\right) \le C_1 N^{-c_1}$ by \eqref{e:optlr} and \Cref{l:good1}.
The claim \eqref{e:t2claim1} then follows after redefining $\eps$. The proof of \eqref{e:t2claim2} follows by the same methods and is slightly easier.
\ep

\subsection{Proof of \texorpdfstring{\Cref{t:main3}}{Theorem 1.4}} \label{s:main3}

In this section, we being by reviewing the reverse heat flow technique, then use it to prove \Cref{t:main3}. Since the argument is standard, our treatment is brisk, and we refer the reader to \cite{bourgade2018extreme} for complete details.

First, we say that a random variable $h=h(N)$ is \emph{smooth on the scale $\sigma = \sigma(N)$} if it has a density $e^{-V(x)}$ that satisfies 
\beq\label{e:sigmasmooth}
\left| V^{(k)}(x) \right|  \le C\sigma^{-k} (1 + |x|)^C
\eeq
for all $x \in \R$, $k \ge 0$, and some constant $C>0$. 

The key technical input is \cite{bourgade2018extreme}*{Lemma 4.1}, which says essentially the following. Suppose $H$ is a generalized Wigner matrix whose entries verify definition \eqref{e:sigmasmooth} on scale $\sigma = N^{-a}$ with a uniform constant $C$. Let $t = N^{-b}$ be a time parameter, and suppose $ 0 < 2a < b < 1$. Then there exists a generalized Wigner matrix $\tilde H$ such that
\beq\label{e:TV}
d_{\mathrm{TV}}  \left( H, \sqrt{1-t}\tilde H + \sqrt{t}\mathrm{GOE}_N \right)\le N^{-D}
\eeq
for $N \ge N_0(a,b,D)$, where $d_{\mathrm{TV}}$ denotes the total variation norm and $\mathrm{GOE}_N$ is independent from $\tilde H$.\footnote{\cite{bourgade2018extreme}*{Lemma 4.1} is actually stated in term of an Ornstein--Uhlenbeck process. Reaching the given form requires a straightforward calculation which may be found, for example, as \cite{huang2015bulk}*{(2.17)}.} The intuition here is that one is finding initial data such that its time evolution under a certain parabolic differential equation is the same as $H$; since parabolic differential equations (e.g. the heat equation) are smoothing, this is impossible unless the target data itself has some degree of smoothness.

By \Cref{c:delocdyson}, we obtain strong delocalization results when $t \gg N^{-1/3}$. This corresponds to $b> 1/3$, which explains the requirement of smoothness on a scale $N^{-1/6 +  \vartheta}$ for some $\vartheta >0$ in \Cref{t:main3}.

\begin{proof}[Proof of \Cref{t:main3}] 
Use \cite{bourgade2018extreme}*{Lemma 4.1} to construct $\tilde H$ satisfying \eqref{e:TV}. From \eqref{e:TV}, it follows that for any $D>0$, we have
\beq
d_{\mathrm{TV}}\left( \sup_{\ell\in\unn{1}{N}}\left\Vert \u_\ell \left(H\right) \right\Vert_\infty, \sup_{\ell\in\unn{1}{N}}  \left\Vert \u_\ell\left(\sqrt{1-t}\tilde H + \sqrt{t}\mathrm{GOE}_N\right)  \right\Vert_\infty \right) \le N^{-D}
\eeq
for $N \ge N_0(a,b,D)$. Then by the definition of total variation distance,
\begin{multline}
\P\left(
		\sup_{\ell\in\unn{1}{N}}\left\Vert \u_\ell \left(H\right) \right\Vert_\infty \geqslant 
		\sqrt{\frac{(4+\varepsilon)\log N}{N}}
	\right)
	\\ \le
\P\left(
		\sup_{\ell\in\unn{1}{N}}
		 \left\Vert \u_\ell\left(\sqrt{1-t}\tilde H + \sqrt{t}\mathrm{GOE}_N\right)  \right\Vert_\infty \geqslant 
		\sqrt{\frac{(4+\varepsilon)\log N}{N}}
	\right) + N^{-D}.
\end{multline}
The conclusion follows after taking $D=1$ and observing that the eigenvectors of $\sqrt{1-t}\tilde H + \sqrt{t}\mathrm{GOE}_N$ satisfy \eqref{heatconclude}, by \eqref{e:gaussdiv3}.
\end{proof}

\begin{rmk}
The hypothesis of subexponential decay in \Cref{d:wigner} is necessary to use \cite{bourgade2018extreme}*{Lemma 4.1} as stated. However, as noted in \cite{bourgade2018extreme}*{Section 1}, this hypothesis could be weakened to assuming that only some large moment exists.
\end{rmk}

\appendix

\section{Eigenvalue regularization}\label{a:eigreg}

We begin by presenting some computations to motivate the definition of the regularized eigenvalues $\tilde \lambda_i$.
Let $i= i(N) \in\unn{1}{N}$ be an eigenvalue index. For $m \in \unn{1}{N}$, recall that $\gamma_m$ is the $m$-th quantile of the semicircle distribution, as defined in \eqref{e:classical}. For $m \le 0$, we define $\gamma_m$ such that $\gamma_{m} = \gamma_1 - (1 - m) N^{-2/3}$. For $m  > N$, we set $\gamma_m = \gamma_N +  (m - N) N^{-2/3} $.

Let $\eps_1 >0$ be a parameter to be fixed later. Let $j$ and $k$ be indices such that $i - 2 N^{\eps_1} \le j \le i - N^{\eps_1}$ and $i + N^{\eps_1} \le k \le i + 2N^{\eps_1}$, and define the interval $I = [ \gamma_j, \gamma_k]$.
We define the eigenvalue counting function for any $E \ge - 10$ by
\beq\label{e:countfunc}
\mathcal{N}(E)=\left\vert
	\left\{
		k\in\unn{1}{N}\mid \lambda_k\in [-10,E]
	\right\}
\right\vert
=
\Tr \one_{[-10,E]}(H).
\eeq
We now suppose that rigidity \eqref{e:rigidity} holds with $\omega = \eps_1/10$. We then have that
\beq\label{e:past0}
\lambda_i - \gamma_j  = \int_{\gamma_j}^{\lambda_i}\, \d E = \int_I \one_{ \{\lambda_i \ge E\} } \, \d E= \int_I \one_{ \{\mathcal{N}(E) \le i-1\}}\, \d E.
\eeq
For each $n \in \unn{1}{N}$, let $r_n \colon \R \rightarrow [0,1]$ be a smooth function such that $r_n(x) = 1$ for $x\le n-1$, $r_n(x) = 0$ for $x \ge n - 1/2$, and $|r'_n| + |r'' _n| + |r'''_n| \le C$ for a constant $C>0$ independent of $i$. We write $r = r_i$. Then, because $\mathcal N (E)$ is integer valued, we have
 \beq\label{e:past1}
 \int_I \one_{\{\mathcal{N}(E) \le i-1\}}\, \d E = \int_I r(\mathcal N(E))\, \d E.
 \eeq

 Let $\delta_1>0$ be a parameter to be chosen later, and set $\eta_1 = N^{-2/3 - \delta_1}\hat i ^{-1/3}$. For every $E \in I$ we define the smoothed indicator function $f_E(x)$ such that $f_E(x) = 1$ for $x \in [ - 10, E]$ and $f_E(x) = 0$ for $x \ge E + \eta_1$ and $x \le -11$. We also demand that $|f^{(k)}_E(x) | \le C \eta_1^{-k}$ for $k \in \unn{1}{5} $ when $x\in [E, E + \eta_1]$, and $|f^{(k)}_E(x) | \le C$ for $k \in \unn{1}{5} $ when $x \in [-11,  -10]$, where $C>0$ is a constant independent of $E$ and $\eta_1$.

 By rigidity, we have
 \beq
 | \mathcal N(E)  - \Tr (f_E) | \le | \{ a : \lambda_a \in [E, E + \eta_1]\} |,
 \eeq
 and we can use rigidity again to deduce
 \begin{align}\label{e:past2}
 \left| \int_I r(\mathcal N(E)) \, \d E - \int_I r(\Tr(f_E)) \, \d E \right| &\le C \int_I | \mathcal N(E) - \Tr(f_E) | \, \d E\\
 &\le C \int_I | \{ a : \lambda_a \in [ E , E + \eta_1 ] \} | \, \d E\\
 &= \sum_a \int_I \one_{\lambda_a \in [E , E + \eta_1] }\, \d E\\
 &\le C \eta_1 | \{ a : \lambda_a \in I + [0 , \eta_1] \} | \le C \eta_1 N^{\eps_1} = C {N^{-2/3+\eps_1 - \delta_1}\hat i ^{-1/3}} .
 \end{align}
 By the Helffer--Sj\"ostrand formula \cite{benaych2016lectures}*{Proposition C.1} and a short computation using integration by parts, we obtain
\begin{align}
\tr (f_E ) &= \frac{N}{ 2 \pi } \int_{\R^2} \left(\I f_E (e) \chi' ( \sigma) - \sigma f'_E (e) \chi' ( \sigma ) \right) m_N (e + \I \sigma ) \, \d e \, \d\sigma \\
&+   \frac{N } { 2\pi } \int_{ | \sigma | > \eta_1}  \int_\R f'_E (e)  \partial_{\sigma} ( \sigma \chi ( \sigma ) )   \Re [ m_N (e + \I \sigma ) ] \, \d e \, \d\sigma \\
&- \frac{ N }{  \pi} \int_0^{\eta_1}  \int_{ \R} f''_E (e) \sigma \Im [ m_N (e + \I \sigma ) ] \, \d e \, \d\sigma \\
& + \frac{ N }{  \pi } \int_{\R} f'_E (e) \eta_1 \Re [ m_N (e + \I \eta_1 ) ]\, \d e.
\end{align}
Details may be found in \cite{landon2018comparison}*{Section 4.1}.
Here $\chi$ is a smooth, symmetric function such that $\chi(x)  = 1$ for $|x| \le 1$ and $\chi(x) = 0$ for $|x| > 2$.
 
We now aim to replace $\tr (f_E )$ with an approximation in order to produce the regularized eigenvalue. Set
\begin{align}
F_E &= \frac{N}{  2\pi }\int_{\R^2} \left(\I f_E (e) \chi' ( \sigma) - \sigma f'_E (e) \chi' ( \sigma ) \right) m_N (e + \I \sigma ) \, \d e \, \d\sigma  \label{e:Aterm} \\ 
& +  \frac{N } { 2\pi } \int_{ | \sigma | > \eta_1}  \int_\R f'_E (e)  \partial_{\sigma} ( \sigma \chi ( \sigma ) )   \Re [ m_N (e + \I \sigma ) ] \, \d e \, \d\sigma \label{e:Bterm} \\
& + \frac{ N }{  \pi } \int_{\R} f'_E (e) \eta_1 \Re [ m_N (e + \I \eta_1 ) ]\, \d e. \label{e:problem}
\end{align}
Note that
\begin{align}
\left| \int_{I} r ( \tr f_E ) \, \d E - \int_I r ( F_E) \, \d E \right|  & \leq C \int_I \int_{\R} \int_{0}^{\eta_1} N \sigma |f''_E ( e) |  \Im [ m_N (e + i\sigma ) ] \, \d\sigma \, \d e \, \d E \label{e:fefe}\\
&\leq \frac{C}{ \eta_1^2} \int_{I} \int_0^{\eta_1} \int_0^{\eta_1} N \sigma \Im [ m_N (E+e + i\sigma ) ] \, \d\sigma \, \d e \, \d E \\
& +   C \int_I \int_{-11}^{-10} \int_{0}^{\eta_1} N \sigma  \Im [ m_N (e + i\sigma ) ] \, \d\sigma \, \d e \, \d E\\
&= \frac{C}{ \eta_1^2} \int_0^{\eta_1} \int_0^{ \eta_1 } \left[ \int_I N \sigma \Im [ m_N (E + e+ i\sigma ) ] \, \d E \right]  \d e \, \d\sigma \label{e:I1}\\
& +   C |I|  \int_{-11}^{-10} \int_{0}^{\eta_1}   N \sigma  \Im [ m_N (e + i\sigma ) ] \,  \d\sigma \, \d e . \label{e:newerr}
\end{align}
We begin by bounding \eqref{e:I1}. Define $E' = E + e$.  Since $\sigma\le \eta_1$ in the integral above, we have
\begin{align}
\int_I N \sigma  \Im [ m_N (E' + i\sigma ) ] \, \d E &\le  \sum_{ a : |i-a| \leq 3 N^{\eps_1} } \eta_1 \int_{I} \frac{\sigma}{ (E' - \lambda_a )^2 + \sigma^2 } \, \d E 
+ \sum_{ a : |i-a| > 3 N^{\eps_1} } \eta_1 \int_{I} \frac{\sigma}{ (E' - \lambda_a )^2 + \sigma^2 } \, \d E. 
\end{align}
We estimate the first term by 
\beq\label{e:ccc2}
 \sum_{ a : |i-a| \leq 3 N^{\eps_1} } \eta_1 \int_{I} \frac{\sigma}{ (E' - \lambda_a )^2 + \sigma^2 } \, \d E \leq  \sum_{ a : |i-a| \leq 3 N^{\eps_1} } \eta_1 \int_{\R} \frac{\sigma}{ (E' - \lambda_a )^2 + \sigma^2 } \, \d E \leq C N^{- 2/3 + \eps_1 - \delta_1}\hat i ^{ - 1/3} .
\eeq
 
 For the second term, we consider only the case $ i < N/4$ explicitly. The case $i > N - N/4$ follows by symmetry. The case $i \in \unn{N/4}{3N/4}$ is similar and was already treated in \cite{landon2018comparison}*{Section 4.1}, so we do not give details here.
 
In the second term, we consider the sum over the terms with $i - a  > 3N^{\eps_1}$ and $a -i   >   3 N^{3\eps_1}$ separately. When $i - a  > 3N^{\eps_1}$, we are considering eigenvalues $\lambda_a$ nearer to the edge $-2$ than $\lambda_i$. For such $a$, we have $N^{-2/3}\hat i ^{-1/3}\leqslant N^{-2/3}\hat a ^{-1/3}$ and $N^{-2/3}\hat a ^{-1/3} \le (i - a + 1)^{1/3} N^{-2/3}\hat i ^{-1/3}$.
 Then by rigidity \eqref{e:rigidity} and the fact that $a < j < i$,
 \begin{align}
 |\gamma_j  - \eta_ 1 - \lambda_a | &\ge | \gamma_j - \gamma_a| - |\gamma_a - \lambda_a|  - \eta_1 \\
 &\ge (j-a) N^{-2/3}\hat i ^{-1/3}  -(i - a + 1)^{1/3} N^{-2/3+\eps_1/2}\hat i ^{-1/3} - N^{-2/3}\hat i ^{-1/3}\\
 &\ge (9/10) ( i -a - 2N^\eps)  N^{-2/3}\hat i ^{-1/3}
 \end{align}
 for $N \ge N_0(\eps_1)$. 
 We also observe that $I$ satisfies 
 \beq\label{e:Ibound}
 | I | \le C N^{-2/3 + \eps_1} {\hat i}^{-1/3}
 \eeq
 for some constant $C >0$. Using \eqref{e:Ibound}, $ \sigma \le \eta_1$, and the definitions of $I$ and $E'$, we find
\begin{align}
\sum_{ a :  i-a > 3 N^{\eps_1} } \eta_1 \int_{I} \frac{\sigma}{ (E' - \lambda_a )^2 + \sigma^2 } \, \d E 
&\leq C\eta_1 \sigma  | I |   \left( N^{2/3}\hat i ^{1/3} \right)^{2} \sum_{ a : i-a  > 3 N^{\eps_1} } \frac{1}{( (9/10)(i- a - 2 N^{\eps_1} ))^2 }  \\
&\label{e:ccc1} \leq C  \eta^2_1 |I| \left(N^{2/3}\hat i ^{1/3}\right)^{2} \leq  C |I| N^{-2\delta_1} \le N^{-2/3+ \eps_1 - 2\delta_1}\hat i ^{-1/3}.
\end{align}

For $a -i   >   3 N^{3\eps_1}$, we recall the definition of $I$ and write 
\begin{align}
\sum_{ a :  a- i  > 3 N^{\eps_1} } \eta_1 \int_{I} \frac{\sigma}{ (E' - \lambda_a )^2 + \sigma^2 } \, \d E \label{e:430}
&\leq \eta_1 \sigma  | I |   \sum_{ a :  a- i  > 3 N^{\eps_1} } \frac{1}{(\lambda_a - \gamma_{k}  - \eta_1 )^2} .
\end{align}
We note that, supposing rigidity \eqref{e:rigidity} holds with $\omega = \eps_1/10$, $a< N/2$, and $a - i > 3N^{\eps_1}$, then
\begin{align}\label{e:citeatend}
 \lambda_a - \gamma_{k}  - \eta_1  \ge \frac{1}{2}(\gamma_a -  \gamma_i) 
& \ge  c \left( N^{-2/3} [ a^{2/3}  - i^{2/3} ] \right)
\end{align}
for $N \ge N_0(\eps_1)$, after possibly increasing $N_0$.
Here we used the definition \eqref{e:classical} to lower bound $\gamma_a -  \gamma_i$.\footnote{As an intermediate step, we use for $j\le N/2$ that $\gamma_j \le -2+C \left( \frac{j}{N} \right)^{2/3}$, which follows from \eqref{e:classical}. The claim follows from using \eqref{e:classical} again and the mean value theorem to write \begin{equation*}\frac{1}{N} \le C (( 2 + \gamma_{j+1})^{3/2} - (2 + \gamma_j)^{3/2} ) \le 2 C (\gamma_{j+1}- \gamma_j ) \sqrt{2 + \gamma_{j+1}} \le   4 C (\gamma_{j+1}- \gamma_j ) \sqrt{2 + \gamma_{j}}\leqslant C(\gamma_{j+1}-\gamma_j)\left(\frac{j}{N}\right)^{1/3}. \end{equation*}} We also used the fact that $| \lambda_a - \gamma_a | \ll | \gamma_a -  \gamma_i|$, which follows from the inequality
\beq
N^{-2/3} \left( a^{2/3}  - i^{2/3} \right) \ge  \frac{N^{\eps_1/9}}{2 N^{2/3} { i}^{1/3}} \gg  \frac{N^{\eps_1/10}}{N^{2/3} { i}^{1/3}}.
\eeq
Then
\begin{align}\label{e:bigevsum}
\sum_{ a : a - i  >  3N^{\eps_1}, \, a < N/2} \frac{1}{(\lambda_a - \gamma_{k}  - \eta_1 )^2}  & \le C N^{4/3} \sum_{ a : a - i  >  3N^{\eps_1}, \, a < N/2} \frac{1}{(a^{2/3}  - k^{2/3}  )^2}  \\
&\le C N^{4/3} \int_{i+ N^{\eps_1}}^\infty \frac{\d x}{(x^{2/3} - i^{2/3})^2} \le C N^{4/3} i^{2/3}.
\label{e:theintegralbound}
\end{align}
Here we used the substitution $x = i y$ to bound the integral in \eqref{e:theintegralbound}:
\beq
 \int_{i+ N^\eps}^\infty \frac{\d x}{(x^{2/3} - i^{2/3})^2} =  i^{-1/3} \int_{1 + \frac{N^{\eps_1}}{i} } \frac{\d y }{(y-1)^2} = i^{-1/3} \int_{\frac{N^{\eps_1}}{i} } \frac{\d y }{y ^2} \le C i^{2/3}.
\eeq
Since we assumed $i < N/4$, we have $|\gamma_m - \lambda_a | > 1/10$ for $a > N/2$ when rigidity \eqref{e:rigidity} holds. This implies
\beq
\sum_{ a : a - m  >  N^{\eps_1}, \, a > N/2} \frac{1}{(\lambda_a - \gamma_{k}  - \eta_1 )^2}  \le CN \le N^{4/3}\hat i ^{2/3}.
\eeq
We conclude that
\beq\label{e:ccc0}
\eqref{e:430} \le N^{\eps_1-2\delta_1} |I| = N^{- 2/3+ 3 \eps_1 - \delta_1}\hat i ^{- 1/3}.
\eeq
Collecting \eqref{e:ccc2}, \eqref{e:ccc1}, and \eqref{e:ccc0}, we have shown that the term  $\eqref{e:I1}$ is  bounded by $C N^{- 2/3+ 3 \eps_1 - \delta_1}\hat i ^{- 1/3}$ when $i \le N/4$. 

Finally, we consider \eqref{e:newerr}. We observe that $\Im [ m_N (e + i\sigma ) ] \le C N^{-1}$ for $\sigma \in (0, \eta_1)$ and $e \in [-11, -10]$ by the local law outside of the spectrum \cite{benaych2016lectures}*{(10.2)} and \eqref{e:msclower2}, and we obtain
\beq
\left| \eqref{e:newerr}\right|  \le  C N^{- 2/3 + \eps_1 - \delta_1}\hat i ^{- 1/3}.
\eeq
 Therefore, $F_E$ is a good approximation to $\tr (f_E )$, which we saw in \eqref{e:past0}, \eqref{e:past1}, and \eqref{e:past2} controls the difference $\lambda_i - \gamma_j$. This leads to the following definition.

\bed \label{d:regev}
Let $M \in \matn$ be a symmetric matrix, and fix $\delta, \eps$. We set $\eps_1 = \eps/3$, $\delta_1 = \delta$ in the previous computation and define the regularized eigenvalues $\tilde \lambda_i$ for $i \in \unn{1}{N}$ by 
\beq
\tilde \lambda_i (M) = \tilde \lambda_i(M, \delta, \eps)  = \int_I r(F_E) \, \d E + \gamma_j.\label{e:lambdatildedef}
\eeq
\eed

\begin{proof}[Proof of \Cref{p:regulareigval}]
Given $\delta, \eps >0$, we define $\tilde \lambda_i \colon \matn \rightarrow \R$ by \Cref{d:regev}. 
We let $\mathcal A_0 = \mathcal A_0(\eps, \delta)$ be the set from \Cref{l:goodset} with $\omega = \eps /100$, where in particular rigidity and delocalization \eqref{e:rigidity} hold, following our choice before \eqref{e:past0} (and actually choosing a smaller parameter $\omega$). By \Cref{l:goodset}, the claimed probability bound on $\mathcal A_0$ holds. 

First, the remarks following \eqref{e:ccc0} show the first equation in \eqref{e:regev} when $i < N/4$ or (by symmetry) when $ i > N - N /4$. The remaining case, where $i$ is in the bulk, is similar, and the necessary bound for \eqref{e:I1} is provided in the proof of \cite{landon2018comparison}*{Lemma 3.2}.

Next, we estimate the derivatives of $\tilde \lambda_i(H)$. We first estimate the derivatives of $F_E$. We recall the Green's function differentiation formula\footnote{This follows from a resolvent expansion to first order. See \cite{benaych2016lectures}*{(2.3)}.}
\beq
\frac{\partial G_{ab}(z)}{\partial H_{cd}} = - G_{ac} (z)G_{db}(z)
\eeq
for $a,b,c,d \in \unn{1}{N}$. Using it, we find for $z = E  + \I \eta$ that
\beq\label{e:wardexample}
 N \left| \partial_{bc}  m_N (E + \I \eta) \right|  \le C \sum_a | G_{ab} G_{ca} |   \leq C \sum_a | G_{ab} |^2 + |G_{ca} |^2 \leq \frac{C}{ \eta} \left( | G_{bb} (z) | + |G_{cc} (z) | \right),
\eeq
where in the last line we used the Ward identity \cite{benaych2016lectures}*{(3.16)}. In the same way, we have
\beq
N | \partial_{bc}^k  m (z) | \leq \frac{ C}{ \eta} \left( | G_{bb} (z) | + |G_{cc} (z) | + |G_{bc} (z) | \right)^{k}.
\eeq

We now claim that that $|G_{ij}(E + \I \eta_1 ) | \le  CN^{ \eps/4 + \delta} $  for $|E| \le \eps^{-1}$ on the event $\mathcal A_0$, for some $C>0$. This follows, for example, from \Cref{l:goodset} and \cite{benaych2016lectures}*{Lemma 10.2} applied with $M = N^{\delta + \eps/8 }$ and  $\eta = N^{-2/3 + \eps/8 }\hat i ^{-1/3}$.
Therefore 
\begin{align}
\left| \partial_{bc} ^k \left( \frac{2 N }{  \pi } \int_{\R}  f'_E (e) \eta_1 \Re [ m_N (e + \I \eta_1 ) ]\, \d e\right) \right| &\le  C \int_{\R} \left|f'_E (e)\right| \left( | G_{bb} (e + \I \eta_1) | + | G_{bc} (e + \I \eta_1) | + |G_{cc} (e + \I \eta_1) | \right)^k\,\d e \\
& \le C \int_{\R}  \left| f'_E (e)\right|  N^{ k (\delta + \eps/4) }\, \d e \le C N^{ k (\delta + \eps/4) } .
\end{align}
The derivatives of the other terms may be bounded similarly, and we obtain $\left| \one_{\mathcal A_0}  \partial^k_{bc} F_E  \right| \le  C (\log N) N^{ k (\delta + \eps/4) }$ for $k \in \unn{1}{5}$ (with the $\log N$ factor coming from the integration in $\sigma$ for the derivative of \eqref{e:Bterm}). Since $r$ is bounded with bounded derivatives, we have
\begin{align}
\one_{\mathcal A_0}\left| \partial^k_{ab} \tilde \lambda_i (H) \right| \le  \left| \int_I \partial^k_{ab}r(F_E) \, \d E  \right| &\le C \int_I  (\log N) N^{ k (\delta + \eps/4) }  \, \d E \\  &\le C N^{- 2/3 + \eps/3  }{\hat i}^{-1/3} N^{ k (\delta + \eps/3) } \le  C N^{- 2/3 + k (\delta + \eps) }{\hat i}^{-1/3}
\end{align}
for all $k \in \unn{1}{5}$ and $a,b \in \unn{1}{N}$, and $N \ge N_0(\eps)$. We used here $|I| \le N^{- 2/3 + \eps/3}{\hat i}^{-1/3}$, which follows from our definition of $I$ about \eqref{e:countfunc} and our choice $\eps_1 = \eps/3$ in \Cref{d:regev}. This proves the second inequality in \eqref{e:regev} for $H$. The proof of \eqref{e:detbounds} for $H$ is similar and uses the trivial inequality $\left| G_{ij } ( E  + \I \eta ) \right| \le \eta^{-1}$.

We next consider the statements involving rank-one perturbations of $H$. By a resolvent expansion to high order it is straightforward to prove that
\beq\label{e:conc1}
\sup_{a,b \in \unn{1}{N} } \sup_{  0 \leq w \leq 1 } \one_{\A} \left| \frac{1}{N} \tr \frac{1}{H-z} - \frac{1}{N} \tr \frac{1}{ \Theta^{(a,b)}_w H  -z} \right| \leq  \frac{C N^{\eps/6}}{N \eta}
\eeq
and
\beq\label{e:conc2}
\sup_{a,b \in \unn{1}{N} } \sup_{i,j \in \unn{1}{N}} \sup_{ 0 \leq w \leq 1 } \one_{\A} \left| \left( \frac{1}{H-z} \right)_{ij} - \left( \frac{1}{\Theta^{(a,b)}_w H -z } \right)_{ij} \right| \leq  \frac{CN^{\eps/6}}{\sqrt{N}}
\eeq
for any $z \in \mathcal D_{\eps/100}$, where $C = C(\eps) >0$ and where $\A = \A_0 \cap \left(  \cup_{i,j \in \unn{1}{N} }\left \{ | h_{ij} | \le N^{-1/2 + \eps/12 } \right\} \right)$.  For \eqref{e:conc1}, we also used the Ward identity, as in \eqref{e:wardexample}. Hence, by standard arguments rigidity and delocalization \eqref{e:rigidity} hold simultaneously for $\Theta^{(a,b)}_w H$ for all choices of $w\in [0,1]$ and $a,b \in \unn{1}{N}$  on $\A$ (see \cite{benaych2016lectures}*{Theorem 2.10} and \cite{benaych2016lectures}*{Section 9}).  With these estimates, we may follow the previous reasoning to obtain \eqref{e:regev} and \eqref{e:detbounds} for all such $\Theta^{(a,b)}_w H$, where the probability bound \eqref{e:regevpbound} holds for $\mathcal A$ by the subexponential decay hypothesis \eqref{e:subexphypo}, after decreasing $C_1$ and $c_1$, if necessary. \end{proof}

\section{Level repulsion estimates}\label{sec:lre}

In this section, we provide proofs of \Cref{p:lrmain} and \Cref{l:lr2}, which were necessary for the proofs of the our main results in \Cref{s:main}. We begin in Section~\ref{s:gaussdiv} by proving a level repulsion estimate for the $\mathrm{GOE_N}$, and then extending it to generalized Wigner matrices with small additive Gaussian perturbations. In Section~\ref{s:lrwigner}, we use a comparison argument to extend these estimates to arbitrary generalized Wigner matrices.

\subsection{Level repulsion for Gaussian divisible ensembles} \label{s:gaussdiv}

In this section only, we consider a matrix-valued stochastic process slightly different from the one defined in \eqref{e:dysondyn}. Given a generalized Wigner matrix $H_0$, we define $H_s$ for times $s > 0$ through the Ornstein--Uhlenbeck dynamics
\beq\label{eq:defdyn}
\d H_s = \frac{1}{\sqrt{N}}\, \d B_s - \frac{1}{2}H_s.
\eeq 
We denote by $\lambda_1(s)\leqslant \dots\leqslant  \lambda_N(s)$ the eigenvalues of $H_s$. 

For any energy $E \in \R$, we define $\kappa(E)={\max(N^{-2/3}, \min (\vert E+2\vert, | E- 2| ))}$. 
Given $\delta>0$ and $\fa >0$, we define the interval
\beq\label{eq:definte}
I = I(E) = I(\delta, \fa , E) =\left[
	E - \frac{\fa N^{-\delta}}{N\sqrt{\kappa(E)}},\,E+ \frac{\fa N^{-\delta}}{N\sqrt{\kappa(E)}}
\right].
\eeq

Let $H$ be a generalized Wigner matrix with eigenvalues $\bm{\lambda}=(\lambda_1,\dots\lambda_N)$, labeled in increasing order. Given a choice of interval $I$, we denote the counting functions for the eigenvalues $\bm{\lambda}$ and corresponding regularized eigenvalues $\tilde {\bm{\lambda}}$ of $H$, defined in \Cref{p:regulareigval}, by
\beq
\mathcal{N}_{\bm{\lambda}}(I)=
\left\vert
	\left\{
		i \in\unn{1}{N}\mid \lambda_i \in I
	\right\}
\right\vert
\quad\text{and}\quad
\mathcal{N}_{\tilde{\bm{\lambda}}}(I)=
\left\vert
	\left\{
		i \in\unn{1}{N}\mid \tilde \lambda_i \in I
	\right\}
\right\vert.
\eeq

\subsubsection{Bulk}

In the following theorem, the cases $k=1,2$ were established for energies bounded away from the spectral edges $\pm 2$ in \cite{landon2017convergence}*{Theorem 5.1}. The extension to $k > 2$ and energies slowly tending to $\pm 2$ is straightforward, and we give the details for completeness. Note that level repulsion estimates for any $k$ were also proved in the bulk in \cite{bourgade2016fixed} for larger times $s$.

\bel\label{l:deformedpacking}
Let $H$ be a generalized Wigner matrix, and let $\bm \lambda = \bm \lambda(s)$ be the vector of eigenvalues for $H_s$, as defined in \eqref{eq:defdyn}. For all integers $k\ge 1$,  all $\fa , \delta, \sigma>0$, and $s\in [N^{-1/2}, 1]$, there exists $N_0 = N_0(\fa , k , \delta, \sigma) >0$ such that for all $E \in [-2 + N^{-\sigma/2}, 2 - N^{-\sigma/2}]$ and $N \ge N_0$,
\beq\label{e:gaussianbulklr}
\P \left( \mathcal N_{\bm \lambda}(I) \ge k  \right) \le N^\sigma N^{ -  \delta k/2 },
\eeq
where $I = I(\delta, \fa , E)$.
\eel

\begin{proof}
We first consider the case $k = 1$. Let $m_N(z)$ and $G(z)$ be the Stieltjes transform and Green's function for $H_s$.
Set $\eps = \fa N^{-\delta}$ and define $\eta$ by the equality $\eps = N \eta \sqrt{\kappa(E)} $. 
Following \cite{landon2017convergence}*{(5.4)}, we have
\begin{align}
\P(\mathcal N_{\bm \lambda}( I ) \ge 1 ) \le \frac{4 \eps^2}{\kappa(E)} \E\left[( \Im m_N ( E + \I \eta) )^2\right] &\le \frac{4 \eps^2}{N^2 \kappa(E)} \sum_{i,j} \E\left[ |G_{ii}| |G_{jj} | \right]\\ &\le \frac{4 \eps^2}{N^2 \kappa(E)} \sum_{i,j} \E\left[ |G_{ii}|^2\right]^{1/2} \E \left[ |G_{jj} |^2 \right]^{1/2}.\label{e:lr1}
\end{align}
These inequalities follow directly from the definition of $\Im m_N$. Next, the computations leading to \cite{landon2017convergence}*{(5.36)} show that for any $r>0$, there exists $N_0(\fa, \delta, \eps, r)$ such that, for $N \ge N_0$,
\beq\label{e:lr2a}
\frac{\eps^2}{N^2} \sum_{i,j} \E\left[ |G_{ii}|^2\right]^{1/2} \E \left[ |G_{jj} |^2 \right]^{1/2} \le N^{\sigma/2} N^{ - \delta ( 1 - r) }.
\eeq
Since $E \in [-2 + N^{-\sigma/2}, 2 - N^{-\sigma/2}]$, we therefore obtain from \eqref{e:lr1} and \eqref{e:lr2a} that 
\beq
\P(\mathcal N_{\bm \lambda}( I ) \ge 1 ) \le 4 \kappa(E)^{-1} N^{\sigma/2} N^{ - \delta ( 1 - r) } \le 4 N^{\sigma} N^{ - \delta ( 1 - r) }.
\eeq
Taking $r=1/3$, this proves the claim for $k=1$. 

We now proceed by induction and suppose the claim holds for $k-1$. We will prove it for $k$. Following \cite{landon2017convergence}*{(5.37)}, we have
\begin{align}
\P(\mathcal N_{\bm \lambda}( I ) \ge k ) &\le  \frac{\eps^2}{\kappa(E)} \E\left[\one_{ \{\mathcal N_{\bm \lambda} \ge k\}}( \Im m_N ( E + \I \eta) )^2\right]\\ 
&\label{e:bulklrinductstep} \le \frac{\eps^2}{N^2 \kappa(E)} \sum_{i,j} \E\left[ \one_{ \{\mathcal N^{(i)}_{\bm \lambda}(I) \ge k -1 \}}|G_{ii}|^2\right]^{1/2} \E \left[ \one_{ \{\mathcal N^{(j)}_{\bm \lambda}(I) \ge k -1\} }|G_{jj} |^2 \right]^{1/2}.
\end{align}
Here we used $\mathcal N^{(i)}_{\bm \lambda}$ to denote the eigenvalue counting function for the $i$-th minor of $\GOE_N$, where the $i$th row and column are removed. 
We also used the Cauchy interlacing theorem to see that between any two eigenvalues of $H$ lies an eigenvalue of the $i$-th minor. 
Let $r>0$ be a parameter.
It then follows from the first inequality in \cite{landon2017convergence}*{(5.14)} and the bound \cite{landon2017convergence}*{(5.39)} that
\beq
\E\left[ \one_{\{ \mathcal N^{(i)}_{\bm \lambda}(I) \ge k -1\} }|G_{ii}|^2\right] \le N^{-k \delta} + (( \lambda_i(0) -E )^2 + s^2 )^{-1} N^{\sigma/8}  N^{\delta r} \left( \P(  \mathcal N^{(i)}_{\bm \lambda}(I) \ge k - 1)  \right)^{1/(1+r) } 
\eeq
for $N \ge N_0(\fa , k ,r ,\delta, \sigma)$. We use \eqref{e:gaussianbulklr} with the parameters $k-1$, $\sigma/8$, and $\delta$, and set $r = \min(1/(k-2), 1/100)$ so that $\frac{1}{1+r} \ge \frac{k-2}{k-1}$ for $k\ge 3$. Then\footnote{Note that the minor in the definition of $\mathcal N^{(i)}_{\bm \lambda}(I)$ is not of the form $M_s$ for a generalized Wigner matrix $M$, but $[(N-1)/N]^{1/2} M_s$. This means we must bound the probability eigenvalues of $M_s$ fall in the interval $[N/(N-1)^{1/2}] I$. This causes no problems in our application of $\eqref{e:gaussianbulklr}$ if we use $2\fa$ in place of $\fa$ to accommodate this scaling.}
\begin{align}
\E\left[ \one_{ \{\mathcal N^{(i)}_{\bm \lambda}(I) \ge k -1\} }|G_{ii}|^2\right] & \le N^{-k \delta} + 
 (( \lambda_i(0) -E )^2 + s^2 )^{-1} N^{\sigma/8}  N^{\delta r} \left(N^{\sigma/8} N^{ -  (k-1) \delta/2} \right)^{1/(1+r) } \\ 
&\le N^{-k \delta} + (( \lambda_i(0) -E )^2 + s^2 )^{-1} N^{\sigma/4}  N^{-(\frac{k}{2} - 2)\delta}, 
\label{e:putit2}
\end{align}   
when $N \ge N_0(\fa ,\delta, \sigma ,k)$, where we adjusted $N_0$ upward if necessary.
We now recall from \cite{landon2017convergence}*{(5.35)} that
\beq\label{e:sqput}
\frac{1}{N} \sum_{1 \le i \le N} \frac{1}{\sqrt{(\lambda_i(0) - E)^2 + s^2 } } \le N^{\sigma/16}
\eeq
for large enough $N \ge  N_0 (\sigma)$.

Putting \eqref{e:putit2} into \eqref{e:bulklrinductstep} and using \eqref{e:sqput} together with $\sqrt{a+b} \le \sqrt{a} + \sqrt{b}$ gives
\beq
\P(\mathcal N_{\bm \lambda}( I ) \ge k )\le \frac{\eps^2}{\kappa(E)} N^{3\sigma/8}  N^{-(\frac{k}{2} - 2)\delta} \le \fa^2 N^{-2\delta} N^{ 7\sigma/ 8 }  N^{-(\frac{k}{2} - 2)\delta} \le  N^{\sigma }  N^{ -  \delta k/2 } .
\eeq
when $N \ge N_0(\fa ,\delta, \sigma ,k)$. This completes the proof.
\end{proof}

\subsubsection{Edge}


Before proceeding the Gaussian divisible case, we state a level repulsion estimate at the edge for the GOE. Its proof is given below, after a series of preliminary results.

\bep\label{p:lregoe}
Let $\bm \lambda$ be the vector of eigenvalues of $\GOE_N$.
Fix $k\in\mathbb{N}$ and $\fa, \fb>0$. Let $E\in\mathbb{R}$ be such that $\kappa(E)\leqslant N^{-\fb}$. Then for any $\delta > 0$, there exists a constant $C=C(\fa, \fb, \delta)$ such that
	\beq
	\P \left(
		\mathcal{N}_{\bm{\lambda}}(I )\geqslant k
	\right)
	\leqslant C N^{-\delta k/2},
	\eeq
	where  $I = I(\delta, \fa, E)$.

\eep

Our argument for the GOE is inspired by the methods of \cite{gustavsson}, and in particular, we provide an improvement on \cite{gustavsson}*{Lemma 2.2}.
The proof of \Cref{p:lregoe} uses the determinantal structure of the eigenvalue distribution for a related ensemble and a slight generalization of Wegner estimates in \cites{gustavsson, o2010gaussian}. Since the eigenvalue structure of the GOE is not determinantal but given in terms of a Pfaffian, we first obtain the result for the Gaussian Unitary Ensemble and then use the relationships between the Orthogonal and Unitary ensemble from \cite{forrester2001inter} to obtain the final result. We first recall the definition of the Gaussian Unitary Ensemble.

\bed\label{d:gue}
A $N\times N$ Hermitian matrix $H$ is distributed as $\operatorname{GUE_N}$ if for all $i,j\in\unn{1}{N}$ such that $i\leqslant j$,
\beq
H_{ij}\overset{(d)}{=} \frac{1}{\sqrt{2N}}\mathcal{N}_{ij}+\frac{\I}{\sqrt{2N}}\mathcal{N}_{ij}
\quad\text{if }i\neq j\quad\text{and }
H_{ii}\overset{(d)}{=}\frac{1}{\sqrt{N}}\mathcal{N}_{ii},
\eeq
where $(\mathcal{N}_{ij})_{1\leqslant i\leqslant j\leqslant N}$ is a family of independent standard Gaussian random variables.
\eed

\begin{rmk}
All of our level repulsion results have a corresponding Hermitian equivalent. In particular, \Cref{p:lregoe} can be stated for the $\GUE_N$ defined in \Cref{d:gue}, and this result follows from equation \eqref{eq:lregue} in the proof of \Cref{p:lregoe}. Note that 
we actually obtain a stronger bound than the one for the $\GOE_N$.
\end{rmk}

For the $\GUE_N$, the eigenvalue distribution is exactly computable  and  can be expressed as the following distribution on the simplex $\Sigma =\{\lambda_1\leqslant \lambda_2\leqslant\dots\leqslant \lambda_n\}$:
\beq\label{eq:eiggue}
p_N(\lambda_1,\dots,\lambda_N)=\frac{1}{Z_N}\prod_{i<j}(\lambda_i-\lambda_j)^2\exp\left({-\frac{N}{2}\sum_{i=1}^N\lambda_i^2}\right),
\eeq
where $Z_N$ is a constant factor that may be determined explicitly \cite{anderson2010introduction}*{Theorem 2.5.2}. It was observed by Dyson \cite{dyson1970correlations} that this joint eigenvalue density can be written as a determinant in terms of orthogonal polynomials:
\beq
p_N(\lambda_1,\dots,\lambda_N)=\frac{1}{Z'_N}\det\left(
	K_N(\sqrt{N}\lambda_i,\sqrt{N}\lambda_j)
\right)_{i,j=1}^N	,
\eeq
where $Z'_N$ is a constant, and
\beq K_N(x,y)
=
\sum_{k=0}^{N-1} \psi_k(x)\psi_k(y),\quad
\psi_k(x)=\frac{H_k(x)\e^{-{x^2}/{4}}}{(2\pi)^{\frac{1}{4}}(k!)^{\frac{1}{2}}},
\quad 
H_k(x)=(-1)^k\e^{x^2/2}\frac{\d^k}{\d x^k}\e^{- x^2/2},
\eeq
where $H_k$ the $k$-th Hermite polynomial.

We do not use the explicit definition of the kernel $K_N(x,y)$, only that the eigenvalue process for $\GUE_N$ has a determinantal structure. We have the following theorem from \cite{hough2006determinantal}, which says that the number of points from a determinantal process in a bounded domain is given by a sum of independent Bernoulli random variables. We state it here in the case of the $\GUE_N$. 

\bet[\cite{hough2006determinantal}*{Theorem 7}] \label{theo:bernoulli}
Let $I\subset \mathbb{R}$ be a closed interval and et $\bm{\lambda}$ be the vector of eigenvalues of a $\GUE_N$. Then there exists $p_1,\dots,p_N >0$ such that 
\beq
\mathcal{N}_{\bm{\lambda}}(I)=\sum_{i=1}^N\mathrm{Ber}(p_i)
\quad\text{with}\quad
\mathrm{Ber}(p)=\left\{
	\begin{array}{ll}
		1&\text{ with probability }p,\\
		0&\text{ with probability } 1-p,
	\end{array}
\right.
\eeq
and the $\mathrm{Ber}(p_i)$ random variables are independent.
\eet
This gives as a direct corollary a level repulsion estimate for eigenvalues of a $\operatorname{GUE_N}$.
\bec\label{coro:lredet}
Let $\bm{\lambda}$ be the vector of eigenvalues of a $\GUE_N$ and  $I\subset \mathbb{R}$ be a closed interval. Then for all $k\in\mathbb{N}$,
\beq
\P \left(
	\mathcal{N}_{\bm{\lambda}}(I)\geqslant k
\right)\leqslant \e^{-\lambda k+(\e^\lambda-1)\E[\mathcal{N}_{\bm{\lambda}}(I)]}
\eeq
\eec
\begin{proof}
Using Markov's inequality and \Cref{theo:bernoulli}, we have for any $\lambda>0$ that
\beq
\P\left(
	\mathcal{N}_{\bm{\lambda}}(I)\geqslant k
\right)
\leqslant \e^{-\lambda k}\E\left[
	\e^{\lambda \E{\mathcal{N}_{\bm{\lambda}}(I)}}
\right]
=
\e^{-\lambda k}\prod_{i=1}^N\E\left[\e^{\lambda \mathrm{Ber}(p_i)}\right]
=
\e^{-\lambda k}\prod_{i=1}^N \left(1+p_i(\e^\lambda -1)\right)
\leqslant
\e^{-\lambda k+(\e^\lambda-1)\sum_{i=1}^N p_i}.
\eeq
The result is obtained by noting that $\E[\mathcal{N}_{\bm{\lambda}}(I)]= \sum_{i=1}^N p_i$.
\end{proof}

From \Cref{coro:lredet}, we see that we now only need to control $\E[\mathcal{N}_{\bm{\lambda}}(I)]$ in order to bound the probability that  there are more than $k$ eigenvalues in $I$. For a sub-microscopic interval, the probability of seeing even one eigenvalue is small, so the expectation of the number of eigenvalues tends to zero. We make this claim precise in the following proposition.

\bep\label{prop:gustavsub}
Let $\bm{\lambda}$ be the vector of eigenvalues of a $\GUE_N$. Fix $
\fa, \fb >0$ and let $E\in\mathbb{R}$ be such that $\kappa(E)\leqslant N^{-\fb}$. Then for any $\delta > 0$,
\beq
\E \left[
	\mathcal{N}_{\bm{\lambda}}(I(\delta, \fa, E))
\right]\leqslant CN^{-\delta},
\eeq
where $I(\delta, \fa, E)$ is defined in \eqref{eq:definte} and $C=C(\fa, \fb, \delta)$ is uniform in $E$ with $\kappa(E)\leqslant N^{-\fb}$.
\eep
\begin{rmk}
	This result is a slight modification of \cite{gustavsson}*{Lemma 2.2}. We unfortunately cannot use  \cite{gustavsson}*{Lemma 2.2} directly, since it gives a non-explicit error of order $\mathcal{O}(1)$, which is larger than the leading term of order $N^{-\delta}$ for sub-microscopic intervals. 
\end{rmk}
\begin{proof}
To match the scaling of \cite{gustavsson}, we consider the measure $\widetilde{\operatorname{GUE}}_N=\sqrt{\frac{N}{2}}\operatorname{GUE}_N$, which is defined so that the largest eigenvalue of a  $\widetilde{\operatorname{GUE}}_N$ is approximately $\sqrt{2N}$.
Define
\beq
E_R=  E+\frac{\fa N^{-\delta}}{N\sqrt{{\kappa}(E)}}, \qquad I_1 =\left[
	\sqrt{\frac{N}{2}} E,\sqrt{\frac{N}{2}}  E_R \right].
\eeq
We first compute the expected number of eigenvalues in the interval $I_1$,
which is the right half of the rescaled version of the interval $I(\delta, \fa, E)$.
A similar argument applies to the left half.
We study here the right edge of the spectrum to use the same asymptotics as in \cite{gustavsson}, but the statement at the left edge follows by symmetry. 

If $\rho_N(x)=K_N(\sqrt{2N}x,\sqrt{2N}x)$ denotes the density of eigenvalues of the $\widetilde{\operatorname{GUE}}_N$ scaled to have its (limiting) support in $[-1,1]$, then we define 
\beq
g(E) = \int_{E/2}^{E_R/2}N\rho_N(x)\, \d x,
\eeq
which gives the number of eigenvalues in $I_1$.
Let $\Ai$ denote the Airy function. From \cite{ercolani}*{(4.4)} and \cite{ercolani}*{(4.21)} (see also the proof of \cite{gustavsson}*{Lemma 2.2}), we have the following asymptotic in $x$ in some fixed neighborhood $[ 0 , 1 + c]$ of $[0,1]$, where $c>0$ is a constant:
\begin{align}\label{eq:gustavexp}
N\rho_N(x)
&=
\left(
	\frac{\Phi'(x)}{4\Phi(x)}
	-
	\frac{\gamma'(x)}{\gamma(x)}
\right)
\left[
	2\Ai(\Phi(x))\Ai'(\Phi(x))
\right]\\
&+
\Phi'(x)\left[
	(\Ai'(\Phi(x)))^2
	-
	\Phi(x)(\Ai(\Phi(x)))^2
\right]
+
\O{\frac{1}{N\sqrt{\vert 1-x\vert}}},\label{e:thirdO}
\end{align}
where we define 
\beq\label{eq:gustav0}
\gamma(x)=\left(
	\frac{x-1}{x+1}
\right)^{1/4},\quad \Phi(x) = \left\{ \begin{array}{ll}
		-\left(
	3N\int_x^1\sqrt{1-y^2}\, \d y
\right)^{2/3} &\text{ for } x \leqslant 1,\\
		 \left(
	3N\int_1^x \sqrt{y^2 -1 }\, \d y
\right)^{2/3} &\text{ for } x> 1.
	\end{array} \right.
\eeq
The function $\gamma$ is defined by taking the limit from the upper half plane using the principal branch of the function $z \mapsto z^{1/4}$. For $x \in [1/2, 1]$, we observe that there exist constants $C,c>0$ such that
\beq\label{e:phiasym}
c N^{2/3}( x-1) \le \Phi(x) \le C N^{2/3} ( x-1).
\eeq
A similar bound holds for $x \in [ 1 , 3/2]$. 

By the proof of \cite{gustavsson}*{Lemma 2.2}, we have
\beq
\left(
	\frac{\Phi'(x)}{4\Phi(x)}
	-
	\frac{\gamma'(x)}{\gamma(x)}
\right)
\left[
	2\Ai(\Phi(x))\Ai'(\Phi(x))
\right]
=
\O{1}
\eeq
on $[0,1+c]$.
Then through integrating this bound we obtain
\beq\label{eq:gustav1}
\int_{E/2}^{E_R/2} \left(
	\frac{\Phi'(x)}{4\Phi(x)}
	-
	\frac{\gamma'(x)}{\gamma(x)}
\right)
\left(
	2\Ai(\Phi(x))\Ai'(\Phi(x))
\right) \d x = \O{E_R-E}=\O{\frac{ \fa N^{-\delta}}{N\sqrt{{\kappa}(E)}}}.
\eeq
For the second term in \eqref{eq:gustavexp},  we can again follow an exact computation
 given in the proof of \cite{gustavsson}*{Lemma 2.2} and find that
  \beq
\int_{E/2}^{E_R/2}
\Phi'(x)\left(
	(\Ai'(\Phi(x)))^2
	-
	\Phi(x)(\Ai(\Phi(x)))^2
\right)\d x
=
\frac{2}{3}\left(
	\Psi_1(E/2)-\Psi_1(E_R/2)
\right)
-\frac{1}{3}\left(\Psi_2(E/2)-\Psi_2(E_R/2)\right),
\label{e:b33}
\eeq
 where we defined
\beq
\Psi_1(x)=\Phi^2(x)\Ai^2(\Phi(x))-\Phi(x) \left[ \Ai'(\Phi(x)) \right]^2, \qquad
\Psi_2(x)=\Ai(\Phi(x))\Ai'(\Phi(x)).
\eeq
We observe that $\Ai(x)$ and $\Ai'(x)$ are smooth functions on $\R$. We may then suppose that $|1 -E| > N^{-2/3}/100$ and therefore using \eqref{e:phiasym} that $\Phi(E)> c$ and $\Phi(E_R) > c$ for some small $c$. Otherwise, $\Phi(E)< c$ and $\Phi(E_R) < c$, and the proposition follows from a Taylor expansion of $\Ai(x)$ and $\Ai'(x)$ about $0$.

In the case that $\Phi(E)$ and $\Phi(E_R)$ are distance at least $c$ from $0$, we may apply the Airy asymptotics given immediately before the proof of \cite{gustavsson}*{Lemma 2.1}, which give 
\beq
c N(1-x)^{3/2} \le \left| \Phi^2(x) \Ai^2(\Phi(x)) \right| \le C N(1-x)^{3/2}
\eeq
 as $x\to 1$ from below. A similar asymptotic holds as $x$ tends to $1$ from above. We obtain that
\beq\label{eq:gustav2}
\frac{2}{3}\left\vert
	\Psi_1(E)-\Psi_1(E_R)
\right\vert
\leqslant C N(E_R-E)\sqrt{\vert 1-E\vert}\leqslant CN^{-\delta}.
\eeq

Next, by definition of $\Psi_2$ we have that 
\begin{align}
\frac{1}{3}\left\vert
\Psi_2(E)-\Psi_2(E_R)
\right\vert 
\leqslant C(E_R-E)\Phi'(E)\left(
	\Ai'(\Phi(E))^2+\Ai(\Phi(E))\Ai''(\Phi(E))
\right)\\
=
	C(E_R-E)\Phi'(E)\left(
		\Ai'(\Phi(E))^2+\Phi(E)\Ai^2(\Phi(E))
	\right),
\end{align}
where we used the fact that $\Ai(x)\Ai''(x)=x\Ai^2(x).$ Using the Airy asymptotics from \cite{gustavsson}, we obtain that both $\Ai'(\Phi(E))^2$ and $\Phi(E)\Ai^2(\Phi(E))$ are bounded in absolute value by $C\sqrt{\Phi(E)}$, which gives
\beq
\frac{1}{3}\vert \Psi_2(E)-\Psi_2(E_R)\vert \leqslant C(E_R-E)\Phi'(E)\sqrt{\Phi(E)}.
\eeq
Finally, by definition of $\Phi(E)$ in \eqref{eq:gustav0}, we have that $\Phi'(E)=\frac{2N\sqrt{\vert 1-E^2\vert}}{\sqrt{\vert \Phi(E)}\vert}$ and thus 
\beq\label{eq:gustav3}
\frac{1}{3}\vert \Psi_2(E)-\Psi_2(E_R)\vert \leqslant
C(E_R-E)N\sqrt{\vert 1-E^2\vert}\leqslant C(E_R-E)N\sqrt{\vert 1-E\vert}\leqslant CN^{-\delta}.
\eeq
Finally, we remark that $(1-x)^{1/2}$ is integrable on intervals containing $1$, so that the third term in \eqref{e:thirdO} contributes $\mathcal O( N^{-1})$.
We obtain the  result by combining \eqref{eq:gustavexp} with \eqref{eq:gustav1}, \eqref{eq:gustav2} and \eqref{eq:gustav3}. The bounds can be made uniform in $E$ since, by assumption of the theorem, $E\in[1/2, 1 + c ]$ for $N > N_0(\fb)$.
\end{proof}

Combining the bound on the expected number of eigenvalues on a sub-microscopic interval from \Cref{prop:gustavsub} and \Cref{coro:lredet} we obtain a bound like \Cref{p:lregoe} for the GUE instead of the GOE. To obtain the bound on the GOE, we use the following relationship between these two ensembles from \cite{forrester2001inter}. 
This theorem states that to obtain the same distribution as eigenvalues of the GUE, one can take two independent matrices distributed as $\GOE_N$ and $\GOE_{N+1}$, put the eigenvalues of these matrices on the real line to obtain a distribution of $2N+1$ points, and extract every other eigenvalue.
\bet[\cite{forrester2001inter}*{(5.9)}]\label{theo:relationgoe}
	We have the following relation in distribution between independent eigenvalue distributions for $\GUE_N$, $\GOE_{N}$, and $\GOE_{N+1}$: 
	\beq
	\operatorname{GUE}_N=\mathrm{Even}(\GOE_N\cup\GOE_{N+1}).
	\eeq
\eet
We are now ready to prove \Cref{p:lregoe}.
\begin{proof}[Proof of \Cref{p:lregoe}]
Setting $\lambda=\delta \log N>0$ in \Cref{coro:lredet} and using the conclusion of  \Cref{prop:gustavsub}, we obtain the bound
\beq\label{eq:lregue}
\P_{\GUE_N}\left(
	\mathcal{N}_{\bm{\lambda}}(I (E))\geqslant k
\right)\leqslant
\exp\left({-\delta k\log N+(N^{\delta}-1)CN^{-\delta}}\right)
\leqslant CN^{-\delta k},
\eeq
 for some $C = C( \fa, \fb, \delta)$. Now, by \Cref{theo:relationgoe}, for any $k,N\in\N$ we have that
\beq\label{e:uppergoesquare}
\P_{\GUE_N}\left(
	\mathcal{N}_{\bm{\lambda}}(I (E))\geqslant k
\right)
\geqslant \P_{\GOE_N}
\left(	
	\mathcal{N}_{\bm{\lambda}}(I (E))\geqslant k 
\right)
\P_{\GOE_{N+1}}\left(\mathcal{N}_{\bm{\lambda}}(I(E))\geqslant k 
\right).
\eeq
There are two cases: either $\P_{\GOE_{N+1}}\left(\mathcal{N}_{\bm{\lambda}}(I(E))\geqslant k 
\right) \ge  \P_{\GOE_N}
\left(	
	\mathcal{N}_{\bm{\lambda}}(I (E))\geqslant k 
\right)$, or the reverse inequality holds. In the first case,
\beq
\P_{\GUE_N}\left(
	\mathcal{N}_{\bm{\lambda}}(I (E))\geqslant k
\right)
\geqslant \P_{\GOE_N}
\left(	
	\mathcal{N}_{\bm{\lambda}}(I (E))\geqslant k 
\right)^2.
\eeq
Then there exists $C = C(\fa, \fb, \delta)>0$ such that
\beq
\P_{\GOE_N}
\left(	
	\mathcal{N}_{\bm{\lambda}}(I_\delta(E))\geqslant k
\right)
\leqslant C N^{-\delta k/2 }.
\eeq
In the second case,
\beq
\P_{\GUE_N}\left(
	\mathcal{N}_{\bm{\lambda}}(I (E))\geqslant k
\right)
\geqslant \P_{\GOE_{N+1}}
\left(	
	\mathcal{N}_{\bm{\lambda}}(I (E))\geqslant k 
\right)^2,
\eeq
and
\beq
\P_{\GOE_N}
\left(	
	\mathcal{N}_{\bm{\lambda}}(I(E))\geqslant k
\right)
\leqslant C (N-1)^{-\delta k /2} \le C N^{-\delta k/2},
\eeq
where we modified the value of $C$. Combining these cases completes the proof.
\end{proof}


Having established the desired level repulsion estimate for the GOE, we now turn to the case of Gaussian divisible matrices. We consider another matrix dynamics which consists of the same stochastic differential equation as \eqref{eq:defdyn}, but with a different initial condition and matrix Brownian motion $\tilde B_s$. Let $\tilde{H}_0$ be a random symmetric matrix sampled from the $\GOE_N$, and denote by $\nu_1(s)\leqslant \dots\nu_N(s)$ the eigenvalues of $\tilde{H}_s$, defined by the dynamics
\beq\label{e:oudynamics}
\tilde{H}_s = \frac{1}{\sqrt{N}}\, \d \tilde B_s - \frac{1}{2}\tilde{H}_s.
\eeq
We then have the following theorem from \cite{bourgade2018extreme}. The coupling in the theorem is given by choosing the Brownian motion $\tilde B_s$ so that the Brownian motions driving the eigenvalue dynamics \eqref{eq:dysonval} for $\tilde H_s$ are the same as those driving the dynamics for $H_s$.
\bet[\cite{bourgade2018extreme}*{Theorem 2.8}] \label{t:bourgadeloc}
There exists a coupling of $H_s$ from \eqref{eq:defdyn} and $\tilde H_s$ from \eqref{e:oudynamics} such that the following holds.
 For any $\varepsilon > 0 $ there exist $C_1(\eps), c_1(\eps), c_2(\eps) > 0$ such that 
 \beq
 \P\left(
 	\vert \lambda_k(s)-\nu_k(s)\vert
 	\geqslant \frac{N^\varepsilon}{Ns}
 	\text{ for all }k\in\unn{1}{N}\text{ and }s\in[0,1]
 \right)
 \leqslant C_1\exp\left( {-c_1(\log N)^{c_2\log\log N}}\right).
 \eeq
\eet
\begin{rmk}
This theorem was proved with a weaker probability bound in \cite{bourgade2013eigenvector} for the case where the matrix entry distributions have finite moments. The proof relies on local law and rigidity estimates, for example \cite{bourgade2013eigenvector}*{Lemma 2.3}, which are the source of this probability loss. However, our assumption that the entry distributions are uniformly subexponential allows the stronger estimates given in \Cref{l:goodset}, which in turn allow us to state \cite{bourgade2018extreme}*{Theorem 2.8} with exponentially high probability bounds, as above.
\end{rmk}

Since $\GOE_N$ is invariant for the dynamics \eqref{eq:defdyn}, for all times $s\in[0,1]$, the $\nu_k(s)$ are distributed as eigenvalues of a $\GOE_N$. We now apply \Cref{t:bourgadeloc} to extend our level repulsion estimate at the edge from the $\GOE_N$ to Gaussian divisible ensembles.

\bel\label{l:lrgaussdiv}
Let $H_0$ be a generalized Wigner matrix and $H_s$ be the dynamics defined in \eqref{eq:defdyn}. Fix $k\in\mathbb{N}$ and  $\fa, \delta>0$. 
 Then there exists a constant $C=C( \fa, \delta)$  and an event $\A = \mathcal A(\delta)$ such that for any $E$ satisfying $\kappa(E) \le N^{- 100 \delta}$ and $s \in [ N^{- 40 \delta} , 1]$,
	\beq
	\P \left(
		\one_{\mathcal A }\mathcal{N}_{\bm{\lambda} (s) }(I )\geqslant k
	\right)
	\leqslant  C N^{-\delta k /2},
	\qquad \P(\mathcal A^c) \le C_1\exp\left( {-c_1(\log N)^{c_2\log\log N}}\right),
	\eeq
	where  $I = I(\delta, \fa, E)$ and $C_1(\delta), c_1(\delta), c_2(\delta) > 0$ are constants. 
\eel
\begin{proof}
 For any $s \in [ 0 , 1]$, $\eps>0$, and $\ell \in\unn{1}{N}$, \Cref{t:bourgadeloc} gives 
\beq\label{e:coupling}
| \lambda_\ell(s)  - \nu_\ell(s) | \le \frac{N^\eps}{Ns}
\eeq
uniformly in $\ell$ on the event $\mathcal A = \mathcal A(\eps)$ of that theorem.

Set $s_0 = N^{3\delta/2} \sqrt{\kappa(E)}$. Then the hypotheses of the lemma imply $N^{-1} \ll s_0 \ll 1$. For $s\in [s_0 , 1]$ and $\eps = \delta/8$, we have that 
\beq\label{e:intervalsize2}
\frac{N^\eps}{Ns } \le \frac{\fa N^{-\delta} }{ N \sqrt{\kappa(E)}}
\eeq
on $\mathcal A$, for $N \ge N_0(\fa, \delta)$. 
We deduce from \eqref{e:coupling} and \eqref{e:intervalsize2} that if $\lambda_\ell (s)  \in I(\delta, \fa, E)$ and $s \in [s_0, 1]$, then $\nu_\ell(s) \in I(\delta, 2 \fa, E)$. It follows that

\beq
\P\left( \one_{\mathcal A}
	\mathcal{N}_{\bm{\lambda}(s)}(I(\delta, \fa, E)) 
	\geqslant k
\right)
\leqslant 
\P\left(
	\mathcal{N}_{\bm{\nu}(s)}(I(\delta, 2\fa, E) )
	\geqslant k
\right)
\eeq
when $N$ is large enough.
Since the eigenvalues $\bm{\nu}_s$ are distributed as the eigenvalues of $\GOE_N$, applying \Cref{p:lregoe} completes the proof. 
\end{proof}

\subsubsection{Uniform estimate}
Observe that the following estimate is effective only for $k$ large.
\bep\label{p:goecombined}
Let $H_0$ be a generalized Wigner matrix and $H_s$ be the dynamics defined in \eqref{e:oudynamics}.
Fix $k\in\mathbb{N}$ and $\delta, \fa >0$.
Suppose $ E \in [ - 2 -  N^{-100\delta} , 2 + N^{- 100\delta} ]$ and $s \in [ N^{- 40\delta} , 1]$. Then there exists a constant $C=C( \fa, \delta)$  and event $\mathcal A(\delta)$ such that
	\beq
	\P \left(
		\one_{\A} \mathcal{N}_{\bm{\lambda}}(I )\geqslant k
	\right)
	\leqslant  C N^{200 \delta } N^{ -  \delta k/2 }  \qquad \P(\mathcal A^c) \le C_1\exp\left( {-c_1(\log N)^{c_2\log\log N}}\right)
	\eeq
	where  $I = I(\delta, \fa, E)$ and $C_1(\delta), c_1(\delta), c_2(\delta) > 0$ are constants. 
	
\eep	

\begin{proof}
We apply \Cref{l:deformedpacking} with $\sigma =  200 \delta$ and \Cref{l:lrgaussdiv}.
This yields 
\beq
\P \left(
		\mathcal{N}_{\bm{\lambda}}(I )\geqslant k
	\right) \le N^{200 \delta } N^{ -  \delta k/2 } + C N^{-\delta k /2 } \le C N^{200 \delta } N^{ -  \delta k/2 } , 
	\eeq
	from which the claim follows. 
\end{proof}

\subsection{Level repulsion for generalized Wigner matrices}\label{s:lrwigner}

Before proceeding to the following proof, we first make some definitions and present some preliminary computations.

For $n\in \N$, $E \in R$, and $\delta>0$, set $I_n = I(\delta, n \fa, E)$. 
Let $q = q_N \in C^\infty(\R)$ be a function such that $q(x) = 1$ for $x \in I_2$, $q(x)=0$ for $x\in I_3^c$, $|q(x)| \le 1$ for all $x\in \R$,
and $|q^{(d)}(x)|  \le C N^{d(1+ \delta )} \kappa(E)^{d/2}$ for $d \in \unn{1}{5}$ and some constant $C= C(\fa)>0$. 
Let $\chi  \in C^\infty(\R)$ be a weakly increasing function such that $\chi (x) = 0$ for $x \le 450 $, $\chi(x)=1$ for $x \ge 500 $, $|\chi(x)| \le 1$ for all $x\in \R$,
and $|\chi^{(d)}(x)|  \le C$ for $d \in \unn{1}{5}$ and some constant $C>0$. 

Let $\omega >0$ be a parameter to be chosen later. We define the smoothed eigenvalue counting function $F(M)$ on $N\times N$ symmetric matrices by 
\beq
g(M) = \sum_{j: | \gamma_j -  E|  \le N^{-1 + \omega} \kappa(E)^{-1/2}} q\left( \widetilde \mu_j\left (M \right) \right), \qquad F(M) =\chi\left( g(M) \right) g(M),
\eeq
where $\widetilde \mu_j(M)$ are the regularized eigenvalues of $M$ defined in \Cref{p:regulareigval} with parameters $\eps, \delta_1 >0$ chosen such that 
\beq \label{e:nuparamchoice}
\omega = \eps =  \frac{\delta_1}{2}.\eeq
 We also define $F_m(M) = F(M)^m$ for $M \in \matn$. 

Recall $\hat j = \min(j, N - j +1)$. For any random matrix $M$ taking values in $\matn$, with corresponding eigenvalues $\{\mu_j\}_{j=1}^N$, we let $\mathcal A_1(M)$ be the event where 
\beq
\sup_{j\in \unn{1}{n}}
\left\vert
	\mu_j -\gamma_j
\right\vert
\leqslant
N^{-2/3+\omega/2}{\hat j }^{-1/3}.
\eeq
We let $\mathcal A_2(M)$ be the event $\mathcal A(\delta_1, \eps)$ from \Cref{p:regulareigval}, and set $\mathcal A = \mathcal A(M) = \mathcal A_1(M) \cap \mathcal A_2(M)$.

\bel
For any $M \in \matn$, we have
\beq\label{e:Fbdd}
\one_{\mathcal A} \left| \partial_{ab}^d F(M)  \right |\le  C N^{7(\omega + \delta + \delta_1 + \eps)} ,\qquad  \left| \partial_{ab}^d F(M)  \right| \le C N^{2\omega + C}
\eeq
for $d \in \unn{1}{5}$.
\eel
\bp
On $\mathcal A$, the $\tilde \mu_j$ satisfy $|\partial_{ab}^d \widetilde \mu_j | \le C N^{-2/3  + d(\delta_1 + \eps)} \hat{j}^{-1/3}$ for $d \le 5$ and a constant $C (\eps, \delta_1)>0$.
Therefore, for $d\le 5$, we have
\begin{align}\label{e:derivativeF}
\one_{\mathcal A} \left| \partial_{ab}^d g(M)  \right|
& \le  C \sum_{j: | \gamma_j -  E|  \le N^{-1 + \omega} \kappa(E)^{-1/2} } \left| \partial_{ab}^d q\left( \widetilde \lambda_j \right) \right|\\ &\le \sum_{m=1}^5 C N^{2\omega + m/3 + d(\delta_1 + \eps) +  d\delta } \kappa(E)^{m/2} \hat{j}^{-m/3} \le C N^{7(\omega + \delta + \delta_1 + \eps)}.
\end{align}
Here, we used the fact that there are at most $2 N^{2 \omega}$ indices $j$ such that $| \gamma_j -  E|  \le N^{-1 + \omega} \kappa(E)^{-1/2}$ on $\mathcal A$.
We also used the fact that, for these $j$, we have $(N / \hat{j})^{1/3} \le C N^{\omega} \kappa(E)^{-1/2}$.
Additionally, we have the trivial bound
\beq\label{e:derivativeFtrivial}
 \left| \partial_{ab}^k g(M)  \right|
 \le C \sum_{j: | \gamma_j -  E|  \le N^{\omega - 1} } \left| \partial_{ab}^d q\left( \widetilde \lambda_j \right) \right| \le C N^{2\omega + Cd}.
\eeq
From \eqref{e:derivativeF} and \eqref{e:derivativeFtrivial}, we deduce the conclusion from the definition of $F(M)$ and the fact that $\chi$ is bounded with bounded derivatives.
\ep

Finally, we note that when $H$ is a generalized Wigner matrix,
\beq\label{e:simulA}
\mathcal B(\q, \omega, \delta_1, \eps, \delta_2, \eps_2) \subset \mathcal A\left(\Theta^{(a,b)}_w H \right)
\eeq
for any $w \in [0,1]$ and $a,b \in \unn{1}{N}$, where $\q, \delta_2, \eps_2$ may be chosen arbitrarily.
We conclude by \Cref{l:good1} and \eqref{e:nuparamchoice} that

\beq\label{e:Aexphigh}
\P \left( \bigcap_{w\in [ 0,1]} \bigcap_{ a,b \in \unn{1}{N} }\A\left( \Theta^{(a,b)}_w H \right)  \right)\ge1 - C_1 \exp\left({-c_1(\log N)^{c_1 \log\log N}}\right)
 \eeq
for some constants $c_1 = c_1(\omega)>0$, $C_1 = C_1(\omega) >0$.

\bel\label{lem:lrewigner} Let $H$ be a generalized Wigner matrix, and let $\bm{\lambda}$ be the vector of its eigenvalues. Fix $k\in \N$ such that $k\ge 500$, $\delta \in (0,1/100)$ and $\fa>0$. Then there exist a constant $C = C( \fa , \delta) > 0$ and an event $\mathcal F = \mathcal F(\delta)$ such that for any $E \in [ -2 - N^{-200\delta}, 2 + N^{-200\delta}]$,

	\beq
	\P \left(
		\one_{\mathcal F} \mathcal{N}_{\bm{\lambda}  }(I )\geqslant k
	\right)
	\leqslant  C N^{-\delta \log k },
	\qquad \P(\mathcal F^c) \le C_1\exp\left( {-c_1(\log N)^{c_2\log\log N}}\right),
	\eeq
	where  $I = I(\delta, \fa, E)$ and $C_1(\delta), c_1(\delta), c_2(\delta) > 0$ are constants. 

\eel

\begin{proof} 

Recall the dynamics $H_s$ defined in \eqref{e:oudynamics} and set $s_0 = N^{-40\delta}$, in preparation for the use of \Cref{p:goecombined}. By \cite{erdos2017dynamical}*{Lemma 16.2} there exists a generalized Wigner matrix $H_0$ such that the matrix $R= H_{s_0}$ satisfies $\E[h_{ij}^k] = \E[r_{ij}^k]$ for $k\in\unn{1}{3}$ and $\left| \E[h_{ij}^4 ]  -  \E[r_{ij}^4 ] \right| \le C N^{-2} s_0$ for some constant $C>0$ depending only on the constants used to verify \Cref{d:wigner} holds for $H_0$. We observe that $R$ is also a generalized Wigner matrix, and in particular it has subexponential entries.

Let $\bm{\nu} = (\nu_1, \nu_2, \dots, \nu_N)$ denote the vector of eigenvalues of $R$ arranged in increasing order. 
By the definition of $F(M)$ and \eqref{e:regev}, 
\beq \label{e:lrintervaldom}
\one_{\mathcal A(R)}\one_{\{\mathcal   N_{\bm{\nu}}(I_4 ) \ge 450 \}} \mathcal N_{\bm{\nu}}( I_4 ) \ge\one_{\mathcal A(R)} \one_{\{\mathcal   N_{\bm{\nu}}(I_4 ) \ge 450 \}} \mathcal N_{\tilde {\bm{\nu}} }( I_3 ) \ge \one_{\mathcal A(R)} F( R )
\eeq
where we recall that $I_n = I(\delta, n\fa, E)$. We compute using \Cref{l:lrgaussdiv} and \eqref{e:lrintervaldom} that
\begin{align}\label{e:fhmoments}
\E \left[\one_{\mathcal A(R)}F(R)^{\lceil \delta \log N \rceil} \right] &\le   
\E \left[\one_{\mathcal A(R)} \one_{\mathcal  \mathcal N_{\bm{\nu}}(I_4 ) \ge 450 }  \mathcal N_{\bm{\nu}}(I_4 )^{\lceil \delta \log N \rceil} \right]  \\
&\le \sum_{k=450}^\infty k^{\lceil \delta \log N \rceil} \P\left( \mathcal N_{\bm{\nu}}(I_4 ) \ge k \right) 
 \le  C \sum_{k=450}^\infty k^{\lceil \delta \log N \rceil}  C N^{200 \delta } N^{ -  \delta k/2 } 
 \le C_1,
\end{align}
for some constant $C_1=C_1(\delta, \fa ) > 1$. Since $\mathcal A(R)$ holds with exponentially high probability by \eqref{e:Aexphigh}, and $F(R) \le N$, we see $\E \left[\one_{\mathcal A^c (R) }F(R)^{\lceil \delta \log N \rceil} \right] \le C_1$ and therefore $\E \left[F(R)^{\lceil \delta \log N \rceil} \right] \le C_1$, after increasing the value of $C_1$. 

Fix any bijection
\beq
\phi \colon \{ ( i, j) : 1 \le i \le j \le N\} \rightarrow \unn{1}{\gamma_N},
\eeq
where $\gamma_N = N ( N + 1) /2$, and define the matrices $H^1, H^2, \dots , H^{\gamma_N}$ by 
\beq
h_{ij}^\gamma = 
\begin{cases}
h_{ij} & \text{if } \phi(i,j) \leq \gamma
\\
r_{ij} & \text{if } \phi(i,j) > \gamma
\end{cases}
\eeq
for $i \le j$. 

Fix some $\gamma \in \unn{1}{\gamma_N}$ and consider the indices $(i,j)$ such that $\phi(i,j) = \gamma$. For any $m \ge 1$, we may Taylor expand $F_m \left(H^\gamma\right)$ in the $(i,j)$ entry, write $\partial = \partial_{ij}$, and find
\begin{align}
\label{e:taylora1} F_m \left( H^\gamma \right) - F_m \left(  \Theta^{(i,j)}_0 H^\gamma \right) &= \partial F_m \left(  \Theta^{(i,j)}_0 H^\gamma \right) h_{ij} + \frac{1}{2!} \partial^2 F_m \left(  \Theta^{(i,j)}_0 H^\gamma \right) h_{ij}^2 
+ \frac{1}{3!}\partial^3 F_m \left(  \Theta^{(i,j)}_0 H^\gamma \right) h_{ij}^3 \\
& \label{e:taylora2} + \frac{1}{4!}\partial^4 F_m \left(  \Theta^{(i,j)}_0 H^\gamma \right) h_{ij}^4 + \frac{1}{5!}\partial^5 F_m \left(  \Theta^{(i,j)}_{w_1(\gamma)} H^\gamma \right) h_{ij}^5,
\end{align}
where $w_1(\gamma) \in [0,1]$ is a random variable depending on $h_{ij}$. Similarly, we can expand $F_m \left( H^{\gamma -1} \right)$ in the $(i,j)$ entry, and after subtracting this expansion from \eqref{e:taylora1} and \eqref{e:taylora2} and taking expectation, we find
\begin{align}
 \label{e:fourthorder} \E \left[ F_m \left( H^{\gamma } \right)  \right]- \E \left[ F_m \left(   H^{\gamma -1} \right) \right] &= \frac{1}{4!} \E \left[ \partial^4 F_m \left(  \Theta^{(i,j)}_0 H^\gamma \right) h_{ij}^4 \right] - \frac{1}{4!}\E \left[ \partial^4 F_m \left(  \Theta^{(i,j)}_0 H^\gamma \right) r_{ij}^4 \right] \\
  \label{e:fifthorder} &+ \frac{1}{5!} \E \left[ \partial^5 F_m \left(  \Theta^{(i,j)}_{w_1(\gamma)} H^\gamma \right) h_{ij}^5\right]
  - \frac{1}{5!} \E \left[ \partial^5 F_m \left(  \Theta^{(i,j)}_{w_2(\gamma)} H^\gamma \right) r_{ij}^5 \right].\end{align}
Here $w_2(\gamma) \in [0,1]$ is a random variable depending on $r_{ij}$, and we used that $\E [ h^k_{ij} ] = \E [ r^k_{ij} ]$ for $k \in \unn{1}{3}$. 

We now use this expansion to show that $\E F(H)^{\lceil \delta \log N \rceil} \le C$ for a constant $C$. Our argument proceeds by induction, with the induction hypothesis at step $m \in \N$ being that 
\beq\label{e:inductionhypo}
 \E \left[ F_n \left(  \Theta^{(a,b)}_{\kappa } H^\gamma \right)\right]  \le K_n
\eeq
holds for a constant $K_n$ depending on $n$, for all $n \le m \le \lceil \delta \log N \rceil$ and choices of  $\kappa \in [ 0 ,1 ]$ and $(a,b)\in\unn{1}{N}^2$. We may assume, by increasing the constants $K_n$ if necessary, that $K_n\ge 1$ and $K_n$ is increasing in $n$. We will fix $K_n$ later.

 The base case $m = 0$ is trivial. Assuming the induction hypothesis holds for $m -1$, we will derive a bound for $m$. Using the independence of $h_{ij}$ and $r_{ij}$ from $ \Theta^{(i,j)}_0 H^\gamma$, we may rewrite the terms on the right side of \eqref{e:fourthorder} as 
\beq\label{e:b72}
\E \left[ \partial^4 F_m \left(  \Theta^{(i,j)}_0 H^\gamma \right) h_{ij}^4 \right] - \E \left[ \partial^4 F_m \left(  \Theta^{(i,j)}_0 H^\gamma \right) r_{ij}^4 \right] =  \E \left[ \partial^4 F_m \left(  \Theta^{(i,j)}_0 H^\gamma \right)  \right]  \E \left[  h_{ij}^4 - r_{ij}^4  \right].
\eeq
For the second factor, we recall that $\left|  \E \left[ h_{ij}^4\right] - \E \left[ r_{ij}^4 \right]  \right|  \le C N^{-2}s_0 =  C N^{-2 - 40\delta}$. For the first, we compute 
\begin{align}\label{e:4thcompute}
\partial^4 F_m  = \partial^4 \left( F^m \right) &= m F_{m-1} F^{(4)}  + 3 m (m - 1 ) F_{m-2} (F^{(2)})^2 + m (m-1) (m-2) (m-3) F_{m-4} (F')^4 \\
& + 4 m (m-1) F_{m-2} F^{(1)} F^{(3)} + 6 m (m-1) (m-2) F_{m-3} (F')^2 F^{(2)}.
\end{align}
Using the induction hypothesis \eqref{e:inductionhypo} for $ n \le m -1$, $m \le \lceil \delta \log N \rceil$, the fact that $F_m \ge 0$, and the first inequality of \eqref{e:Fbdd}, we find from \eqref{e:4thcompute} that 
\beq\label{e:4analogy}
\left| \one_{\mathcal A(H^\gamma)} \E \left[ \partial^4 F_m \left(  \Theta^{(i,j)}_0 H^\gamma \right)  \right]  \right| \le C K_{m-1} (\log N)^4 N^{28(\omega + \delta + \delta_1 + \eps)}.
\eeq
Further, by the second inequality in \eqref{e:Fbdd}, and because $\mathcal A^c(H^\gamma)$ holds with exponentially high probability  by \eqref{e:Aexphigh}, we find\footnote{We remark that constants in the probability bound \eqref{e:Aexphigh} do not depend on the choice of $\gamma$, since the $H^\gamma$ verify \Cref{d:wigner} simultaneously for the appropriate choice of constants. Therefore, the $C$ in \eqref{e:4analogytrivial} is uniform in $\gamma$.}
\beq\label{e:4analogytrivial}
\left| \one_{\mathcal A^c(H^\gamma)} \E \left[ \partial^4 F_m \left(  \Theta^{(i,j)}_0 H^\gamma \right)  \right]  \right| \le C N^{-2- 10\delta}.
\eeq

 It follows from \eqref{e:nuparamchoice}, \eqref{e:b72}, \eqref{e:4analogy}, and \eqref{e:4analogytrivial} that if $\omega$ is chosen small enough relative to $\delta$, so that $28(\omega + \delta + \delta_1 + \eps)< 29\delta $, then there exists some $C = C(\delta, \fa)$ such that the bound
\beq\label{e:4final}
\eqref{e:fourthorder} \le  C K_{m-1} N^{-2-10 \delta} .
\eeq
holds for the fourth order terms for all $m \le \log N$. 

For the terms in \eqref{e:fifthorder}, we first observe that $\sup_{i,j} | r_{ij} | + | h_{ij} |\le C N^{-1/2 + \delta}$ on a set $\mathcal C$ of exponentially high probability. Therefore, we find similarly to our computation of \eqref{e:4final} that
\beq
\E \left[ \partial^5 F_m \left(  \Theta^{(i,j)}_{w_1(\gamma)} H^\gamma \right) h_{ij}^5\right] \le
 C N^{-5/2 + 5 \delta} \left( \E \left[ \left| \partial^5 F_m \left(  \Theta^{(i,j)}_{w_1(\gamma)} H^\gamma \right) \right| \right]  + 1\right),
\eeq
and likewise for the second term in \eqref{e:fifthorder}. Then we obtain, analogously to the bound for the fourth order term in \eqref{e:4analogy}, and using that $\delta < 1/100$,
\beq\label{e:5final}
\eqref{e:fifthorder} \le  C K_{m-1} N^{-2- 10\delta}.
\eeq
Therefore
\beq
\left| \E \left[ F_m \left( H^{\gamma } \right)  \right]- \E \left[ F_m \left(   H^{\gamma -1} \right) \right] \right| \le C K_{m-1} N^{-2-10 \delta},
\eeq
and summing over all $\mathcal O(N^2)$ pairs $(i,j)$, we find
\beq\label{e:4momentconclude}
\left| \E \left[ F_m \left( R \right)  \right]- \E \left[ F_m \left( H^\gamma \right) \right] \right| \le C K_{m-1} N^{- 10 \delta}.
\eeq
for any $\gamma$.

By \eqref{e:fhmoments}, $\E \left[  F_m(R) \right] \le C_1$ for all $m \le \log N$. Together with \eqref{e:4momentconclude}, we deduce that 
\beq
\E\left[ F_m\left( H^\gamma \right) \right] \le C K_{m-1} N^{-10\delta} + C_1.
\eeq
 This bounds the quantity in \eqref{e:inductionhypo} when $\kappa = 1$.

To address other values of $\kappa$, we consider the following expansion:
\begin{align}
\label{e:taylorc1} F_m \left( H^\gamma \right) - F_m \left(  \Theta^{(a,b)}_\kappa H^\gamma \right) &= \partial F_m \left(  \Theta^{(a,b)}_0 H^\gamma \right) h_{ij} + \frac{1}{2!} \partial^2 F_m \left(  \Theta^{(a,b)}_0 H^\gamma \right) h_{ij}^2 
+ \frac{1}{3!}\partial^3 F_m \left(  \Theta^{(a,b)}_0 H^\gamma \right) h_{ij}^3 \\
& \label{e:taylorc2} + \frac{1}{4!}\partial^4 F_m \left(  \Theta^{(a,b)}_0 H^\gamma \right) h_{ij}^4 + \frac{1}{5!}\partial^5 F_m \left(  \Theta^{(a,b)}_{w(\kappa)} H^\gamma \right) h_{ij}^5,
\end{align}
Here $w(\kappa) \in [ 0, 1]$ is a random variable. The same argument that gave the bound \eqref{e:4final} shows that the right side of \eqref{e:taylorc1} and the first term of \eqref{e:taylorc2} may be bounded in absolute value by $C K_{m-1} N^{-2 - 10\delta}$. The second term of \eqref{e:taylorc2} is also bounded by $CK_{m-1} N^{-2-10\delta}$ by the reasoning leading to \eqref{e:5final}.  We  conclude
\beq
\sup_{\kappa \in [0,1]} \sup_{a,b \in \unn{1}{n}} \E \left[ F_m \left(  \Theta^{(a,b)}_{\kappa } H^\gamma \right)\right]  \le C_2 K_{m-1} N^{-10\delta} + C_1,
\eeq
for some constant $C_2 = C_2(\delta, \fa)$.

We may therefore take
\beq
K_m = C_2 K_{m-1} N^{-10\delta} + C_1,  \qquad K_0 = C_1
\eeq
where we recall that we assumed $C_1 > 1$. For $N\ge N_0(\delta, \fa)$, we have 
\beq
K_m \le \frac{K_{m-1}}{2 } + C_1,
\eeq 
which implies
$K_m \le 2 C_2$
for all $m \le \lceil \delta \log N \rceil$ when $N\ge N_0(\delta, \fa)$. This implies the existence of a constant $C_3$ such that $K_m \le C_3$
for all $m \le \lceil \delta \log N \rceil$ for all $N \in \N$.

The conclusion of our induction argument is that $\E F(H )^{\lceil \delta \log N \rceil} \le C$ for a constant $C=C(\delta, \fa)$ independent of $N$. This implies, using Markov's inequality, that
\beq
\P( \one_{\A} \mathcal N_{\tilde \lambda }( I_2) \ge k) \le \P( \one_{\A} F(H ) \ge k) \le \frac{\E F(R)^{\lceil \delta \log N\rceil} }{ k^{\lceil \delta \log N \rceil }} \le \frac{C}{N^{\delta \log k}}
\eeq
for all $k\in \mathbb N$ such that $k\ge 500$.
Finally, using that \eqref{e:regev} holds on $\mathcal A(H)$, we have
\begin{equation}
\P( \one_{\mathcal A} N_{\bm{\lambda} } (I_1 ) \ge k) \le \P(\one_{\mathcal A} \mathcal N_{\tilde {\bm \lambda}}( I_2) \ge k)  \le \frac{C}{N^{\delta \log k}}.
\end{equation}
This completes the proof after setting $\mathcal F = \mathcal A(H)$.
\end{proof}

Given the previous lemma, the proof of the next proposition is similar to that of \cite{landon2017convergence}*{Theorem 3.6}. 
\bep\label{p:lr}
 Let $H$ be a generalized Wigner matrix, and let $\bm{\lambda}$ be the vector of its eigenvalues. Fix $k\in \N$ such that $k \ge 500$, $\fa>0$, and $\delta \in (0, 1/100)$. Then there exist a constant $C = C( \fa , \delta) > 0$ and an event $\mathcal A = \mathcal A(\delta)$ such that for any $i \in \unn{1}{N} $,

	\beq
	\P \left(
		\one_{\mathcal A} \mathcal{N}_{\bm{\lambda}  }(I(\lambda_i)  )\geqslant k
	\right)
	\leqslant  C N^{\delta/5  + \delta( 1-  \log k)/2 },
	\qquad \P(\mathcal A^c) \le C_1\exp\left( {-c_1(\log N)^{c_2\log\log N}}\right),
	\eeq
	where  $I(\lambda_i) = I(\delta, \fa, \lambda_i)$ and $C_1(\delta), c_1(\delta), c_2(\delta) > 0$ are constants. 

\eep

\begin{proof}[Proof of \Cref{p:lr}]
Let $\mathcal A_1 = \mathcal A_1 ( \delta) $ be the set from \Cref{l:goodset}, and let $\mathcal A_2 = \mathcal A_2(\delta)$ be the set from \Cref{lem:lrewigner}. From these lemmas, we know that $\mathcal A = \mathcal A_1 \cap \mathcal A_2$ holds with exponentially high probability in the sense of \eqref{e:exphigh}.

Observe that  the event $\{ N(I(\lambda_i)  )\geqslant k  \}$ satisfies 
\beq
\{ \mathcal N( I(\lambda_i)  )\geqslant k  \} \subset \{ N( I(\lambda_i)  )\geqslant k, |\lambda_i - \gamma_i| \le N^{-2/3 + \delta/10} {\hat i}^{-1/3} \} \cup \mathcal A^c.
\eeq
Set 
\beq I_j = \left[ \gamma_i + (j -2) \frac{\fa N^{-\delta} }{ N \sqrt{\kappa( \gamma_i )}}, \gamma_i + (j +2) \frac{\fa N^{-\delta} }{ N \sqrt{\kappa(\gamma_i )}}\right].
\eeq
Note that 
\beq
\{ \mathcal N (I(\lambda_i)  )\geqslant k, |\lambda_i - \gamma_i| \le N^{-2/3 + \delta/10} \hat i \} \subset \bigcup_{|j| <  \frac{2 N^{\delta/5} }{\fa} } \{ \mathcal N(I_j) \ge k \}.
\eeq
We used here that $\sqrt{\kappa(\gamma_i) } \le 2 ({\hat i} / N)^{1/3}$. 
By \Cref{lem:lrewigner}, 
\beq
\P\left( \one_{\mathcal A} \bigcup_{|j| < 2 N^{\delta/5} /\fa } \{ \mathcal N(I_j) \ge k \} \right) \le \sum_{{|j| < 2 N^{\delta/5} /\fa } }  \P\left(  \one_{\mathcal A} \mathcal N(I_j) \ge k \right) \le C(\fa, \delta)   N^{\delta/5  + \delta (1 -  \log k)/2} . 
\eeq
In the last inequality, we used the fact that $I_j$ may be slightly larger than $I(E)$ for $E = \gamma_i + j \fa N^{-\delta-1} \kappa(\gamma_i)^{-1/2}$, so \Cref{lem:lrewigner}, must be applied with $\delta' = \delta/2$. This finishes the proof after recalling $\mathcal A$ holds with exponentially high probability.
\end{proof}

\bp[Proof of \Cref{p:lrmain}]
This is an immediate consequence of \Cref{p:lr} and \Cref{p:regulareigval}.
\ep

\subsubsection{The case $k=2$}
%

The previous estimate control the probability of having $k$ eigenvalues in a sub-microscopic interval for large $k$. For $k=2$, we have the following more precise estimate following from gap universality.
\bep\label{prop:lrek2}
	Let $H$ be a generalized Wigner matrix. There exists $\delta_0$ such that the following holds. Let $E \in [-2 - N^{-\delta_0}, 2 + N^{-\delta_0}]$. If $E < \gamma_1$ take $i=1$; if $i > \gamma_N$, take $i = N$. Otherwise, let $i\in\unn{1}{N-1}$ be such that $\gamma_i\leqslant E <  \le \gamma_{i+1}$. For any $0<\delta<\delta_0$, there exist $\alpha = \alpha(\delta)>0$ and $C=C(\delta)>0$ such that 
	\beq\label{e:lr2}
	\P{\left(\mathcal{N}(E-N^{-2/3-\delta}\hat{i}^{-1/3},E+N^{-2/3-\delta}\hat{i}^{-1/3}\geqslant 2\right)}\leqslant CN^{-\alpha-\delta}.
	\eeq
\eep
\begin{proof}
First, for $i\in\unn{\alpha N^{1-c}}{N-\alpha N^{1-c}}$, where $c>0$ is a small constant, the proof was given in the discussion following \cite{erdos2012gap}*{(6.32)}; it combines gap universality with the level repulsion induced by Gaudin distribution for Gaussian ensembles.\footnote{The proof was written in the bulk $(-2+\kappa,2-\kappa)$ corresponding to $i\in\unn{\alpha N}{(1-\alpha)N}$ for any $\alpha>0$, but it directly translates to an estimate on this slightly extended bulk for $c>0$ small enough.} For $i\in\unn{1}{N^{1/4}}\cup\unn{N-N^{1/4}}{N}$, level repulsion follows in a similar way from gap universality at the edge, which was proved in \cite{bourgade2014edge}*{Theorem 2.7}.\footnote{One must first establish level repulsion for the GOE and GUE for this argument to go through. See \cite[Remark 1.5]{knowles2013eigenvector} for a sketch of the proof of this fact.}
 It therefore remains to handle the intermediate regime $i\in\unn{N^{1/4}}{\alpha N^{1-c}}\cup\unn{N-\alpha N^{1-c}}{N-N^{1/4}}$. Gap universality here is a consequence of the proof of \cite{bourgade2018extreme}*{Theorem 1.6}, as noted after the theorem statement there. For completeness, we describe how to convert this universality result to a level repulsion estimate.

Denote $I_i(E)= [E-N^{-2/3-\delta}\hat{i}^{-1/3},E+N^{-2/3-\delta}\hat{i}^{-1/3}]$. Let $\omega>0$ be a parameter chosen later and denote $\mathcal C = \mathcal{C}(\omega)$ the event where rigidity \eqref{e:rigidity} holds with parameter $\omega$. For any $D>0$, there exists $C = C(D)>0$ such that
\beq\label{e:lre1}
	\P{\left(\mathcal{N}(I_i(E))\geqslant 2\right)}\leqslant
	 \P{\left(\one_\mathcal{C}\mathcal{N}(I_i(E))\geqslant 2\right)}
	 +
	 C N^{-D},
\eeq
Consider any $\theta>0$. There exist $c>0$ and $C>0$ such that
\begin{align}
\P{\left(\one_\mathcal{C}\mathcal{N}(I_i(E))\geqslant 2\right)}
&\leqslant \sum_{j:\vert i-j\vert \leqslant N^{2\omega}} 
\P\left(
	\lambda_{i+1}-\lambda_i\leqslant N^{-2/3-\delta}\hat{i}^{-1/3}
\right)
\\
&\leqslant 
C\sum_{j:\vert i-j\vert \leqslant N^{2\omega}} 
\left(
	\P_{\mathrm{GOE}_N}\left(
		\lambda_{i+1}-\lambda_i\leqslant N^{-2/3-\delta}\hat{i}^{-1/3}
	\right)+N^{-c+\delta}
\right)
\\
&\leqslant CN^{2\omega-2\delta+C\theta}+CN^{-c+\delta+2\omega},
\end{align}
where the level repulsion estimates for the Gaussian ensembles used in the last inequality can be deduced from \cite{bourgade2014edge}*{Theorem 3.2}. In second line, we used gap universality with an observable $O$ defined as a smoothed indicator function on the scale $N^{-\delta}$ which obeys $\Vert O\Vert_\infty \leqslant CN^{\delta}$, and considered $O(\lambda_{i+1} - \lambda_i)$. Gap universality for this regime of $i$ follows from the proof of \cite{bourgade2018extreme}*{Theorem 1.6}, as mentioned previously. We now set $\delta_0=c/4$ and $\omega=C\theta=\delta/10$. For $\delta<\delta_0$, we have
\beq\label{e:lre2}
	\P{\left(\one_\mathcal{C}\mathcal{N}(I_i(E))\geqslant 2\right)}
	\leqslant
	CN^{-\delta}\left(
		N^{3\omega-\delta}+N^{-c + 2\delta + 2\omega}	
	\right)
	\leqslant CN^{-\delta-7\delta/10},
\eeq
and we obtain the final result by taking $\alpha=7\delta/10$ and combining \eqref{e:lre2} with \eqref{e:lre1}.
\end{proof}

\bp[Proof \Cref{l:lr2}]
The proof is the same as that of \Cref{p:lr}, after using rigidity to exclude the possibility that $\lambda_i$ is far from the spectrum $[-2, 2]$. Transposing the argument exactly loses a factor of $N^{\delta/5}$ in comparison to \Cref{prop:lrek2}, but this may be replaced by $N^{\delta/n}$ for any $n \in \N$, so the loss can be made arbitrarily small, and in particular less than the $\alpha$ from \Cref{prop:lrek2}.
\ep

As discussed in the introduction and following \cite{bourgade2013eigenvector}*{Definition 5.1}, by combining \Cref{prop:lrek2} with \cite{bourgade2013eigenvector}*{Theorem 1.2}, we obtain asymptotic Gaussianity of \emph{all} eigenvectors of generalized Wigner matrices. We state this result here for completeness.

\bec\label{c:asymptoticallynormal}
	Let $H$ be a real symmetric generalized Wigner matrix and fix $\q\in\S^{N-1}$. For any $m\in\N$ and $I\subset\unn{1}{N}$, we have
	\beq
		(\sqrt{N}\vert \scp{\q}{u_k}\vert)_{k\in I}\rightarrow (\vert \mathcal{N}_i)_{i=1}^m
	\eeq
	in the sense of convergence in moments, where $(\mathcal{N}_i)_{i=1}^m$ is a family of independent standard Gaussian random variables. In the complex Hermitian case, the $\mathcal{N}_i$ are replaced by complex-valued standard Gaussians.
\eec
\section{Preliminary estimates} \label{a:preliminary}

\begin{proof}[Proof of \Cref{l:good1}]
In the proof of \Cref{p:regulareigval}, in particular the lines following \eqref{e:conc1} and \eqref{e:conc2}, it was shown that if rigidity \eqref{e:rigidity}, delocalization, and  $\cup_{i,j \in \unn{1}{N} }\left \{ | h_{ij} | \le N^{-1/2 + \omega/12 } \right\}$ hold for $H$ with exponentially high probability in the sense of \eqref{e:exphigh}, then the local semicircle law \eqref{e:sclaw}, and rigidity and delocalization \eqref{e:rigidity}, hold uniformly for any rank-one perturbation $\Theta^{(a,b)}_w H$ also with exponentially high probability (with smaller constants $C_1, c_1 >0$, but the same spectral domain $\mathcal D_\omega$ and arbitrary control parameter $\omega >0$). It remains to show the isotropic local law \eqref{e:isolaw} also holds for the rank-one perturbations. For this, we observe that for the first term in a resolvent expansion of $\langle \q, G - \Theta^{(a,b)}_w G \q\rangle$, we have
\beq\label{e:c3}
\langle \q, G (H - \Theta^{(a,b)}_w H) G \q \rangle = w\langle \q, G \mathbf e_a \rangle h_{ab} \langle \mathbf e_b, G \q \rangle + w\langle \q, G \mathbf e_b \rangle h_{ba} \langle \mathbf e_a, G \q \rangle
\eeq
Here $\mathbf e_a$ and $\mathbf e_b$ are the standard basis vectors. Both terms can be bounded using polarization, $| h_{ab} | + | h_{ab} |\le N^{-1/2 + \omega/12}$, and \eqref{e:isolaw} for $H$, and we obtain that \eqref{e:c3} is bounded by $C N^{-1/2 + \omega/12}$ in absolute value for some $C(\omega) >0$.
Similar reasoning can be applied to the higher order terms in the resolvent expansion, and we conclude that \eqref{e:isolaw} holds for all $\Theta^{(a,b)}_w H$ after taking the expansion to a sufficiently high order, as in the proof of \Cref{p:regulareigval}.
Finally, we observe that $\cup_{i,j \in \unn{1}{N} }\left \{ | h_{ij} | \le N^{-1/2 + \omega/12 } \right\}$ holds with exponentially high probability, since the $h_{ij}$ were assumed to be subexponential. We therefore conclude using \Cref{l:goodset} for $H$ and a union bound.

\end{proof}

\begin{proof}[Proof of \Cref{l:good2}]
By \Cref{p:lr}, we have for all $\ell \in \unn{1}{N}$ that
\beq\label{e:unpt2}
\P \left( \mathcal{N}_{\bm{\lambda}}(I_{\delta_2}(\lambda_\ell ))\ge k
\right)
\leqslant N^{-D - 2 }
\eeq
when $N \ge N_0 (\delta_2)$ and $k\ge k_0(\delta_2)$, where $N_0$ and $k_0$ are independent of $\ell$. 
Applying a union bound to \eqref{e:unionthis1} and \eqref{e:unpt2} completes the proof.
\end{proof}

In preparation for the proof of \Cref{l:belowG}, we note the following lemma, which is essentially given in the proof of  \cite{bourgade2013eigenvector}*{Corollary A.2}. We omit the routine modifications necessary to derive the form given here.

\bel[\cite{bourgade2013eigenvector}*{Appendix A}]\label{l:extendedsc}
Fix  $\omega ,\delta_1,\varepsilon_1,\delta_2, \eps_2>  0$ and $\q\in\S^{N-1}$.  Let $H$ be a generalized Wigner matrix, and let $\B=\B(\q,\omega,\delta_1,\varepsilon_1,\delta_2, \eps_2)$ be the set from \Cref{d:goodset}. 
For all $z=E+\I \eta\in \D_{\eps_2/8}$ and $ 0 < y \le \eta$, 
\beq\label{e:extendedsc}
\sup_{a,b \in \unn{1}{N}} \sup_{w \in [0,1]} \one_{\B} | \langle  \q , \Theta^{(a,b)}_w G(E + \I y) \q \rangle| \le C \log N \frac{\eta}{y}
\eeq
for some constant $C(\eps_2)>0$.
\eel

\bp[Proof of \Cref{l:belowG}]
We work on the set $\B$ and omit this from the notation.
For $E \in I_{\delta_2}(\lambda_\ell) \cup I_{\delta_2}(\tilde \lambda_\ell)$, we have
\beq\label{e:mscbelow}
\Im \msc(E+\I\eta_\ell N^{2\eps_2} )\leqslant C \sqrt{ \left| |E| - 2\right |+ \eta_\ell N^{2\eps_2} }\leqslant CN^{\eps_2/2} \left(   \left(\frac{\hat \ell}{N}\right)^{1/3}  + N^{\eps_2/2} \sqrt{\eta_\ell} \right) \le C N^{\eps_2/2} \left( \frac{\hat \ell}{N}\right)^{1/3},
\eeq
where in the first inequality we used \Cref{l:msclower}, and in the second we used \eqref{e:regev}, $\sqrt{a + b} \le \sqrt{a} + \sqrt{b}$, and rigidity \eqref{e:rigidity} with parameter $\omega' = \eps_2$. Similarly, we obtain
\beq\label{e:mscbelow2}
\Im \msc(E+\I\eta_\ell )\le C N^{\eps_2/2} \left( \frac{\hat \ell}{N}\right)^{1/3}.
\eeq
When $N^{1 - (21/8)\eps_2} \ge \hat \ell$, then $\eta_{\ell} \ge N^{-1 + \eps_2/8}$ and $z = E + \I \eta_\ell \in \mathcal D_{\eps_2/8}$, so it is permissible to apply the isotropic local law \eqref{e:isolaw} to $\Theta^{(a,b)}_w H$ with $\omega' = \eps_2/8$, which gives
\beq\label{e:appiso}
\left\vert
	\scp{\q}{\Theta^{(a,b)}_w G^s(z)\q}
	-
	\msc(z)
\right\vert
\le C N^{\eps_2} \left(\frac{\hat \ell}{N} \right)^{1/3}.
\eeq
Together, \eqref{e:appiso} and \eqref{e:mscbelow2} imply \eqref{e:belowG}.

In the case that $\hat \ell \ge N^{1 - (21/8)\eps_2}$, we can again use the isotropic local law \eqref{e:isolaw} with $z =E + \I \eta_{\ell} N^{2\eps_2} $ and $\omega' = \eps_2/8$, \eqref{e:extendedsc} with $y = \eta_\ell$, and \eqref{e:mscbelow}, to obtain \eqref{e:belowG}.

Given \eqref{e:belowG} and \eqref{e:extendedsc}, \eqref{e:extendedderiv} follows by a standard argument using the Green's function differentiation formula $\partial_{ab} G_{ij}  = -G_{ia} G_{bj}$ after expanding the inner product $ \langle \q , G(E + \I \eta_{\ell} ) \q \rangle$ and computing the imaginary part. 
Full details can be found in the proof of \cite{bourgade2013eigenvector}*{Corollary A.2}.
\ep

\begin{proof}[Proof of \Cref{l:eigsqsum}]
When $\mathcal B$ holds and $| p  -\ell | > N^{2\omega}$, we have $| \lambda_p - \lambda_\ell| > \frac{1}{2} | \gamma_p - \gamma_\ell|$ for $N \ge N_0 (\omega)$. It then suffices to show
\beq\label{e:eigsqsum2}
\one_{\B} \sum_{p: \vert p-\ell\vert>N^{2\omega}}
\frac{1}{(\gamma_p-\gamma_\ell)^2}
\leqslant
C N^{4/3 + 2\omega} {\hat \ell}^{2/3}.
\eeq
This is a straightforward computation and can be accomplished using, for example, the technique demonstrated in \eqref{e:citeatend} and the following lines.
\end{proof}

\bibliography{deloc}
\bibliographystyle{abbrv}

\end{document}